

\input amstex



\catcode`\X=12\catcode`\@=11

\def\n@wcount{\alloc@0\count\countdef\insc@unt}
\def\n@wwrite{\alloc@7\write\chardef\sixt@@n}
\def\n@wread{\alloc@6\read\chardef\sixt@@n}
\def\r@s@t{\relax}\def\v@idline{\par}\def\@mputate#1/{#1}
\def\l@c@l#1X{\firstpart.#1}\def\gl@b@l#1X{#1}\def\t@d@l#1X{{}}

\def\crossrefs#1{\ifx\all#1\let\tr@ce=\all\else\def\tr@ce{#1,}\fi
   \n@wwrite\cit@tionsout\openout\cit@tionsout=\jobname.cit 
   \write\cit@tionsout{\tr@ce}\expandafter\setfl@gs\tr@ce,}
\def\setfl@gs#1,{\def\@{#1}\ifx\@\empty\let\next=\relax
   \else\let\next=\setfl@gs\expandafter\xdef
   \csname#1tr@cetrue\endcsname{}\fi\next}
\def\m@ketag#1#2{\expandafter\n@wcount\csname#2tagno\endcsname
     \csname#2tagno\endcsname=0\let\tail=\all\xdef\all{\tail#2,}
   \ifx#1\l@c@l\let\tail=\r@s@t\xdef\r@s@t{\csname#2tagno\endcsname=0\tail}\fi
   \expandafter\gdef\csname#2cite\endcsname##1{\expandafter
     \ifx\csname#2tag##1\endcsname\relax?\else\csname#2tag##1\endcsname\fi
     \expandafter\ifx\csname#2tr@cetrue\endcsname\relax\else
     \write\cit@tionsout{#2tag ##1 cited on page \folio.}\fi}
   \expandafter\gdef\csname#2page\endcsname##1{\expandafter
     \ifx\csname#2page##1\endcsname\relax?\else\csname#2page##1\endcsname\fi
     \expandafter\ifx\csname#2tr@cetrue\endcsname\relax\else
     \write\cit@tionsout{#2tag ##1 cited on page \folio.}\fi}
   \expandafter\gdef\csname#2tag\endcsname##1{\expandafter
      \ifx\csname#2check##1\endcsname\relax
      \expandafter\xdef\csname#2check##1\endcsname{}%
      \else\immediate\write16{Warning: #2tag ##1 used more than once.}\fi
      \multit@g{#1}{#2}##1/X%
      \write\t@gsout{#2tag ##1 assigned number \csname#2tag##1\endcsname\space
      on page \number\count0.}%
   \csname#2tag##1\endcsname}}
\def\multit@g#1#2#3/#4X{\def\t@mp{#4}\ifx\t@mp\empty%
      \global\advance\csname#2tagno\endcsname by 1 
      \expandafter\xdef\csname#2tag#3\endcsname
      {#1\number\csname#2tagno\endcsnameX}%
   \else\expandafter\ifx\csname#2last#3\endcsname\relax
      \expandafter\n@wcount\csname#2last#3\endcsname
      \global\advance\csname#2tagno\endcsname by 1 
      \expandafter\xdef\csname#2tag#3\endcsname
      {#1\number\csname#2tagno\endcsnameX}
      \write\t@gsout{#2tag #3 assigned number \csname#2tag#3\endcsname\space
      on page \number\count0.}\fi
   \global\advance\csname#2last#3\endcsname by 1
   \def\t@mp{\expandafter\xdef\csname#2tag#3/}%
   \expandafter\t@mp\@mputate#4\endcsname
   {\csname#2tag#3\endcsname\lastpart{\csname#2last#3\endcsname}}\fi}
\def\t@gs#1{\def\all{}\m@ketag#1e\m@ketag#1s\m@ketag\t@d@l p
   \m@ketag\gl@b@l r \n@wread\t@gsin
   \openin\t@gsin=\jobname.tgs \re@der \closein\t@gsin
   \n@wwrite\t@gsout\openout\t@gsout=\jobname.tgs }
\outer\def\localtags{\t@gs\l@c@l}
\outer\def\globaltags{\t@gs\gl@b@l}
\outer\def\newlocaltag#1{\m@ketag\l@c@l{#1}}
\outer\def\newglobaltag#1{\m@ketag\gl@b@l{#1}}

\newif\ifpr@ 
\def\m@kecs #1tag #2 assigned number #3 on page #4.%
   {\expandafter\gdef\csname#1tag#2\endcsname{#3}
   \expandafter\gdef\csname#1page#2\endcsname{#4}
   \ifpr@\expandafter\xdef\csname#1check#2\endcsname{}\fi}
\def\re@der{\ifeof\t@gsin\let\next=\relax\else
   \read\t@gsin to\t@gline\ifx\t@gline\v@idline\else
   \expandafter\m@kecs \t@gline\fi\let \next=\re@der\fi\next}
\def\pretags#1{\pr@true\pret@gs#1,,}
\def\pret@gs#1,{\def\@{#1}\ifx\@\empty\let\n@xtfile=\relax
   \else\let\n@xtfile=\pret@gs \openin\t@gsin=#1.tgs \message{#1} \re@der 
   \closein\t@gsin\fi \n@xtfile}

\newcount\sectno\sectno=0\newcount\subsectno\subsectno=0
\newif\ifultr@local \def\ultralocal{\ultr@localtrue}
\def\firstpart{\number\sectno}
\def\lastpart#1{\ifcase#1 \or a\or b\or c\or d\or e\or f\or g\or h\or 
   i\or k\or l\or m\or n\or o\or p\or q\or r\or s\or t\or u\or v\or w\or 
   x\or y\or z \fi}

\def\resetall{\global\advance\sectno by 1\subsectno=0
   \gdef\firstpart{\number\sectno}\r@s@t}
\def\resetsub{\global\advance\subsectno by 1
   \gdef\firstpart{\number\sectno.\number\subsectno}\r@s@t}
\def\newsection#1\par{\resetall\vskip0pt plus.3\vsize\penalty-250
   \vskip0pt plus-.3\vsize\bigskip\bigskip
   \message{#1}\leftline{\bf#1}\nobreak\bigskip}
\def\subsection#1\par{\ifultr@local\resetsub\fi
   \vskip0pt plus.2\vsize\penalty-250\vskip0pt plus-.2\vsize
   \bigskip\smallskip\message{#1}\leftline{\bf#1}\nobreak\medskip}

\def\t@gsoff#1,{\def\@{#1}\ifx\@\empty\let\next=\relax\else\let\next=\t@gsoff
   \def\@@{p}\ifx\@\@@\else
   \expandafter\gdef\csname#1cite\endcsname##1{\zeigen{##1}}
   \expandafter\gdef\csname#1page\endcsname##1{?}
   \expandafter\gdef\csname#1tag\endcsname##1{\zeigen{##1}}\fi\fi\next}
\def\verbatimtags{\ifx\all\relax\else\expandafter\t@gsoff\all,\fi}
\def\zeigen#1{\hbox{$\langle$}#1\hbox{$\rangle$}}

\def\(#1){\edef\dot@g{\ifmmode\ifinner(\hbox{\noexpand\etag{#1}})
   \else\noexpand\eqno(\hbox{\noexpand\etag{#1}})\fi
   \else(\noexpand\ecite{#1})\fi}\dot@g}

\newif\ifbr@ck
\def\eat#1{}
\def\[#1]{\br@cktrue[\br@cket#1'X]}
\def\br@cket#1'#2X{\def\temp{#2}\ifx\temp\empty\let\next\eat
   \else\let\next\br@cket\fi
   \ifbr@ck\br@ckfalse\br@ck@t#1,X\else\br@cktrue#1\fi\next#2X}
\def\br@ck@t#1,#2X{\def\temp{#2}\ifx\temp\empty\let\neext\eat
   \else\let\neext\br@ck@t\def\temp{,}\fi
   \def\teemp{#1}\ifx\teemp\empty\else\rcite{#1}\fi\temp\neext#2X}
\def\resetbr@cket{\gdef\[##1]{[\rtag{##1}]}}
\def\references{\resetbr@cket\newsection References\par}

\newtoks\symb@ls\newtoks\s@mb@ls\newtoks\p@gelist\n@wcount\ftn@mber
    \ftn@mber=1\newif\ifftn@mbers\ftn@mbersfalse\newif\ifbyp@ge\byp@gefalse
\def\defm@rk{\ifftn@mbers\n@mberm@rk\else\symb@lm@rk\fi}
\def\n@mberm@rk{\xdef\m@rk{{\the\ftn@mber}}%
    \global\advance\ftn@mber by 1 }
\def\rot@te#1{\let\temp=#1\global#1=\expandafter\r@t@te\the\temp,X}
\def\r@t@te#1,#2X{{#2#1}\xdef\m@rk{{#1}}}
\def\b@@st#1{{$^{#1}$}}\def\str@p#1{#1}
\def\symb@lm@rk{\ifbyp@ge\rot@te\p@gelist\ifnum\expandafter\str@p\m@rk=1 
    \s@mb@ls=\symb@ls\fi\write\f@nsout{\number\count0}\fi \rot@te\s@mb@ls}
\def\byp@ge{\byp@getrue\n@wwrite\f@nsin\openin\f@nsin=\jobname.fns 
    \n@wcount\currentp@ge\currentp@ge=0\p@gelist={0}
    \re@dfns\closein\f@nsin\rot@te\p@gelist
    \n@wread\f@nsout\openout\f@nsout=\jobname.fns }
\def\m@kelist#1X#2{{#1,#2}}
\def\re@dfns{\ifeof\f@nsin\let\next=\relax\else\read\f@nsin to \f@nline
    \ifx\f@nline\v@idline\else\let\t@mplist=\p@gelist
    \ifnum\currentp@ge=\f@nline
    \global\p@gelist=\expandafter\m@kelist\the\t@mplistX0
    \else\currentp@ge=\f@nline
    \global\p@gelist=\expandafter\m@kelist\the\t@mplistX1\fi\fi
    \let\next=\re@dfns\fi\next}
\def\symbols#1{\symb@ls={#1}\s@mb@ls=\symb@ls} 
\def\bigsymbol{\textstyle}
\symbols{\bigsymbol\ast,\dagger,\ddagger,\sharp,\flat,\natural,\star}
\def\ftnumbers{\ftn@mberstrue} \def\ftsymbols{\ftn@mbersfalse}
\def\paginal{\byp@ge} \def\resetftnumbers{\ftn@mber=1}
\def\ftnote#1{\defm@rk\expandafter\expandafter\expandafter\footnote
    \expandafter\b@@st\m@rk{#1}}

\long\def\jump#1\endjump{}
\def\ssum{\mathop{\lower .1em\hbox{$\textstyle\Sigma$}}\nolimits}

\def\qed{\nobreak\kern 1em \vrule height .5em width .5em depth 0em}
\def\newneq{\hbox{\rlap{\hbox to 1\wd9{\hss$=$\hss}}\raise .1em 
   \hbox to 1\wd9{\hss$\scriptscriptstyle/$\hss}}}
\def\subsetne{\setbox9 = \hbox{$\subset$}\mathrel{\hbox{\rlap
   {\lower .4em \newneq}\raise .13em \hbox{$\subset$}}}}
\def\supsetne{\setbox9 = \hbox{$\subset$}\mathrel{\hbox{\rlap
   {\lower .4em \newneq}\raise .13em \hbox{$\supset$}}}}

\def\vbar{\mathchoice{\vrule height6.3ptdepth-.5ptwidth.8pt\kern-.8pt}
   {\vrule height6.3ptdepth-.5ptwidth.8pt\kern-.8pt}
   {\vrule height4.1ptdepth-.35ptwidth.6pt\kern-.6pt}
   {\vrule height3.1ptdepth-.25ptwidth.5pt\kern-.5pt}}
\def\f@dge{\mathchoice{}{}{\mkern.5mu}{\mkern.8mu}}
\def\b@c#1#2{{\rm \mkern#2mu\vbar\mkern-#2mu#1}}
\def\b@b#1{{\rm I\mkern-3.5mu #1}}
\def\b@a#1#2{{\rm #1\mkern-#2mu\f@dge #1}}
\def\bb#1{{\count4=`#1 \advance\count4by-64 \ifcase\count4\or\b@a A{11.5}\or
   \b@b B\or\b@c C{5}\or\b@b D\or\b@b E\or\b@b F \or\b@c G{5}\or\b@b H\or
   \b@b I\or\b@c J{3}\or\b@b K\or\b@b L \or\b@b M\or\b@b N\or\b@c O{5} \or
   \b@b P\or\b@c Q{5}\or\b@b R\or\b@a S{8}\or\b@a T{10.5}\or\b@c U{5}\or
   \b@a V{12}\or\b@a W{16.5}\or\b@a X{11}\or\b@a Y{11.7}\or\b@a Z{7.5}\fi}}

\catcode`\X=11 \catcode`\@=12

\def\sciteu{\sciteerror{undefined}}
\def\sciteuphantom{\complainaboutcitation{undefined}}

\def\sciteerror#1#2{{\mathortextbf{\scite{#2}}}\complainaboutcitation{#1}{#2}}
\def\mathortextbf#1{\hbox{\bf #1}}
\def\complainaboutcitation#1#2{%
\vadjust{\line{\llap{---$\!\!>$ }\qquad scite$\{$#2$\}$ #1\hfil}}}

\sectno=-1   
\localtags
\NoBlackBoxes
\define\mr{\medskip\roster}
\define\sn{\smallskip\noindent}
\define\mn{\medskip\noindent}
\define\bn{\bigskip\noindent}
\define\ub{\underbar}
\define\wilog{\text{without loss of generality}}
\define\ermn{\endroster\medskip\noindent}
\define\dbca{\dsize\bigcap}
\define\dbcu{\dsize\bigcup}
\define \nl{\newline}
\documentstyle {amsppt}
\topmatter
\title{On Arhangelskii's Problem \\
 Sh668} \endtitle
\author {Saharon Shelah \thanks {\null\newline I would like to thank 
Alice Leonhardt for the beautiful typing. \null\newline
 Done Sept. (and Dec.) 1997 \null\newline
 Latest Revision - 98/Oct/21} \endthanks} \endauthor 
\affil{Institute of Mathematics\\
 The Hebrew University\\
 Jerusalem, Israel
 \medskip
 Rutgers University\\
 Mathematics Department\\
 New Brunswick, NJ  USA} \endaffil
\mn
\abstract We prove the consistency (modulo supercompact) of a negative
answer to Arhangelskii's problem (some Hausdorff compact space cannot be
partitioned to two sets not containing a closed copy of Cantor discontinuum).
In this model we have CH.  Without CH we get consistency results using a pcf
assumption, close relatives of which are necessary for such results. 
\endabstract
\endtopmatter
\document  
\def\renewcommand{\newcommand}	       
\edef\cite{\the\catcode`@}%
\catcode`@ = 11
\let\@oldatcatcode = \cite
\chardef\@letter = 11
\chardef\@other = 12
%
%
%
%
\def\@innerdef#1#2{\edef#1{\expandafter\noexpand\csname #2\endcsname}}%
%
%
\@innerdef\@innernewcount{newcount}%
\@innerdef\@innernewdimen{newdimen}%
\@innerdef\@innernewif{newif}%
\@innerdef\@innernewwrite{newwrite}%
%
%
%
\def\@gobble#1{}%
%
%
%
\ifx\inputlineno\@undefined
   \let\@linenumber = \empty 
\else
   \def\@linenumber{\the\inputlineno:\space}%
\fi
%
%
%
\def\@futurenonspacelet#1{\def\cs{#1}%
   \afterassignment\@stepone\let\@nexttoken=
}%
\begingroup 
\def\\{\global\let\@stoken= }%
\\ 
\endgroup
\def\@stepone{\expandafter\futurelet\cs\@steptwo}%
\def\@steptwo{\expandafter\ifx\cs\@stoken\let\@@next=\@stepthree
   \else\let\@@next=\@nexttoken\fi \@@next}%
\def\@stepthree{\afterassignment\@stepone\let\@@next= }%
%
%
%
\def\@getoptionalarg#1{%
   \let\@optionaltemp = #1%
   \let\@optionalnext = \relax
   \@futurenonspacelet\@optionalnext\@bracketcheck
}%
%
%
\def\@bracketcheck{%
   \ifx [\@optionalnext
      \expandafter\@@getoptionalarg
   \else
      \let\@optionalarg = \empty
      \expandafter\@optionaltemp
   \fi
}%
\def\@@getoptionalarg[#1]{%
   \def\@optionalarg{#1}%
   \@optionaltemp
}%
%
%
%
\def\@nnil{\@nil}%
\def\@fornoop#1\@@#2#3{}%
\def\@for#1:=#2\do#3{%
   \edef\@fortmp{#2}%
   \ifx\@fortmp\empty \else
      \expandafter\@forloop#2,\@nil,\@nil\@@#1{#3}%
   \fi
}%
\def\@forloop#1,#2,#3\@@#4#5{\def#4{#1}\ifx #4\@nnil \else
       #5\def#4{#2}\ifx #4\@nnil \else#5\@iforloop #3\@@#4{#5}\fi\fi
}%
\def\@iforloop#1,#2\@@#3#4{\def#3{#1}\ifx #3\@nnil
       \let\@nextwhile=\@fornoop \else
      #4\relax\let\@nextwhile=\@iforloop\fi\@nextwhile#2\@@#3{#4}%
}%
%
%
%
\@innernewif\if@fileexists
\def\@testfileexistence{\@getoptionalarg\@finishtestfileexistence}%
\def\@finishtestfileexistence#1{%
   \begingroup
      \def\extension{#1}%
      \immediate\openin0 =
         \ifx\@optionalarg\empty\jobname\else\@optionalarg\fi
         \ifx\extension\empty \else .#1\fi
         \space
      \ifeof 0
         \global\@fileexistsfalse
      \else
         \global\@fileexiststrue
      \fi
      \immediate\closein0
   \endgroup
}%
%
%
%
%
\def\bibliographystyle#1{%
   \@readauxfile
   \@writeaux{\string\bibstyle{#1}}%
}%
\let\bibstyle = \@gobble
%
%
\let\bblfilebasename = \jobname
\def\bibliography#1{%
   \@readauxfile
   \@writeaux{\string\bibdata{#1}}%
   \@testfileexistence[\bblfilebasename]{bbl}%
   \if@fileexists
      \nobreak
      \@readbblfile
   \fi
}%
\let\bibdata = \@gobble
%
%
\def\nocite#1{%
   \@readauxfile
   \@writeaux{\string\citation{#1}}%
}%
\@innernewif\if@notfirstcitation
%
%
\def\cite{\@getoptionalarg\@cite}%
%
%
\def\@cite#1{%
   \let\@citenotetext = \@optionalarg
   \printcitestart
   \nocite{#1}%
   \@notfirstcitationfalse
   \@for \@citation :=#1\do
   {%
      \expandafter\@onecitation\@citation\@@
   }%
   \ifx\empty\@citenotetext\else
      \printcitenote{\@citenotetext}%
   \fi
   \printcitefinish
}%
\def\@onecitation#1\@@{%
   \if@notfirstcitation
      \printbetweencitations
   \fi
   \expandafter \ifx \csname\@citelabel{#1}\endcsname \relax
      \if@citewarning
         \message{\@linenumber Undefined citation `#1'.}%
      \fi
      \expandafter\gdef\csname\@citelabel{#1}\endcsname{%
\strut
\vadjust{\vskip-\dp\strutbox
\vbox to 0pt{\vss\parindent0cm \leftskip=\hsize 
\advance\leftskip3mm
\advance\hsize 4cm\strut\openup-4pt 
\rightskip 0cm plus 1cm minus 0.5cm ?  #1 ?\strut}}
         {\tt
            \escapechar = -1
            \nobreak\hskip0pt
            \expandafter\string\csname#1\endcsname
            \nobreak\hskip0pt
         }%
      }%
   \fi
   \csname\@citelabel{#1}\endcsname
   \@notfirstcitationtrue
}%
%
%
\def\@citelabel#1{b@#1}%
%
%
\def\@citedef#1#2{\expandafter\gdef\csname\@citelabel{#1}\endcsname{#2}}%
%
%
%
\def\@readbblfile{%
   \ifx\@itemnum\@undefined
      \@innernewcount\@itemnum
   \fi
   \begingroup
      \def\begin##1##2{%
         \setbox0 = \hbox{\biblabelcontents{##2}}%
         \biblabelwidth = \wd0
      }%
      \def\end##1{}
      %
      %
      \@itemnum = 0
      \def\bibitem{\@getoptionalarg\@bibitem}%
      \def\@bibitem{%
         \ifx\@optionalarg\empty
            \expandafter\@numberedbibitem
         \else
            \expandafter\@alphabibitem
         \fi
      }%
      \def\@alphabibitem##1{%
         \expandafter \xdef\csname\@citelabel{##1}\endcsname {\@optionalarg}%
         \ifx\biblabelprecontents\@undefined
            \let\biblabelprecontents = \relax
         \fi
         \ifx\biblabelpostcontents\@undefined
            \let\biblabelpostcontents = \hss
         \fi
         \@finishbibitem{##1}%
      }%
      \def\@numberedbibitem##1{%
         \advance\@itemnum by 1
         \expandafter \xdef\csname\@citelabel{##1}\endcsname{\number\@itemnum}%
         \ifx\biblabelprecontents\@undefined
            \let\biblabelprecontents = \hss
         \fi
         \ifx\biblabelpostcontents\@undefined
            \let\biblabelpostcontents = \relax
         \fi
         \@finishbibitem{##1}%
      }%
      \def\@finishbibitem##1{%
         \biblabelprint{\csname\@citelabel{##1}\endcsname}%
         \@writeaux{\string\@citedef{##1}{\csname\@citelabel{##1}\endcsname}}%
         \ignorespaces
      }%
      %
      %
      \let\em = \bblem
      \let\newblock = \bblnewblock
      \let\sc = \bblsc
      \frenchspacing
      \clubpenalty = 4000 \widowpenalty = 4000
      \tolerance = 10000 \hfuzz = .5pt
      \everypar = {\hangindent = \biblabelwidth
                      \advance\hangindent by \biblabelextraspace}%
      \bblrm
      \parskip = 1.5ex plus .5ex minus .5ex
      \biblabelextraspace = .5em
      \bblhook
      \input \bblfilebasename.bbl
   \endgroup
}%
%
%
\@innernewdimen\biblabelwidth
\@innernewdimen\biblabelextraspace
%
%
%
\def\biblabelprint#1{%
   \noindent
   \hbox to \biblabelwidth{%
      \biblabelprecontents
      \biblabelcontents{#1}%
      \biblabelpostcontents
   }%
   \kern\biblabelextraspace
}%
%
%
%
\def\biblabelcontents#1{{\bblrm [#1]}}%
%
%
\def\bblrm{\rm}%
%
%
\def\bblem{\it}%
%
%
\def\bblsc{\ifx\@scfont\@undefined
              \font\@scfont = cmcsc10
           \fi
           \@scfont
}%
%
%
\def\bblnewblock{\hskip .11em plus .33em minus .07em }%
%
%
\let\bblhook = \empty
%
%
%
\def\printcitestart{[}
\def\printcitefinish{]}
\def\printbetweencitations{, }
\def\printcitenote#1{, #1}
%
%
%
\let\citation = \@gobble
%
%
%
\@innernewcount\@numparams
%
%
\def\newcommand#1{%
   \def\@commandname{#1}%
   \@getoptionalarg\@continuenewcommand
}%
%
%
\def\@continuenewcommand{%
   \@numparams = \ifx\@optionalarg\empty 0\else\@optionalarg \fi \relax
   \@newcommand
}%
%
%
\def\@newcommand#1{%
   \def\@startdef{\expandafter\edef\@commandname}%
   \ifnum\@numparams=0
      \let\@paramdef = \empty
   \else
      \ifnum\@numparams>9
         \errmessage{\the\@numparams\space is too many parameters}%
      \else
         \ifnum\@numparams<0
            \errmessage{\the\@numparams\space is too few parameters}%
         \else
            \edef\@paramdef{%
               \ifcase\@numparams
                  \empty  No arguments.
               \or ####1%
               \or ####1####2%
               \or ####1####2####3%
               \or ####1####2####3####4%
               \or ####1####2####3####4####5%
               \or ####1####2####3####4####5####6%
               \or ####1####2####3####4####5####6####7%
               \or ####1####2####3####4####5####6####7####8%
               \or ####1####2####3####4####5####6####7####8####9%
               \fi
            }%
         \fi
      \fi
   \fi
   \expandafter\@startdef\@paramdef{#1}%
}%
%
%
%
%
\def\@readauxfile{%
   \if@auxfiledone \else 
      \global\@auxfiledonetrue
      \@testfileexistence{aux}%
      \if@fileexists
         \begingroup
            \endlinechar = -1
            \catcode`@ = 11
            \input \jobname.aux
         \endgroup
      \else
         \message{\@undefinedmessage}%
         \global\@citewarningfalse
      \fi
      \immediate\openout\@auxfile = \jobname.aux
   \fi
}%
%
%
\newif\if@auxfiledone
\ifx\noauxfile\@undefined \else \@auxfiledonetrue\fi
%
%
%
%
\@innernewwrite\@auxfile
\def\@writeaux#1{\ifx\noauxfile\@undefined \write\@auxfile{#1}\fi}%
%
%
%
\ifx\@undefinedmessage\@undefined
   \def\@undefinedmessage{No .aux file; I won't give you warnings about
                          undefined citations.}%
\fi
%
%
\@innernewif\if@citewarning
\ifx\noauxfile\@undefined \@citewarningtrue\fi
%
%
%
\catcode`@ = \@oldatcatcode

\def\widestnumber#1#2{}
\def\rm{\fam0 \tenrm}
\def\fakesubhead#1\endsubhead{\bigskip\noindent{\bf#1}\par}
\newpage

\head {Anotated Content} \endhead  \resetall 
\bn
\S1 $\quad$ General spaces: consistency from strong assumptions
\mr
\item "{{}}"  [We define $X^* \rightarrow (Y^*)^1_\theta$ for topological
spaces $X^*,Y^*$.  Then starting with a Hausdorff space $Y^*$ with $\theta$
points such that any set of $< \sigma$ members is discrete and $\kappa = 
\kappa^{< \kappa} \in (\theta,\lambda)$ and appropriate ${\Cal A} 
\subseteq [\lambda]^\theta$ such that any two members has intersection 
$< \sigma$, we force
appropriate $X^*$.  We then show that the assumption holds under appropriate
pcf assumption and finish with some improvements.]
\endroster
\bn
\S2 $\quad$ Consistency from supercompact, with clopen basis
\mr
\item "{{}}"  [We deal here with the set theoretic assumption.  We show that
the assumptions can be gotten from supercompact for the case we agree to
have CH, relying on earlier 
consistency results.]
\endroster
\bn
\S3 $\quad$ Equi-consistency
\mr
\item "{{}}"  [We show that some versions of the topological question and 
suitable combinatorial questions are equi-consistents.  See
\cite{Sh:108}, \cite{HJSh:249}, \cite{Sh:460}, \cite{Sh:F276}.  
\ub{Saharon} We then 
indicate the changes
needed for the not necessarily closed subspace case colouring by more colours
and other spaces.  For discussion see \cite{Sh:666},\S1.]
\endroster
\bn
\S4 $\quad$  Helping equi-consistency
\newpage

\head {\S1 General spaces: consistency from strong assumptions} \endhead  \resetall 
\bigskip

In our main theorem, \scite{t.2}, we give set theoretic sufficient conditions
for being able to force counterexamples to Arhangelskii's problem, possibly
replacing the cantor discontinuum by any other space.  It has a version for
spaces with clopen basis.  Then (in claim \scite{t.3}) we connect this to
pcf theory: after easy forcing the assumptions of Theorem \scite{t.2} can be
proved, if we start with a suitable (strong) pcf assumption (whose status is
not known).  Then in claim \scite{t.4} we deal with variants of the theorem,
weakening the topological and/or set theoretic assumptions.  Further variants
are discussed in the end ($T_3$ spaces without clopen basis and variants of
\scite{t.3}).
\bigskip

\definition{\stag{t.1} Definition}  Let $n \in [1,\omega)$ (though we
concentrate on $n=1$) ${{}}$. \nl
1) We say $X^* \rightarrow (Y^*)^n_\theta$, \ub{if} $X^*,Y^*$ are topological
spaces and for every $h:[X^*]^n \rightarrow \theta$ there is a closed
subspace $Y$ of $X^*$, homeomorphic to $Y^*$ such that $h \restriction [Y]^n$
is constant (if $n=1$ we may write $h:X^* \rightarrow \theta,h \restriction
Y$). \nl
2) If we omit the ``closed", we shall write $\rightarrow_w$ instead of
$\rightarrow$.  We write $(Y^*)^n_{< \theta}$ meaning: for every $h:[X^*]^n
\rightarrow \gamma < \theta$.  We use $\nrightarrow,\nrightarrow_w$ for the
negations.
\enddefinition
\bigskip

\proclaim{\stag{t.2} Theorem}  Assume
\mr
\item "{$(A)$}"  $\quad (i) \quad \lambda > \kappa > \theta > 
\sigma \ge \aleph_0$ and $\kappa = \kappa^{< \kappa}$ 
\sn
\item "{${{}}$}" $\quad (ii) \quad (\forall \alpha < \kappa)(|\alpha|^\sigma <
\kappa)$ and $\kappa > \theta^* \ge \theta$
\sn
\item "{$(B)_1$}"  ${\Cal A} \subseteq [\lambda]^\theta$ and \nl
$A_1 \ne A_2 \in {\Cal A} \Rightarrow |A_1 \cap A_2| < \sigma$
\sn
\item "{$(B)_2$}"  ${\Cal A}$ is $(< \kappa)$-free which means:
if ${\Cal A}' \subseteq {\Cal A},|{\Cal A}'| < \kappa$
\ub{then} for some list \nl
$\{A_\varepsilon:\varepsilon < \zeta\}$ of ${\Cal A}'$, for each
$\varepsilon < \zeta$  we have \nl
$|A_\varepsilon \cap \dbcu_{\xi < \varepsilon} A_\xi| < \sigma$
\sn
\item "{$(C)$}"  if $F:\lambda \rightarrow [\lambda]^{\le \kappa}$, \ub{then}
some $A \in {\Cal A}$ (or just some $A$ such that $(\exists A')(A \subseteq
A' \and |A| = \theta \and A' \in {\Cal A})$ is $F$-free which means
\sn
{\roster
\itemitem{ $(*)$ }  for $\alpha \ne \beta$ from $A$ we have $\alpha \notin
F(\beta)$
\endroster}
\item "{$(D)$}"  $Y^*$ is a Hausdorff space with set of points $\theta$ and
a basis ${\Cal B} = \{b_i:i < \theta^*\}$
\sn
\item "{$(E)$}"  if $Y$ is a subset of $Y^*$ with $< \sigma$ points,
\ub{then} $Y$ is a discrete subset \nl
(if $\sigma = \aleph_0$ this follows from Hausdorff), i.e. there is a
sequence of open (for $Y^*$) pairwise disjoint sets $\langle {\Cal U}_y:
y \in Y \rangle$, such that $y \in {\Cal U}_y$.
\ermn
\ub{Then} \nl
1) for some $\kappa$-complete $\kappa^+$-c.c. forcing notion $P$, in $V^P$
there is $X^*$ such that:
\mr
\item "{$(a)$}"  $X^*$ is a Hausdorff topological space with $\lambda$
points and basis of size $|{\Cal A}| + \theta^*$
\sn
\item "{$(b)$}"  $X^* \rightarrow (Y^*)^1_{< \text{ cf}(\theta)}$ (that is if
$X^* = \dbcu_{i < i(*)} X_i$ where $i(*) < \text{ cf}(\theta)$ \ub{then} some 
closed subspace $Y$ of $X^*$ homeomorphic to $Y^*$ is included in some single
$X_i$ (i.e. $(\exists i)(Y \subseteq X_i)$). 
\ermn
2) If in addition $Y^*$ has a clopen basis ${\Cal B}$ of cardinality
$\le \theta^*$ such that the union of $< \sigma$ members of ${\Cal B}$ 
is clopen, \ub{then} we can require that $X^*$ has a clopen basis.
\endproclaim
\bigskip

\remark{\stag{t.2A} Remark}  We may define the conditions historically (see
\cite{ShSt:258}, \cite{RoSh:599}, so put only the required conditions).
Then we can allow $\theta^* = \kappa$, but see \scite{t.4}.
\endremark
\bigskip

\demo{Proof}  We write the proof for part (1) and indicate the changes for
part (2).  Without loss of generality
\mr
\item "{$\bigotimes_1$}"  $(\forall \alpha < \beta < \lambda)
(\forall B \in [\lambda]^{< \lambda})
(\exists^{\kappa^+}A \in {\Cal A})[\{\alpha,\beta\} \subseteq A \and A \cap
B \subseteq \{\alpha,\beta\}]$. \nl
\mn
[Why?  As we can use $\bigl\{ \{2 \alpha:\alpha \in A\}:A \in {\Cal A}
\bigr\}$, \wilog \, $\bigcup\{A:A \in {\Cal A}\} = \{2\alpha:\alpha <
\lambda\}$ and choose $A_{\alpha,\beta,\gamma} \in [\lambda]^\theta$ for
$\alpha < \beta < \gamma < \lambda$ such that $\{\alpha,\beta\} \subseteq
A_{\alpha,\beta,\gamma}$ and $\langle A_{\alpha,\beta,\gamma} \backslash
\{\alpha,\beta\}:\alpha < \beta < \gamma < \lambda \rangle$ are pairwise
disjoint subsets of $\{2 \alpha +1:\alpha < \lambda\}$, each of cardinality
$\theta$ and replace ${\Cal A}$ by ${\Cal A}^* =: {\Cal A} \cup \{
A_{\alpha,\beta,\gamma}:\alpha < \beta < \gamma < \lambda\}$.
Now clause (A), (D), (E) are not affected.
Clearly clause $(B)_1$ holds (i.e. ${\Cal A}^* \subseteq [\lambda]^\theta$
and $A \ne B \in {\Cal A}^* \Rightarrow |A \cap B| < \sigma$).  Also clause
(C) is inherited by any extension of the original ${\Cal A}$.  Lastly for
clause $(B)_2$, if ${\Cal A}' \subseteq {\Cal A}^*,|{\Cal A}'| < \kappa$,
let $\langle A_\zeta:\zeta < \zeta^* \rangle$ be a list of ${\Cal A}' \cap
{\Cal A}$ as guaranteed by $(B)_2$ and let $\langle A_\zeta:\zeta \in
[\zeta^*,\zeta^* + |{\Cal A}' \backslash {\Cal A}|) \rangle$ list with no
repetitions ${\Cal A}' \backslash {\Cal A}$, now check.]
\sn
\item "{$\bigotimes_2$}"   ${\Cal B}$ is a basis of $Y^*$ of cardinality
$\theta^*$, and for part (2), ${\Cal B}$ is as there. \nl
\mn
[Why?  Straight.]
\ermn
Let ${\Cal A} = \{A_\zeta:\zeta < \lambda^*\}$ and ${\Cal B} = \{b_i:i <
\theta^*\}$. \nl
We define a forcing notion $P$:
\sn
$p \in P$ has the form $p = (u,u_*,v,v_*,\bar w) = (u^p,u^p_*,v^p,v^p_*,
\bar w^p)$ such that:
\mr
\item "{$(\alpha)$}"  $u_* \subseteq u \in [\lambda]^{< \kappa}$
\sn
\item "{$(\beta)$}"  $v_* \subseteq v \in [\lambda^*]^{< \kappa}$
\sn
\item "{$(\gamma)$}"  $\bar w = \bar w^p = 
\langle w_{\zeta,i}:\zeta \in v_* \text{ and } i < \theta^* \rangle =
\langle w^p_{\zeta,i}:\zeta \in v_*,i < \theta^* \rangle$
\sn
\item "{$(\delta)$}"  $w_{\zeta,i} \subseteq u_*$ and \nl
$b_i \cap b_j = \emptyset \Rightarrow w_{\zeta,i} \cap w_{\zeta,j} =
\emptyset$; this is toward being Hausdorff
\sn
\item "{$(\varepsilon)$}"  $\zeta \in v_* \Rightarrow A_\zeta \subseteq u$
\sn
\item "{$(\zeta)$}"  letting $A^p_\zeta = \cup\{w_{\zeta,i}:i < \theta^*\}
\cap A_\zeta$ for $\zeta \in v^p_*$ it has cardinality $\theta$ and for
simplicity even order type $\theta$ and for some
$\langle \gamma^p_{\zeta,j}:j < \theta \rangle$ list its members with no
repetitions we have \nl
$w^p_{\zeta,i} \cap A^p_\zeta = \{\gamma^p_{\zeta,j}:j < \theta \text{ and }
j \in b_i\}$
\sn
\item "{$(\eta)$}"  if $\zeta \in v^p_*,i < \theta^*$ and $\xi \in v^p_*$
\ub{then} the set ${\Cal U}^p_{\zeta,\xi,i}$ is an open subset (for part 
(2), clopen subset) of the space $Y^*$ where ${\Cal U}^p_{\zeta,\xi,i} =:
\{j < \theta:\gamma^p_{\xi,j} \in w^p_{\zeta,i}\}$.
\ermn
$\bigoplus \qquad$  \ub{convention} if $\zeta \in \lambda^* \backslash v^p_*$ 
we stipulate $w^p_{\zeta,i} = \emptyset$.
\mn
The \ub{order} is: $p \le q$ iff $u^p \subseteq u^q,u^p_* =u^q_* \cap u^p,
v^p \subseteq v^q,v^p_* = v^q_* \cap v^p$ and $\zeta \in v^p_* \Rightarrow
w^p_{\zeta,i} = w^q_{\zeta,i} \cap u^p$.
\mn
Clearly
\mr
\item "{$(*)_0$}"  $P$ is a partial order. 
\ermn
What is the desired space in $V^P$?  We define a $P$-name
${\underset\tilde {}\to X^*}$ as follows: \nl
set of points $\bigcup\{u^p_*:p \in {\underset\tilde {}\to G_P}\}$ \nl
The topology is defined by the following basis: \nl
$\{ \dbca_{\ell < n} {\underset\tilde {}\to {\Cal U}_{\zeta_\ell,i_\ell}}:
n < \omega,\zeta_\ell < \lambda^*,i_\ell < \theta^*\}$ where \nl
${\underset\tilde {}\to {\Cal U}_{\zeta,i}}
[{\underset\tilde {}\to G_P}] = \cup\{w^p_{\zeta,i}:p \in 
{\underset\tilde {}\to G_P},\zeta \in v^p_*\}$ \nl
(for part (2), also their compliments and hence their Boolean combinations)
\mr
\item "{$(*)_1$}"  for $\alpha < \lambda$ and $p \in P$ will have
$p \Vdash ``\alpha \in {\underset\tilde {}\to X^*}"$ iff $\alpha 
\in u^p_*$ and $p \Vdash ``\alpha \notin {\underset\tilde {}\to X^*}"$ 
iff $\alpha \in u^p_\alpha \backslash u^p_*$
\sn
\item "{$(*)_2$}"  $P$ is $\kappa$-complete, in fact if $\langle
p_\varepsilon:\varepsilon < \delta \rangle$ is increasing in $P$ and $\delta
< \kappa$ \ub{then} 
$p = \dbcu_{\varepsilon < \delta} p_\varepsilon$ is an upper
bound where $u^p = \dbcu_{\varepsilon < \delta} u^{p_\varepsilon},
u^p_* = \dbcu_{\varepsilon < \delta} u^{p_\varepsilon}_*,
v^p = \dbcu_{\varepsilon < \delta} v^{p_\varepsilon},v^p_* = 
\dbcu_{\varepsilon < \delta} v^{p_\varepsilon}_*$
and $w^p_{\zeta,i} = \cup\{w^{p_\varepsilon}_{\zeta,i}:\zeta \in
v^{p_\varepsilon}_*,\varepsilon < \delta\}$ \nl
[why?  straight]
\sn
\item "{$(*)_3$}"  $P' = \{p \in P:\text{if } \zeta < \lambda^* \text{ and }
|A_\zeta \cap u^p| \ge \sigma$ then $\zeta \in v^p\}$ is a dense subset of $P$
\nl
[why?  for any $p \in P$ we define by induction on $\varepsilon \le
\sigma^+:p_\varepsilon \in P$, increasingly continuous with $\varepsilon$.
Let $p_0 = p$, if $p_\varepsilon$ is defined, we define $p_{\varepsilon +1}$
by

$$
v^{p_{\varepsilon+1}} = \{\zeta < \lambda^*:\zeta \in v^{p_\varepsilon}
\text{ or } |A_\zeta \cap u^{p_\varepsilon}| \ge \sigma\}
$$

$$
v^{p_{\varepsilon +1}}_* = v^{p_\varepsilon}_*
$$

$$
u^{p_{\varepsilon +1}} = u^{p_\varepsilon} \cup \bigcup\{A_\zeta:
\zeta \in v^{p_{\varepsilon +1}}\}
$$

$$
u^{p_{\varepsilon +1}}_* = u^{p_\varepsilon}_* (= u^p_*)
$$

$$
w^{p_{\varepsilon +1}}_{\zeta,i} \text{ is: } w^{p_\varepsilon}_{\zeta,i}
\text{ if } \zeta \in v^{p_\varepsilon}_*,i < \theta^*
$$
(and there are no other cases). \nl
By assumption $(A)(ii)$, the set $v^{p_{\varepsilon +1}}$ has cardinality
$< \kappa$, so $p_{\varepsilon +1}$ belongs to $P$.
\ermn
Clearly $p_\varepsilon \le p_{\varepsilon +1} \in P$.
\mn
Now for $\varepsilon$ limit let $p_\varepsilon = \dbcu_{\xi < \xi} p_\xi$.
So we can carry the definition.  Now $p_{\sigma^+} = \dbcu_{\varepsilon <
\sigma} p_\varepsilon$ is as required because if $A_\zeta \in {\Cal A},
|A_\zeta \cap u^{p_{\sigma^+}}| \ge \sigma$ then for some $\varepsilon <
\sigma^+,|A_\zeta \cap u^{p_\varepsilon}| \ge \sigma$ hence $\zeta \in
v^{p_{\varepsilon +1}}$ hence $A_\zeta \subseteq u^{p_{\varepsilon +1}}
\subseteq u^{p_{\sigma^+}}$. \nl
Note that we use here $\sigma^+ < \kappa$.]
\mr
\item "{$(*)_4$}"  $P$ satisfies the $\kappa^+$-c.c. \nl
[Why?  Let $p_j \in P$ for $j < \kappa^+$, \wilog \,
$p_j \in P'$ for $j < \kappa^+$.  Now by the $\Delta$-system lemma for some
unbounded $S \subseteq \kappa^+$ and $v^\otimes \in [\lambda^*]^{< \kappa},
u^\otimes \in [\lambda]^{< \kappa}$ we have: \nl
$j \in S \Rightarrow v^\otimes \subseteq v^{p_j} \and u^\otimes \subseteq
u^{p_j} \text{ and } \langle v^{p_j} \backslash v^\otimes:j \in S \rangle$
are pairwise disjoint and $\langle u^{p_j} \backslash u^\otimes:j \in S
\rangle$ are pairwise disjoint.  Without loss of generality otp$(v^{p_j})$,
otp$(u^{p_j})$ are constant for $j \in S$ and any two $p_i,p_j$ are 
isomorphic over $v^\otimes,u^\otimes$ (if not clear see \scite{t.4}). \nl
Now for $j_1,j_2 \in S$ the condition $p_{j_1},p_{j_2}$ are 
compatible because of the following $(*)_5$]
\sn
\item "{$(*)_5$}"   assume $p^1,p^2 \in P$ satisfies
{\roster
\itemitem{ $(i)$ }  $v^{p^1}_* \cap (v^{p^2} \backslash v^{p^2}_*) 
= \emptyset$ and 
$u^{p^1}_* \cap (u^{p^2} \backslash u^{p^2}_*) = \emptyset$
\sn
\itemitem { $(ii)$ }  $v^{p^2}_* \cap (v^{p^1} \backslash v^{p^1}_*) = 
\emptyset$ and $u^{p^2}_* \cap (u^{p^1} \backslash u^{p^1}_*) =
\emptyset$
\sn
\itemitem{ $(iii)$ }  if $\zeta \in v^{p^1}_* \cap v^{p^2}_*$ then
$A^{p^1}_\zeta = A^{p^2}_\zeta$ and \nl
$i < \theta^* \Rightarrow w^{p^1}_{\zeta,i} \cap (u^{p^1} \cap 
u^{p^2}) = w^{p^2}_{\zeta,i} \cap (u^{p^1} \cap u^{p^2})$
\sn
\itemitem{ $(iv)_1$ }  if $\zeta \in v^{p^1}_* \backslash v^{p^2}_*$ then
$|A_\zeta \cap u^{p^2}| < \sigma$ or just $|A^{p^1}_\zeta \cap u^{p^2}| <
\sigma$
\sn
\itemitem{ $(iv)_2$ }  similarly \footnote{note that if $p^1,p^2 \in P'$, then
clauses $(iv)_1,(iv)_2$ holds automatically.} for 
$\zeta \in v^{p^2}_* \backslash v^{p^1}_*$ 
\endroster}
\ermn
\ub{then} there is $q \in P$ such that: 
\mr
\item "{$(a)$}"  $v^q = v^{p^1} \cup v^{p^2}$ 
\sn
\item "{$(b)$}"  $v^q_* = v^{p^1}_* \cup v^{p^2}_*$ 
\sn
\item "{$(c)$}"  $u^q = u^{p^1} \cup u^{p^2}$ 
\sn
\item "{$(d)$}"  $u^q_* = u^{p^1}_* \cup u^{p^2}_*$
\sn
\item "{$(e)$}"  $p^1 \le q$ and $p^2 \le q$.
\ermn
[Why?   To define the condition $q$ we just have to define
$w^q_{\zeta,i}$ (for $\zeta \in v^q_* = v^{p^1}_* \cup v^{p^2}_*$ 
and $i < \theta^*$).  If $\zeta \in v^{p^1}_* \cap v^{p^2}_*$ we let 
$w^q_{\zeta,i} = w^{p^1}_{\zeta,i} \cup w^{p^2}_{\zeta,i}$ for $i < \theta^*$.
\nl
Now for $\ell = 1,2$, let $v^{p^\ell}_* \backslash v^{p^{3 - \ell}}_*$
be listed as $\langle \Upsilon(\varepsilon,\ell):\varepsilon < 
\varepsilon^*_\ell \rangle$ with no repetitions such that 
$B^\ell_\varepsilon =: A^{p^\ell}_{\Upsilon(\varepsilon,\ell)} \cap 
(\dbcu_{\xi < \varepsilon} A^{p^\ell}_{\Upsilon(\xi,\ell)} \cup 
u^{p^{3 - \ell}})$ is of cardinality $< \sigma$. \nl
\sn
[Why possible?  By the assumption $(B)_2$ and clause $(iv)$ above.] \nl
Now for each
$\zeta \in v^{p^{3 - \ell}}_* \backslash v^{p^\ell}_*$ we choose by induction
on $\varepsilon < \varepsilon^*_\ell$ the sequence 
$\langle w^{\ell,\varepsilon}_{\zeta,i}:i < \theta^* \rangle$ such that
\nl
1)  $w^{\ell,\varepsilon}_{\zeta,i} \subseteq u^{p^{3 - \ell}} \cup 
\dbcu_{\xi < \varepsilon} A^{p^\ell}_{\Upsilon(\xi,\ell)}$.
\nl
2)  $w^{\ell,\varepsilon}_{\zeta,i}$ is increasingly continuous with
$\varepsilon$.
\nl
3) $w^{\ell,0}_{\zeta,i} = w^{p^{3 - \ell}}_{\zeta,i}$. \nl
4) $\varepsilon' < \varepsilon \Rightarrow w^{\ell,\varepsilon}_{\zeta,i}
\cap (u^{p^{3 - \ell}} \cup \dbcu_{\xi < \varepsilon_1} 
A^{p^\ell}_{\Upsilon(\xi,\ell)}) = w^{\ell,\varepsilon'}_{\zeta,i}$.
\nl
5) if $i < j < \theta^*$ and $b_i \cap b_j = \emptyset$ (hence
$w^{p^\ell}_{\zeta,i} \cap w^{p^\ell}_{\zeta,j} = \emptyset)$ then 
$w^{\ell,\varepsilon}_{\zeta,i} \cap w^{\ell,\varepsilon}
_{\zeta,j} = \emptyset$.
\nl
6) $\{j < \theta:\gamma^{p^\ell}_{\Upsilon(\varepsilon,\ell),j} \in
w^{\ell,\varepsilon +1}_{\zeta,i}\}$ is an open set in $Y^*$
(for part (2): clopen).]
\sn
For $\varepsilon = 0$ use clause (3) and for limit $\varepsilon$ take unions
(see clause (2)).  Suppose we have defined for $\varepsilon$ and let us
define for $\varepsilon +1$.  By an assumption above $B^\ell_\varepsilon$ has
cardinality $< \sigma$ and so $Z^\ell_\varepsilon = \{j < \theta:
\gamma^{p^\ell}_{\Upsilon(\varepsilon,\ell),j} \in B^\ell_\varepsilon\}$
is a subset of $\theta$ of cardinality $< \sigma$.  Hence, by assumption (E),
we can find a sequence $\langle t_j(\varepsilon,\ell):j \in Z^\ell
_\varepsilon \rangle$ such that: $t_j(\varepsilon,\ell) < \theta^*$ and
$j \in b_{t_j(\varepsilon,\ell)}$ for $j \in Z^\ell_\varepsilon$ and
$\langle b_{t_j(\varepsilon,\ell)}:j \in Z^\ell_\varepsilon \rangle$ is a
sequence of pairwise disjoint subsets of $Y^*$.  \nl
Lastly, we let

$$
\align
w^{\ell,\varepsilon +1}_{\zeta,i} = w^{\ell,\varepsilon}_{\zeta,i}
\cup \bigl\{ \gamma^{p^\ell}_{\Upsilon(\varepsilon,\ell),s}:
&\text{ for some } j \in Z^\ell_\varepsilon \text{ we have}: \\
  &\,\gamma^{p^\ell}_{\Upsilon(\varepsilon,\ell),t_j(\varepsilon,\ell)}
\in w^{\ell,\varepsilon}_{\zeta,i} \text{ and} \\
  &\,s \in b_{t_j(\varepsilon,\ell)} \bigr\}.
\endalign
$$
\mn
Clearly this is O.K. and we are done.  Remember that the union of
$< \sigma$ set from ${\Cal B}$ is clopen for part (2).]
\mr
\item "{$(*)_6$}"  in $(*)_5$ if in addition for $\ell =1,2$ we have
$Z_\ell \subseteq u^{p^\ell} \backslash u^{p^{3 - \ell}}$ such that
$(\forall \zeta \in v^{p^\ell}_*)[|A^{p^\ell}_\zeta \cap Z_\ell| < \sigma]$
\ub{then} we may add to the conclusion: \nl
$\ell \in \{1,2\},\zeta \in v^{p^{3 - \ell}}_* \backslash v^{p^\ell}_*,
i < \theta^* \Rightarrow w^q_{\zeta,i} \cap Z_\ell = \emptyset$. \nl
More generally if $g_\ell:(v^{p^{3 - \ell}}_* \backslash v^{p^\ell}_*) \times
\theta^* \times Z_\ell \rightarrow \{0,1\}$ we can add \nl
$\ell \in \{1,2\},\zeta \in v^{p^{3 - \ell}}_* \backslash v^{p^\ell}_*,
i < \theta^*,\gamma \in Z_\ell \Rightarrow [\gamma \in w^q_{\zeta,i} 
\leftrightarrow g_\ell(\zeta,i,\gamma) = 1]$. \nl
[Why?  During the proof of $(*)_5$ when for 
$\zeta \in v^{p^{3 - \ell}}_* \backslash v^{p^\ell}_*$, we
define $\langle w^{\ell,\varepsilon}_{\zeta,i}:i < \theta^* \rangle$ by
induction on $\varepsilon$ we add \nl
(7)  $i < \theta^*,\gamma \in Z_\ell \cap (u^{p^{3 - \ell}} \cup \dbcu
_{\xi < \varepsilon} A^{p^\ell}_{\Upsilon(\xi,\ell)})$ implies $\gamma \in
w^{\ell,\varepsilon}_{\zeta,i} \leftrightarrow g_\ell(\zeta,i,\gamma)=1$.  In
the proof when we use clause (E), instead of using $B^\ell_\varepsilon = 
A^{p^\ell}_{\zeta(\varepsilon,\ell)} \cap (\dbcu_{\xi < \varepsilon} 
A^{p^\ell}_{\zeta(\xi,\ell)} \cup u^{p^{3 - \ell}})$ we use 
$B^\ell_\varepsilon = A^{p^\ell}_{\zeta(\varepsilon,
\ell)} \cap (\dbcu_{\xi < \varepsilon} A^{p^\ell}_{\zeta(\xi,\ell)} \cup
u^{p^{3 - \ell}} \cup Z_\ell)$ which still has cardinality $< \sigma$.]
\mn
Now we come to the main point
\sn
\item "{$(*)_7$}"  in $V^P$, if $i(*) < \text{ cf}(\theta)$ and
$X^* = \dbcu_{i < i(*)} X_i$ \ub{then} some closed $Y \subseteq X^*$ 
is homeomorphic to $Y^*$.
\ermn
[Why?  Toward contradiction assume $p^* \in P$ and $p^* \Vdash_P ``\langle
{\underset\tilde {}\to X_i}:i < i(*) \rangle$ is a counterexample to
$(*)_7$". \nl
Without loss of generality 
$p^* \Vdash_P ``\langle {\underset\tilde {}\to X_i}:i < i(*) \rangle$ is a
partition of $X^*$, i.e. of $\bigcup\{u^p_*:p \in
{\underset\tilde {}\to G_P}\}$". \nl
For each $\alpha < \lambda$ let $\langle(p_{\alpha,j},i_{\alpha,j}):j <
\kappa \rangle$ be such that:
\mr
\widestnumber\item{$(iii)$}
\item "{$(i)$}"  $\langle p_{\alpha,j}:j < \kappa \rangle$ is a maximal
antichain of $P$ above $p^*$
\sn
\item "{$(ii)$}"  $p_{\alpha,j} \Vdash_P ``\alpha \in X_{i_{\alpha,j}}"$,
so $i_{\alpha,j} < i(*)$ and $\alpha \in u^{p_{\alpha,j}}_*$
\sn
\item "{$(iii)$}"  $p^* \le p_{\alpha,j}$.
\ermn
Now choose a function $F$, Dom$(F) = \lambda$ as follows:

$$
F(\alpha) \text{ is } \bigcup\{u^{p_{\alpha,j}}:j < \kappa\}.
$$
\mn
So we can find $\zeta(*) < \lambda^*$ and $A \subseteq A_{\zeta(*)}$ of
order type $\theta$ such that: if $\alpha \ne \beta$ are from $A$ then 
$\alpha \notin F(\beta)$.  Let $A = \{\beta_\varepsilon:\varepsilon < 
\theta\}$ with no repetitions.  Now we 
shall choose by induction on $\varepsilon \le \theta,p_\varepsilon,
g_\varepsilon$ and if
$\varepsilon < \theta$ also $j_\varepsilon < \kappa$ such that:
\roster
\item "{$(a)$}"  $p_\varepsilon \in P$ \nl
$u^{p_\varepsilon} = 
u^{p^*} \cup \dbcu_{\varepsilon(1) < \varepsilon} 
u^{p_{\beta_{\varepsilon(1)},j_{\varepsilon(1)}}}$ \nl
$u^{p_\varepsilon}_* = 
u^{p^*}_* \cup \dbcu_{\varepsilon(1) < \varepsilon} 
u^{p_{\beta_{\varepsilon(1)},j_{\varepsilon(1)}}}_*$ \nl
$v^{p_\varepsilon} = v^{p^*} \cup \dbcu_{\varepsilon(1) < \varepsilon}
v^{p_{\beta_{\varepsilon(1)},j_{\varepsilon(1)}}}$ \nl
$v^{p_\varepsilon}_* = v^{p^*}_* \cup \dbcu_{\varepsilon(1) < \varepsilon}
v^{p_{\beta_{\varepsilon(1)},j_{\varepsilon(1)}}}_*$ \nl
(so $p_0 = p^*$)
\sn
\item "{$(b)$}"  $j_\varepsilon = 
\text{ Min}\{j < \kappa:p_{\beta_\varepsilon,j}$ is compatible with 
$p_\varepsilon\}$
\sn
\item "{$(c)$}"  $g_\varepsilon$ is a function, increasing with
$\varepsilon$, from $v^{p_\varepsilon}_*
\times \theta^*$ into the family  of open subsets of $Y^*$ (for part (2),
clopen)
\sn
\item "{$(d)$}"  if $b_{i_1} \cap b_{i_2} = \emptyset$ then $g_\varepsilon
(\zeta,i_1) \cap g_\varepsilon(\zeta,i_2) = \emptyset$ (if defined)
\sn
\item "{$(e)$}"  letting $\Upsilon_\varepsilon = \text{ otp}\{\xi <
\varepsilon:i_{\beta_\varepsilon,j_\varepsilon} = i_{\beta_\xi,j_\xi}\}$ we
have for every $\zeta \in v^{p_\varepsilon}_*$ and $i < \theta^*$ and $\xi <
\varepsilon$:
$$
\beta_\xi \in w^{p_\varepsilon}_{\zeta,i} \Leftrightarrow
\Upsilon_\xi \in g_\varepsilon(\zeta,i)
$$
\sn
\item "{$(f)$}"  $p_\varepsilon$ is increasing continuous.
\ermn
No problem to carry the definition.  As for $\varepsilon$ successor, for 
this $(*)_6$ was prepared.  In limit $\varepsilon$ take union.  In all 
cases $j_\varepsilon$ is well defined by clause $(i)$ above.  Let 
$i^* < i(*)$ be minimal such that the set $Z = \{\varepsilon < \theta:
i_{\beta_\varepsilon,j_\varepsilon} = i^*\}$ has cardinality $\theta$.  
Note: $\zeta(*) \notin
v^{p_{\beta_\varepsilon,j}}$ as $A \cap F(\beta_\varepsilon)$ is a
singleton so $|A \cap u^{p_{\beta_\varepsilon,j}}| \le 1$ and 
$p_{\beta_\varepsilon,j} \in P'$.  Now we define $p$:

$$
u^p = u^{p_\theta}
$$

$$
u^p_* = u^{p_\theta}_*
$$

$$
v^p = v^{p_\theta} \cup \{\zeta(*)\}
$$

$$
v^p_* = v^{p_\theta}_* \cup \{\zeta(*)\}
$$

$$
A^p_{\zeta(*)} = \{\beta_\varepsilon:\varepsilon \in
Z\} \text{ and } \gamma^p_{\zeta(*),\varepsilon} \text{ is the } 
\varepsilon \text{-th member of } A^p_{\zeta(*)}
$$
\mn
$w^p_{\zeta,i}$ is
\mr
\item "{$(\alpha)$}"  $w^{p_\theta}_{\zeta,i}$ \ub{if} $\zeta \in 
v^{p_\theta}$
\sn
\item "{$(\beta)$}"  $\{\beta_\varepsilon:\varepsilon \in Z \text{ and otp}
(Z \cap \varepsilon) \in b_i\}$ \ub{if} $\zeta = \zeta(*)$.
\ermn
We can easily check that $p \in P$ and $p^* \le 
p_{\beta_\varepsilon,j_\varepsilon} \le p \in P$ (but we do not ask
$p_\varepsilon \le p$).
Clearly $p$ forces that $\{\beta_\varepsilon:\varepsilon \in Z\}$ is
included in one $X_i$. \nl
Let $g:\theta \rightarrow \lambda$ be $g(\xi) = \beta_\varepsilon$ when
$\xi < \theta,\varepsilon \in Z$, otp$(Z \cap \varepsilon) = \xi$.
Now $p \ge p^*$ and we are done by $(*)_8$ below.]
\mr
\item "{$(*)_8$}"  if $p \in P$ and $\zeta \in v^p_*$ then \nl
$p \Vdash$ ``the mapping $j \mapsto \gamma^p_{\zeta,j}$ for $j < \theta$ is
a homeomorphism from $Y^*$ onto the closed subspace
$\underset\tilde {}\to X \restriction \{\gamma^p_{\zeta,j}:j < \theta\}$ of
$\underset\tilde {}\to X$" \nl
[Why?  Let $p \in G,G \subseteq P$ is generic over $V$.
{\roster
\itemitem{ $(\alpha)$ } If $b \in {\Cal B}$, then for 
some open set ${\Cal U}$ of $\underset\tilde {}\to X$ (clopen for
part (2)) we have
$$
{\Cal U} \cap \{\gamma^p_{\zeta,j}:j < \theta\} = \{\gamma^p_{\zeta,j}:j \in
b\}
$$
[Why?  As $b = b_i$ for some $i < i(*)$ and $p$ forces that \nl
$w_{\zeta,i} \cap \{\gamma^p_{\zeta,j}:j < \theta\} = \{\gamma^p_{\zeta,j}:
j \in b_i\}$.]
\sn
\itemitem{ $(\beta)$ }  If $b$ is an open set for $Y^*$, then for some open
subset ${\Cal U}$ of $\underset\tilde {}\to X$ we have
$$
{\Cal U} \cap \{\gamma^p_{\zeta,j}:j < \theta\} = \{\gamma^p_{\zeta,j}:
j \in b\}
$$
[Why?  As $b = \dbcu_{i \in Z} b_i$ for some $Z \subseteq \theta^*$ 
and apply clause $(\alpha)$]
\sn
\itemitem{ $(\gamma)$ }  if ${\Cal U}$ is an open subset of
$\underset\tilde {}\to X$ and $\gamma^p_{\zeta,j(*)} \in {\Cal U}$
(so $\zeta \in u^p_*$), then for some $i(*) < \theta^*$ we have
$$
\gamma^p_{\zeta,j(*)} \in w^p_{\zeta,i(*)} \cap
\{\gamma^p_{\zeta,j}:j < \theta\} \subseteq
{\underset\tilde {}\to {\Cal U}_{\zeta,i(*)}} \cap
\{\gamma^p_{\zeta,j}:j < \theta\} \subseteq {\Cal U}.
$$
[Why?  By the definition of the topology $\underset\tilde {}\to X$ we can
find $n < \omega,\xi_\ell < \lambda^*$ and $i_\ell < \theta^*$ for $\ell < n$
such that
$\gamma^p_{\zeta,j(*)} \in \dbca_{\ell < n} 
{\underset\tilde {}\to {\Cal U}_{\xi_\ell,i_\ell}}[G] \subseteq {\Cal U}$.  
We can find $q \in P$ such that $p \le q$ and $\xi_\ell \in v^q_*$ for
$\ell < n$.  For each $\ell < n$,
by clause $(\eta)$ in the definition of $P$ we have
${\Cal U}^q_{\zeta,\xi_\ell,j_\ell}$ is an open set
for $Y^*$, and necessarily $j(*) \in {\Cal U}^q_{\zeta,\xi_\ell,j_\ell}$.  Let
$i(*)$ be such that $j(*) \in b_{i(*)} \subseteq \dbca_{\ell < n}
{\Cal U}^q_{\zeta,\xi_\ell,j_\ell}$ hence
$\gamma^p_{\zeta,j(*)} \in {\underset\tilde {}\to {\Cal U}_{\zeta,i(*)}}[G]
\cap \{\gamma^p_{\zeta,j}:j < \theta\} \subseteq \dbca_{\ell < n}
{\underset\tilde {}\to {\Cal U}_{\xi_\ell,j_\ell}}[G] \subseteq {\Cal U}$ as
required.  So $i(*)$ is as required.]
\sn
\itemitem{ $(\delta)$ }  $\{\gamma^p_{\zeta,j}:j < \theta\}$ is a closed
subset of $\underset\tilde {}\to X$ \nl
[Why?   Let $\beta \in \lambda \backslash \{\gamma^p_{\zeta,j}:j < \theta\}$
and let $p \le q \in P$; it suffices to find $q^+,q \le q^+ \in P$ and 
$\xi \in v^{q^+}_*$ and $i < \theta^*$ such that 
$\beta \in u^{q^+} \backslash u^{q^+}_*$ or $\beta \in w^{q^+}_{\xi,i}$ and
$w^{q^+}_{\xi,i} \cap \{\gamma^p_{\zeta,j}:j < \theta\} = \emptyset$.
If $\beta \notin u^q_*$ define $q^+$ like $q$ except that $u^{q^+} = u^q \cup
\{\beta\}$ (but $u^{q^+}_* = u^q_*$).  So without loss of generality 
$\beta \in u^q_*$. \nl
We can find a set $u \subseteq u^q_*$ such that $\beta \in u,A^q_\zeta \cap
u = \emptyset$ and $\zeta' \in v^q_* \Rightarrow \{j < \theta:
\gamma^q_{\zeta',j} \in u\}$ is an open subset of $Y$, (just as in the 
proof of $(*)_5$; for part (2) we ask ``clopen subset of $Y$").
By $\otimes_1$ we can find $\xi \in \lambda^* \backslash v^q$ such that
$\{\emptyset\} = A_\xi \cap u^q_*$ (why? apply $\otimes_1$ with $\alpha <
\beta \in \lambda \backslash u^q$ and $B = u^q$) and let
$\gamma_{\varepsilon,i} \in A_\xi$ for $i < \theta$ be increasing.  We define
$q^+$ as follows.
\endroster}
\ermn

$$
v^{q^+} = v^q \cup \{\xi\}
$$

$$
v^{q^+}_* = v^q_* \cup \{\xi\}
$$

$$
u^{q^+} = u^q \cup A_\xi
$$

$$
u^{q^+}_* = u^q_* \cup \{\gamma_{\xi,j}:j < \theta\}
$$

$w^{q^+}_{\zeta,i} \text{ is } w^q_{\zeta,i} \text{ if } \zeta \in v^q_*
\text{ and is } \{\gamma_{\xi,j}:j \in b_i\} \cup u$ if $\zeta = \xi \and 
0 \in b_i \text{ and is } \{\gamma_{\xi,j}:j \in b_i\}$ if $\zeta = \xi 
\and 0 \notin b_i$.]
\mn
Lastly, we would like to know that 
$\underset\tilde {}\to X$ is a Hausdorff space.  We prove more
\mr
\item "{$(*)_9$}"  In $V^P$ if $u_1 \subseteq u_2 \in [\lambda]^{< \sigma}$
then for some $\zeta,i$ we have
$$
w_{\zeta,i} \cap u_2 \cap \underset\tilde {}\to X = u_1
\cap \underset\tilde {}\to X
$$
[Why?  Let $p_0 \in P$ force that ${\underset\tilde {}\to u_1} \subseteq
{\underset\tilde {}\to u_2}$ form a counterexample, as 
$P$ is $\kappa$-complete
some $p_1 \ge p_0$ forces ${\underset\tilde {}\to u_1} = u_2,
{\underset\tilde {}\to u_2} = u_2$ and $p_1 \in P'$.  Necessarily $u_2
\subseteq u^{p_1}_*$, as in the proof of $(*)_8(\delta)$. \nl
Let $\zeta(*) \in \lambda^* \backslash v^{p_1}$ be such that $A_{\zeta(*)}
\cap u^{p_1} = \emptyset$ (as in the proof of $(*)_8(\delta)$).  Let
$\gamma_{\zeta(*),j} \in A_{\zeta(*)}$, for $j < \theta$ be increasing.  Let
$u \subseteq u^{p_1}_*$ be such that $u \cap u_2 = u_1$ and $\zeta' \in
v^{p_1}_* \Rightarrow \{j < \theta:\gamma^{p_1}_{\zeta',j} \in u\}$ is clopen
in $Y$ (exists as in the proof of $(*)_5$) and define $q \in P$:
\endroster

$$
u^q = u^{p_1} \cup u_2
$$

$$
u^q_* = u^{p_1}_* \cup (u_2 \backslash u^{p_1})
$$

$$
v^q = v^p \cup \{\zeta(*)\}
$$

$$
v^q_* = v^{p_1} \cup \{\zeta(*)\}
$$

$w^q_{\zeta,i} \text{ is: } w^{p_1}_{\zeta,i} \text{ if } \zeta \in v^q,
\text{ is } \{\gamma_{\zeta(*),j}:j \in b_i\} \cup u \text{ if } 
\zeta = \zeta(*) \and 0 \in b_i \text{ and is } \{\gamma_{\zeta(*),j}:j \in
b_i\} \text{ if } \zeta = \zeta(*) \and 0 \notin b_i$.]
\sn
Together all is done.   \hfill$\square_{\scite{t.2}}$
\enddemo
\bn
\centerline {$* \qquad * \qquad *$}
\bn
Now when are the assumptions of \scite{t.2} hold?
\proclaim{\stag{t.3} Claim}  Assume
\mr
\item "{$(a)$}"   ${\frak a} \in [\text{Reg } \cap \mu \backslash 
\kappa]^\theta,J = [{\frak a}]^{< \sigma}$
\sn
\item "{$(b)$}"   $\Pi{\frak a}/J$ is $(\lambda^*)^+$-directed,
\sn
\item "{$(c)$}"   $\lambda \ge \mu$ is singular,
\sn
\item "{$(d)$}"  $\lambda^* > \lambda > \kappa^{< \kappa} = 
\kappa > \theta$; 
\sn
\item "{$(e)$}"  $\lambda^* < 2^\lambda$ is regular.
\ermn
\ub{Then}
\mr
\item "{$(f)$}"  In $V_1 = V^{\text{Levy}(\lambda^*,2^\lambda)}$ we have 
(a),(c),(d) and (e) and $2^\lambda = \lambda^*$ and
\sn
\item "{$(g)$}"   the assumptions $(A)(i),(B_1),(B_2),(C)$ of Theorem 
\scite{t.2} hold (recall $(A)(i)$ means we omit $\theta^*$ and $(\forall
\alpha < \kappa)(|\alpha|^\sigma < \kappa)$.
\endroster
\endproclaim
\bigskip

\demo{Proof}  Let ${\frak a} = \{\lambda_\varepsilon:\varepsilon < \theta\}$
without repetitions; \wilog \, $\lambda_i > \kappa^{++}$.  
Let $J' = [\theta]^{< \sigma}$.  By \cite[Ch.II,1.4]{Sh:g} (at least the proof,
see below) in $V$ we can find $\langle f_\alpha:\alpha < \lambda^* \rangle$ 
such that (more but irrelevant here)
\mr
\item "{$(*)_0$}"  Assume $J'$ is an ideal on ${\frak a},\Pi{\frak a}/J'$ is
$(\lambda^*)^+$-directed, $\lambda^* > \sup({\frak a})$.  \ub{Then} 
we can find $\langle f_\alpha:\alpha < \lambda^* \rangle$ such that
$f_\alpha \in \dsize \prod_{\varepsilon < \theta}
\lambda_\varepsilon$ and for every $Z \in [\lambda^*]^{< \kappa}$ for
some sequence $\bar a = \langle a_\alpha:\alpha \in Z \rangle$ such that 
$a_\alpha \in J'$ for $\alpha \in Z$ and some well ordering $<^*$ of $Z$
we have
{\roster
\itemitem{ $(i)$ }  $\alpha_1 \in Z \and \alpha_2 \in Z \and \varepsilon_1 <
\theta \and \varepsilon_2 < \theta \and f_{\alpha_1}(\varepsilon_1) = 
f_{\alpha_2}(\varepsilon_2) \rightarrow \varepsilon_1 = \varepsilon_2$
\sn
\itemitem{ $(ii)$ }  $\alpha \in Z \and \beta \in Z \and \alpha <^* \beta \and
\varepsilon \in \theta \backslash a_\beta \Rightarrow
f_\alpha(\varepsilon) \ne f_\beta(\varepsilon)$.
\endroster}
[Why?  There we get only: for some $\langle f_\alpha:\alpha < \lambda^*
\rangle$ as above with $(i) + (ii)$ replaced by:
for every $Z \in [\lambda^*]^{< \kappa}$ (even
$Z \in [\lambda]^{\mu'}$ if ${\frak a} \subseteq \mu'$ has order type 
$\le \sigma$)
we can find $\bar a = \langle a_\alpha:\alpha \in Z \rangle$ such that
$a_\alpha \in J'$ and $\alpha \ne \beta \and \alpha \in Z \and \beta \in Z
\and \varepsilon \in \theta \backslash a_\alpha \backslash a_\beta
\Rightarrow f_\alpha(\varepsilon) \ne f_\beta(\varepsilon)$.  Clause (i) is
easy, just replace $f_\alpha$ by $f'_\alpha$ which is defined by 
$f'_\alpha(\lambda_\varepsilon) =
\theta \times f_\alpha(\lambda_\varepsilon) + \varepsilon$.  We shall
prove $\langle f_\alpha:\alpha < \lambda^* \rangle$ is as required.  So let
$Z \in [\lambda^*]^{< \kappa}$.  
We can choose by induction on $\zeta,Z_\zeta \subseteq Z$ 
increasingly continuous
in $\zeta$ such that $Z_0 = \emptyset,|Z_{\zeta +1} \backslash Z_\zeta| \le
\sigma,[Z_\zeta \ne Z \Rightarrow Z_\zeta \ne Z_{\zeta +1}]$ and
$\alpha \in Z_{\zeta +1} \backslash Z_\zeta \and \varepsilon \in a_\alpha
\and \beta \in Z \and \varepsilon \notin a_\beta \and f_\alpha(\varepsilon) =
f_\beta(\varepsilon) \rightarrow \beta \in Z_{\zeta +1}$.
As $|a_\alpha| < \sigma$ and $(\forall \varepsilon <
\theta)(\forall \gamma)(\exists^{\le 1} \beta)(\varepsilon \notin a_\beta
\and f_\beta(\varepsilon) = \gamma)$ there is no problem, (in fact if
$\sigma$ is regular we can ask $< \sigma$).  Now list $Z$ as
$\langle \alpha_\xi:\xi < \xi^* \rangle$ such that
$\{\xi:\alpha_\xi \in Z_{\zeta +1} \backslash Z_\zeta\}$ is a convex set
of order type of its cardinality so $\le \sigma$, which is above 
$\{\xi:\alpha_\xi \in Z_\zeta\}$ and $\alpha_\xi \in Z_{\zeta +1} \backslash
Z_\zeta \Rightarrow \sigma > |\cup\{a_\varepsilon:\varepsilon \in
Z_{\zeta +1} \backslash Z_\zeta$ and $\alpha_\varepsilon \le \alpha_\xi\}|$.
Define a well ordering $<^*$ of $Z$ by $\alpha_\varepsilon <^* 
\alpha_\varepsilon \equiv \xi < \varepsilon$.  Now
we define $a'_\alpha \in J'$ for $\alpha \in Z$ as follows: if $\alpha =
\alpha_\xi \in Z_{\zeta +1} \backslash Z_\zeta$ then $a'_\alpha = \cup
\{a_{\alpha_\varepsilon}:\alpha_\varepsilon \in Z_{\zeta +1} \backslash
Z_\zeta$ and $\varepsilon \le \xi\}$.  Now suppose 
for $\langle a'_\alpha:\alpha \in Z \rangle$ fails clause (ii), so there are
$\varepsilon < \theta$ and $\alpha < \beta$ from $Z$ which exemplifies this.
As $a_\alpha \subseteq a'_\alpha$ and the choice of $\langle a_\alpha:\alpha 
\in Z \rangle$ necessarily $\varepsilon \in a_\alpha$ and we get easy
contradiction.]
\ermn
Clearly in $V_1$ we have (a),(c),(d) and $(*)_0$ above (and we can forget
$V$ and (b), recall (b) says ``$\Pi{\frak a}$ is $(\lambda^*)^+$-directed", 
on the existence of $\bar f$ as in $(*)_0$, see \cite[Ch.VIII,\S5]{Sh:g}). \nl
We can in $V_1$ list $\langle h_\alpha:\alpha < \lambda^* \rangle$ the
functions $h:\lambda \rightarrow [\lambda]^\kappa$.
Now for each $\zeta < \lambda^*$ we define a function 
$g_\zeta:\kappa^{++} \rightarrow [\kappa^{++}]^{\le \kappa}$ by

$$
\align
g_\zeta(\gamma) = \biggl\{ \beta < \kappa^{++}:&\text{ for some }
\varepsilon_1,\varepsilon_2 < \theta \text{ we have} \\
  &f_\alpha(\varepsilon_1) \times \kappa^{++} \times \theta +
\beta \times \theta + \varepsilon_1 \in \\
  &h_\zeta[f_\zeta(\varepsilon_2) \times \kappa^{++} \times \theta +
\gamma \times \theta + \varepsilon_2] \biggr\}.
\endalign
$$
\mn
So we can for each $\zeta < \lambda^*$ find $Z_\zeta \in 
[\kappa^{++}]^{\kappa^{++}}$ such that

$$
\beta_1 \ne \beta_2 \in Z_\zeta \Rightarrow \beta_1 \notin g_\zeta(\beta_2).
$$
\mn
For $\zeta < \lambda^*$ let $A_\zeta = 
\{f_\zeta(\varepsilon) \times \kappa^{++} \times \theta +
\beta \times \theta + \varepsilon:\varepsilon < \theta \text{ and }
\beta < \kappa^{++} \text{ is the } \varepsilon \text{-th member of }
Z_\zeta\}$. \nl
Now we shall check.  
\enddemo
\bn
Let ${\Cal A} = \{A_\zeta:\zeta < \lambda^*\}$.  Clearly
\mr
\item "{$(*)_1$}"  $A_\zeta \in [\lambda]^\theta$ (hence ${\Cal A}
\subseteq [\lambda]^\theta$)
\sn
\item "{$(*)_2$}"  $\zeta_1 \ne \zeta_2 \Rightarrow |A_{\zeta_1} \cap
A_{\zeta_2}| < \sigma$ \nl
[Why?  Let $\alpha \in A_{\zeta_1} \cap A_{\zeta_2}$ so for some $\ell =1,2$
we have $\alpha = f_{\zeta_\ell}(\varepsilon_\ell) \times \kappa^{++} \times
\theta + \beta_\ell \times \theta + \varepsilon_\ell$ with $\beta_\ell <
\kappa^{++},\varepsilon_\ell < \theta$.  Clearly this implies $\varepsilon_0
= \varepsilon_1,\beta_1 = \beta_2,f_{\zeta_1}(\varepsilon_\ell) =
f_{\zeta_2}(\varepsilon_\ell)$, 
and otp$(\beta_\ell \cap Z_{\zeta_\ell}) = \varepsilon_\ell$, so $\beta_\ell$
depends just on $\zeta_\ell$ and $\varepsilon_\ell$ (not on $\alpha$) and by
$(i)$ of $(*)_0$ also $\varepsilon_\ell$ is determined by $\zeta_\ell,
f_{\zeta_\ell}(\varepsilon_\ell)$ hence $|A_{\zeta_1} \cap A_{\zeta_2}| <
|\{\varepsilon < \theta:f_{\zeta_1}(\varepsilon) = f_{\zeta_2}(\varepsilon)
\}| < \sigma$ as $\zeta_1 < \zeta_2 \rightarrow f_{\zeta_1}
<_J f_{\zeta_2}$ recalling $J = [\theta]^{< \sigma}$.]
\sn
\item "{$(*)_3$}"  $|{\Cal A}| = \lambda^*$ \nl
[Why?  By the choice of ${\Cal A}$ and $(*)_1 + (*)_2$.]
\sn
\item "{$(*)_4$}"  if $F:\lambda \rightarrow [\lambda]^{\le \kappa}$,
\ub{then} some $A \in {\Cal A}$ is $F$-free \nl
[Why?  For some $\alpha$ we have $F = h_\alpha$, so $Z_\alpha,A_\alpha$ were
chosen to make this true.]
\sn
\item "{$(*)_5$}"  if ${\Cal A}' \in [{\Cal A}]^{< \kappa}$, then we can
list ${\Cal A}'$ as $\{A_{\zeta_i}:i < i(*)\}$ such that \nl
$|A_{\zeta_i} \cap \dbcu_{j < i} A_{\zeta_j}| < \sigma$ \nl
[Why?  Let ${\Cal A}' = \{A_\zeta:\zeta \in Z\}$ where $Z \subseteq
\lambda^*,|Z| < \kappa$, so by $(*)_0$ we can find $\langle a_\alpha:\alpha
\in Z \rangle,<^*$ as there.  Let 
$Z = \{\zeta_i:i < i(*)\}$ be $<^*$-increasing with $i$ and so
$$
A_{\zeta_i} \cap \dbcu_{j < i} A_{\zeta_j} \subseteq
\{f_{\zeta_i}(\varepsilon) \times \kappa^{++} \times \theta + \theta \times
\beta_{\zeta_i,\varepsilon} + \varepsilon:\varepsilon \in a_{\zeta_i}\}
$$
which has cardinality $< \sigma$ where 
$\beta_{\zeta_i,\varepsilon}$ is the $\varepsilon$-th member of
$Z_\alpha$.]
\ermn
So clause (A)(i) holds by our assumption (note, $\theta^*$ does not appear
here), clause $(B)_1$ holds by $(*)_1 + (*)_2$ and $(B)_2$ holds by $(*)_5$
and lastly (C) holds by $(*)_4$.  \hfill$\square_{\scite{t.3}}$
\bn
\proclaim{\stag{t.4} Claim}  We can weaken the assumption \scite{t.2} 
omitting in part (2) the ``closed under union of $< \sigma$" and by omitting
(A)(ii) and by replacing (E) by (E)$^-$, i.e. having:
\mr
\item "{$(A)(i)$}"  $\lambda > \kappa > \theta \ge \sigma \ge \aleph_0$ and
$\kappa = \kappa^{< \kappa}$ \nl
(i.e. this is $(A)(i)$ without $(A)(ii)$ i.e. omitting 
``$(\forall \alpha < \kappa)(|\alpha|^\sigma <
\kappa),\kappa > \theta^* \ge \theta$")
\sn
\item "{$(E)^-$}"  if $Y_0,Y_1$ are disjoint subsets of $Y^*$ each with
$< \sigma$ points, \ub{then} there are open disjoint sets ${\Cal U}_0,
{\Cal U}_1$ of $Y^*$ such that $Y_0 \subseteq {\Cal U}_0,Y_\theta 
\subseteq Y_1$.
\endroster
\endproclaim
\bigskip

\demo{Proof}  We indicate the changes.

We can further demand from $\langle b_i:i < \theta^* \rangle$ that
\mr
\item "{$\boxtimes_1$}"  $b_{2i} \cap b_{2i+1} = \emptyset$ and if
$b_{i_0} \cap b_{i_1} = \emptyset$ then for some $j$ we have
$(b_{2j},b_{2j+1}) = (b_{i_0},b_{i_1})$.
\ermn
In the defintion of $P$ we replace clause $(\delta)$ by
\mr
\item "{$(\delta)^-$}"  $w_{\zeta,i} \subseteq u_*$ and $w_{\zeta,2i} \cap
w_{\zeta,2i+1} = \emptyset$.
\ermn
However, as we have weakened assumption $(A)$, the $\kappa^+$-c.c. may fail.
So we define: we say $(f,g)$ is an isomorphism from $p \in P$ onto $q \in P$
if:
\mr
\widestnumber\item{$(iii)$}
\item "{$(i)$}"  $f$ is a one-to-one mapping from $u^p$ onto $u^q$
\sn
\item "{$(ii)$}"  $g$ is a one-to-one mapping from $v^p$ onto $v^q$
\sn
\item "{$(iii)$}"  $f$ maps $u^p_*$ onto $u^q_*$
\sn
\item "{$(iv)$}"  $g$ maps $v^p_*$ onto $v^q_*$
\sn
\item "{$(v)$}"  if $\zeta \in v^p_*$ then $A_{g(\zeta)} = \{f(\beta):
\beta \in A_\zeta\}$
\sn
\item "{$(vi)$}"  if $\zeta \in v^p_*$ and $j < \theta$ then
$\gamma^q_{g(\zeta),j} = f(\gamma^p_{\zeta,j})$
\sn
\item "{$(vii)$}"  if $\zeta \in v^p_*$ and $i < \theta^*$ then
$$
w^q_{g(\zeta),i} = \{f(\beta):\beta \in w^p_{\zeta,i}\}.
$$
\ermn
We say $p,q$ are isomorphic if such $(f,g)$ exists.  Clearly being isomorphic
is an equivalent relation.  Let $\chi$ be large enough and ${\frak C}$ be an
elementary submodel of $({\Cal H}(\chi),\in,<^*)$ of cardinality $\kappa$
such that $\lambda,\kappa,\theta^*,\theta,\sigma,Y^*,\langle b_i:i <
\theta^* \rangle,{\Cal A},P$ belong to ${\frak C}$ and 
${}^{\kappa >}{\frak C} \subseteq {\frak C}$.  Let

$$
Q = \{p \in P:\text{for some } q \in P \cap {\frak C} \text{ we have }
p,q \text{ are isomorphic}\}.
$$
\mn
In the rest of the proof $P$ is replaced by $Q$, each time we construct a
condition we have to check if it belongs to $Q$.

The only place we use $(\forall \alpha < \kappa)(|\alpha|^\sigma < \kappa)$ 
is in the proof of $(*)_3$.  So
omit $(*)_3$, and this requires us just to improve the proof of $(*)_4$.  Let
$p_j \in Q$ for $j < \kappa^+$ and let $v_j = \{\zeta < \lambda^*:A_\zeta \cap
{\Cal U}^{p_j}$ has cardinality $\ge \sigma\} \cup v^{p_j}$, so 
clearly $|v_j| \le \kappa$ and $v^{p_j} \subseteq v_j$.

For some stationary $S \subseteq \{\delta < \kappa^+:\text{cf}(\delta) =
\kappa\}$, the conditions $p_j$ for $j \in S$ are pairwise isomorphic and
$j \in S$ implies $v^{p_j} \cap (\dbcu_{i <j} v_i) = v^\otimes$ and
$u^{p_j} \cap (\dbcu_{i<j} (u^{p_i} \cup \cup \dbcu_{\zeta \in v_i} A_\zeta))
= u^\otimes$.  Also \wilog \, for $j_1,j_2 \in S$ the isomorphism $(f,g)$
from $p_{j_1}$ to $p_{j_2}$ satisfies $f \restriction u^\otimes =
\text{ id}_{u^\otimes},g \restriction v^\otimes = \text{ id}_{v^\otimes}$.  
Now for $i,j \in S,p_i,p_j$ are compatible by $(*)_5$ in the proof
of \scite{t.2}.

In the proof of $(*)_5$ and $(*)_6$ (hence $(*)_7$), clause $(E)^-$ gives 
us less but the change in the
definition of $P$ (weakening $(\delta)$ to $(\delta)^-$) demands less and 
they fit.

Lastly, for proving ``$\underset\tilde {}\to X$ is Hausdorff", clause
$(\delta)^-$ is weaker but as $Y^*$ is Hausdorff (and the choice of
$\langle b_i:i < \theta^* \rangle$) there is no problem.
\hfill$\square_{\scite{t.4}}$
\enddemo
\bn
\ub{\stag{t.5} Comment}  1) We could make in \scite{t.2} only some of the
changes from \scite{t.4}, e.g. allow $(A)^-$ and $(E)$. \nl
2) In \scite{t.2}(1) can we make the space regular $(T_3)$?

In view of \scite{t.2}(2) this may be not so interesting, still let
$R_0 \subseteq \{(i,j):b_i \cap b_j = \emptyset\}$ (so to include generalizing
in \scite{t.4} we choose $R_0 \subseteq \{(2i,2i+j):i < \theta^*\}$) and 
$R_1 \subseteq 
\{(i,j):b_i \cup b_j = Y^*\},R_2 \subseteq \{(i,j):b_i \subseteq b_j\}$. \nl
We need: for $i_0 < \theta^*,j < \theta$ such that $j \in b_i$ there are
$i_1,i_2 < \theta^*$ such that $j \in b_{i_1},b_{i_1} \subseteq b_{i_0},
b_{i_0} \cup b_{i_2} = Y^*,b_{i_1} \cap b_{i_2} = \emptyset$ and moreover
$(i_0,i_1) \in R_2,(i_0,i_2) \in R_1,(i_1,i_2) \in R_0$. \nl
Then we should change the definition of $P$, clause $(\delta)$ to
\mr
\item "{$(\delta)^-$}"  $(a) \quad w_{\zeta,i} \subseteq u_*$
\sn
\item "{${{}}$}"  $(b) \quad (i,j) \in R_0 
\Rightarrow w_{\zeta,i} \cap w_{\zeta,j} = \emptyset$;
\sn
\item "{${{}}$}" $(c) \quad (i,j) \in R_1 
\Rightarrow w_{\zeta,i} \cup w_{\zeta,j} = u_*$
\sn
\item "{${{}}$}" $(d) \quad (i,j) \in R_2 
\Rightarrow w_{\zeta,i} \subseteq w_{\zeta,j}$.
\ermn
As $\bar b$ can be with repetition \wilog \, $\langle R_0,R_1,R_2 \rangle$
have a tree structure.  That is \wilog \, there is a partial function
$g^*:\theta^* \rightarrow \theta^*$ such that $g^*(2i) = g^*(2i+1) < i$
(so $2i \in \text{ Dom}(g^*) \leftrightarrow 2i+1 \in \text{ Dom}(g^*)$)
and $R_0 = \{(2i,2i+1):2i \in \text{ Dom}(g^*) \text{ and } i=0
\text{ mod } 3\},R_1 = \{(g^*(2i),2i+1):2i \in \text{Dom}(g^*) \text{ and }
i=1 \text{ mod } 3\},R_2 = \{(g^*(2i),2i):2i \in \text{ Dom}(g^*) \text{ and }
i=2 \text{ mod } 3\}$.
\mn
So we need to have the following property of $Y^*$ (this will be used in 
the proof of $(*)_5,(*)_6$ hence $(*)_7$ dealing with 
$w^{\varepsilon,\ell}_{\zeta,i}$ by induction on $i < \theta^*$)
\mr
\item "{$\boxtimes_{Y^*}$}"  if $k \in \{0,1,2\},{\Cal U}_k$ is open 
for $Y^*,Z \subseteq Y,|Z| < \sigma$ and $Z \cap {\Cal U} = Z \cap$ 
(closure $({\Cal U}))$ and $Z'_\ell \subseteq Z$ for $\ell=0,1,2$ satisfy 
$Z'_1 \subseteq Z'_0,Z'_0 \cup Z'_2 = Z,
Z'_1 \cap Z'_2 = \emptyset$ and $Z'_k = {\Cal U}_k \cap Z$, \ub{then} we can
find open subsets ${\Cal U}_\ell$ of $Y^*$ for $\ell \in \{0,1,2\} \backslash
\{k\}$ such that ${\Cal U}_0
\subseteq {\Cal U}_1,{\Cal U}_0 \cup {\Cal U}_2 = Y^*$ and ${\Cal U}_1 \cap
{\Cal U}_2 = \emptyset$ and $Z'_\ell = {\Cal U}_\ell \cap Z$ for $\ell =1,2$.
\ermn
If $\sigma = \aleph_0$ this requirement on $Y^*$ follows from $Y^*$ being
Hausdorff.  Also clause $(E)$ of \scite{t.2} implies $\boxtimes_{Y^*}$. \nl
3) In \scite{t.3}, as indicated in the proof, we can replace in the assumption
(b) + (e), i.e. ``$\Pi{\frak a}/J$ is $(\lambda^*)^+$-directed, $\lambda^* < 
2^\lambda$" is regular by:
\mr
\item "{$(*)$}"  there is $\bar f = \langle f_\alpha:\alpha < \lambda^*
\rangle,f_\alpha \in \pi {\frak a}$ such that for every $Z \in [\lambda^*]
^{< \kappa}$ we can find $\langle a_\alpha:\alpha \in Z \rangle,a_\alpha \in
J$ such that $\alpha \ne \beta \in Z \and \varepsilon \in \theta \backslash
a_\alpha \backslash a_\beta \Rightarrow f_\alpha(\varepsilon) \ne f_\beta
(\varepsilon)$.
\ermn
4) By the proof of \scite{t.3}, if $({\frak a},J = [{\frak a}]^{< \sigma}$
and $\bar f$ are as in $(*)$ of \scite{t.5}(3)) then
\mr
\item "{$(*)'$}"  there is $\bar f' = \langle f'_\alpha:\alpha < \lambda^*
\rangle,f'_\alpha \in \Pi {\frak a}$ such that for every $Z \in [\lambda^*]
^{< \kappa}$ we can find $\langle a_\alpha:\alpha \in Z \rangle,a_\alpha \in
J$ and well ordering $<^*$ of $Z$ such that $\alpha <^* \beta \in Z \and 
\varepsilon \in \theta \backslash a_\beta \Rightarrow f'_\alpha(\varepsilon)
\ne f'_\beta(\varepsilon)$ (in fact $\bar f' = \bar f$).
\ermn
5) See more \scite{3.10}: for more colours.
\newpage

\head{\S2 Consistency from supercompact} \endhead  \resetall 
\bigskip

In the first section we got consistency results concerning Arhangelskii's
problem but using pcf statement of unclear status (they come from
\scite{t.3}); this is very helpful toward finding the consistency strength,
and unavoidable if e.g. we like CH to fail (see \S3), but it does not 
give a well grounded consistency result.  Here relying on Theorem 
\scite{t.2} of the first section, we get 
consistency results using supercompact cardinals.  First
we give a sufficient condition for clause (C) of Theorem \scite{t.2} which
is reasonable under instances of G.C.H.  We then (\scite{g.6}) quote a
theorem of Hajnal Juhasz Shelah \cite{HJSh:249}, 
\cite{Sh:F267} (for $\sigma = \aleph_0,\sigma
> \aleph_0$, respectively) and from it (in claim \scite{g.6A}), in the
natural cases, prove that the assumptions of \scite{t.2} hold deducing
(in \scite{g.7}) the consistency of CH + there is a $T_3$-space $X$ with
clopen basis with $\aleph_{\omega +1}$ point such that $X \rightarrow$
(Cantor set)$^1_{\aleph_0}$ starting with a supercompact cardinal.  This
gives a (consistent) negative answer to Arhangelskii's problem.  We can even
make it compact.

\bn
\ub{\stag{g.5} Observation}:  If clauses \footnote{(actually from $(B)_1$, only
``$(B)^-_1{\Cal A} \subseteq [\lambda]^\theta$" is used; as we do not change
${\Cal A}$ and the cardinals this is O.K.} $(A)(i) + (B)_1$ of Theorem 
\scite{t.2} holds, \ub{then} clause (C) there follows from
\mr
\item "{$(C)^+$}"  if $\langle Y_i:i < \kappa^+ \rangle$ is a partition of
$\lambda$ then for some $A \in {\Cal A}$ and $i < \kappa^+$\nl
we have $A \subseteq Y_i$.
\endroster
\bigskip

\demo{Proof}  Let $F:\lambda \rightarrow [\lambda]^{\le \kappa}$ be given.
Choose by induction on $\zeta \le \lambda$ a set $U_\zeta \subseteq \lambda$
and $g_\zeta:U_\zeta \rightarrow \kappa^+$, both increasingly continuous
with $\zeta$ such that:
\mr
\item "{$(*)(i)$}"  if $\alpha \in U_\zeta$ then $F(\alpha) \subseteq
U_\zeta$ and
\sn
\item "{$(ii)$}"  if $\alpha \in U_\zeta$ then $F(\alpha) \backslash
\{\alpha\} \subseteq \{\beta \in U_\zeta:g_\zeta(\beta) \ne g_\zeta
(\alpha)\}$.
\ermn
For $\zeta = 0$ let $U_\zeta = \emptyset = g_\zeta$, for $\zeta$ limit
take unions.  If $U_\zeta = \lambda,U_{\zeta + 1} = U_\zeta,g_{\zeta +1} =
g_\zeta$, otherwise let $\alpha_\zeta = \text{ Min}\{\lambda \backslash
U_\zeta\}$ and let $W_\zeta \in [\lambda]^{\le \kappa}$ be such that
$\alpha_\zeta \in W_\zeta$ and 
$(\forall \alpha \in W_\zeta)[F(\alpha) \subseteq
W_\zeta]$.  Let $\varepsilon_\zeta = 
\sup\{g_\zeta(\beta):\beta \in U_\zeta \cap
W_\zeta\}$ so $\varepsilon_\zeta < \kappa^+$ and let 
$U_{\zeta +1} = U_\zeta \cup W_\zeta,g_{\zeta +1}$ extends $g_\zeta$
such that $g_{\zeta +1} \restriction (W_\zeta \backslash U_\zeta)$ is one
to one with range $[\varepsilon,\varepsilon + \kappa)$. \nl
Now applying $(C)^+$ to the partition which $\dbcu_\zeta g_\zeta$ defines,
we get some $A \in {\Cal A}$ on which $\dbcu_\zeta g_\zeta$ is constant
so by $(*)(ii)$ we are done.  \hfill$\square_{\scite{g.5}}$
\enddemo
\bn
By \cite{HJSh:249}, \cite{Sh:F276}
\proclaim{\stag{g.6} Claim}  Assume $V \models GCH$ (for simplicity) and
$\sigma < \chi < \chi_0^{< \chi} \le \kappa < \mu < \mu^+ = \lambda$ and
$\sigma,\chi,\chi_0,\kappa,\lambda$ are regular, {\rm cf\/}$(\mu) = 
\sigma,\chi$ is
a supercompact cardinal $> \sigma$ (or just $\lambda$-supercompact), e.g.
$\mu = \chi^{+ \sigma}_0$. \nl
\ub{Then} for some forcing notion, $\sigma$-complete of cardinality
$\chi_0$, in $V^P,2^\sigma = \sigma^+ = \theta,
2^{\sigma^+} = \chi_0 = \sigma^{++}$
(and GCH holds) and some $\langle B_\delta:\delta \in S \rangle$ satisfies:
\mr
\item "{$(*)$}"  $S \subseteq \{\delta < \lambda:\text{cf}(\delta) =
\sigma^+\}$ is stationary, \nl
$B_\delta \subseteq \delta$, {\rm otp\/}$(B_\delta) = \theta^+$ and 
$\delta_1 \ne \delta_2 \in S \Rightarrow |B_{\delta_1} \cap B_{\delta_2}| 
< \sigma$ and \nl
$\theta^* < \kappa = \kappa^{< \kappa} < \lambda$, (hence
$\diamondsuit_S$, so if $\mu = \chi^{+ \omega}_0$ then $\lambda = 
\aleph_{\omega + 1}$).
\endroster
\endproclaim
\bn
What we need is getting in such model, condition $(C)^+$ of \scite{g.5} which
also is from \cite{Sh:F276} but for completeness we shall prove it.
\proclaim{\stag{g.6A} Claim}  Assume
\mr
\item "{$(a)$}"   $\langle B_\delta:\delta \in S \rangle,
\sigma,\kappa,\mu,\lambda$ are as in the conclusion $(*)$ of the 
previous claim and
\sn
\item "{$(b)$}"   $S$ reflects in no ordinal of cofinality $\le \kappa$
(holds automatically if $\kappa < \sigma^{+\sigma}$, see \cite{Sh:108},
\cite{Sh:88a}), but see \scite{g.8}, \scite{g.9a}.
\ermn
\ub{Then} withought loss of generality
$\sigma,\theta =: \sigma^+,\lambda,{\Cal A} = \{B_\delta:\delta \in S\}$
satisfies the requirements in Theorem \scite{t.2} (in $V^P$).
\endproclaim
\bigskip

\demo{Proof}  Without loss of generality ``$\delta \in S \Rightarrow
\mu^\omega$ divides $\delta$", and as we are assuming GCH and $\delta \in
S \Rightarrow \text{ cf}(\delta) = \sigma^+ \ne \sigma = \text{ cf}(\mu)$ 
we have $\diamondsuit_S$ (\cite{Sh:108}).  
So let $\langle f_\delta:\delta \in S \rangle$
be such that $f_\delta:\delta \rightarrow [\delta]^\kappa$ satisfy
$(\forall f:\lambda \rightarrow [\lambda]^\kappa)(\exists^{\text{stat}}
\delta \in S)(f_\delta = f \restriction \delta)$.  For each $\delta \in S$,
let $B_\delta = \{\alpha_{\delta,\varepsilon}:\varepsilon < \sigma^+\}$
increasing with $\varepsilon$ and let $g_\delta:\kappa^{++} \rightarrow
[\kappa^{++}]^{\le \kappa}$ be defined by

$$
\align
g_\delta(\beta) = \bigl\{ \gamma < \kappa^{++}:&\text{ for some }
\varepsilon_1,\varepsilon_2 < \sigma^+ \text{ we have} \\
  &\alpha_{\delta,\varepsilon_1} \times \kappa^{++} + \gamma \in
f_\delta(\alpha_{\delta,\varepsilon_2} \times \kappa^{++} + \beta) \bigr\}.
\endalign
$$
\mn
So by a variant of the $\Delta$-system lemma (or use $(2^\kappa)^+$ instead
$\kappa^{++}$ if we avoid GCH) there is 
$Z_\delta \in [Z_\delta]^{\kappa^{++}}$ such that
$\gamma_1 \ne \gamma_2 \in Z_\delta \Rightarrow \gamma_1 
\notin g_\delta(\gamma_2)$.
Let $\gamma_{\delta,\varepsilon} \in Z_\delta$ be strictly increasing with
$\varepsilon < \sigma^+$ and let $B'_\delta = \{\alpha_{\delta,\varepsilon}
\times \kappa^{++} + \gamma_{\delta,\varepsilon}:\varepsilon < \sigma^+\}$.
So clauses $(A),(B)_1$ are immediate.  Now clearly $(C)^+$ of \scite{g.5}
holds hence (C) and $(B)_2$ of \scite{t.2} follow from the assumption on $S$
(see \cite{Sh:108}).

Now ${\Cal A} = \{B'_\delta:\delta \in S\}$ are as required in Theorem 
\scite{t.2}. \hfill$\square_{\scite{g.6A}}$
\enddemo
\bn
\ub{\stag{g.7} Conclusion}:  If CON($\exists$ supercompact), \ub{then}
CON(CH + there is a $T_3$-topological space $X$ with clopen basis, even
compact, with $\aleph_{\omega +1}$ members, $\aleph_{\omega+1}$ nodes such
that if we divide $X$ to countably many parts, at least one contains a 
closed copy of the Cantor set).
\bigskip

\demo{Proof}  By \scite{g.6A} + \scite{t.2}.
\hfill$\square_{\scite{g.7}}$
\enddemo
\bn
Instead of using \cite{HJSh:249}, \cite{Sh:F276} we can directly 
use \cite{Sh:108}, recall that there even G.C.H. holds.
\bn
\centerline {$* \qquad * \qquad *$}
\bn
Lastly, we
start to resolve the connection between the various statements around.  Now
\cite{HJSh:249} continue and strengthen \cite{Sh:108}, \cite{Sh:88a}.  We show
that by a ``small nice forcing" (not involving extra large cardinals
assumption) we can get the result of \cite{HJSh:249} used above from the one
in \cite{Sh:108}, \cite{Sh:88a}.  (See also \cite[\S5]{Sh:652} 
on the semi-additive colouring involved). 
\proclaim{\stag{g.8} Claim}  Assume
\mr
\item "{$(a)$}"  cf$(\mu) = \kappa < \mu,(\forall \alpha < \mu)
(|\alpha|^\kappa < \mu)$
\sn
\item "{$(b)$}"  $S \subseteq \{\delta < \mu^+:\text{cf}(\delta) = 
\kappa^+\}$ is stationary, $S \notin I[\lambda]$ and
\sn
\item "{$(c)$}"  $2^{\kappa^+} \le \mu$ and $\kappa = \kappa^{< \kappa}$.
\ermn
\ub{Then} for some forcing notion $Q$ we have:
\mr
\item "{$(a)$}"  $Q$ is $(< \kappa)$-complete, $|Q| = \kappa^+$ and $Q$ is
$\kappa^+$-c.c.
\sn
\item "{$(b)$}"  in $V^Q$, for some stationary $S' \subseteq S$ we have
$\langle A_\delta:\delta \in S' \rangle,A_\delta$ an unbounded subset of
$\delta$ of order type $\kappa^+$ and $\delta_1 \ne \delta_2 \in S'
\Rightarrow |A_{\delta_1} \cap A_{\delta_2}| < \kappa$.
\endroster
\endproclaim
\bigskip

\demo{Proof}  Let $\mu = \dsize \sum_{i < \kappa} \lambda_i$ where
$\lambda_i < \mu$ is 
increasing continuous with $i,\lambda_0 > \kappa$.  Choose
$\bar A = \langle A_\alpha:\alpha \in S \rangle$, with $A_\alpha = 
\{\gamma_{\alpha,\varepsilon}:\varepsilon < \kappa^+\}$ any unbounded subset
of $\alpha$ of order type $\kappa^+$ and $\gamma_{\alpha,\varepsilon}$ 
increasing with $\varepsilon$.
\sn
We can find $\bar a^\alpha = \langle a^\alpha_i:i < \kappa \rangle$ for
$\alpha < \mu^+$ such that
\mr
\item "{$(*)_1$}"  $\alpha = \dbcu_{i < \kappa} a^\alpha_i,a^\alpha_i$ is
increasing continuous in $i,|a^\alpha_i| \le \lambda_i$
\sn
\item "{$(*)_2$}"  if $\alpha \in a^\beta_i$ then $a^\alpha_i \subseteq
a^\beta_i$. \nl
Without loss of generality
\sn
\item "{$(*)_3$}"  $A_\alpha \subseteq a^\alpha_0$.
\ermn
Let $\bold c:[\mu^+]^2 \rightarrow \kappa$ be $\bold c\{\alpha,\beta\} =
\text{ Min}\{i:\alpha \in a^\beta_i\}$ for $\alpha < \beta < \lambda^+$ so
\mr
\item "{$\boxtimes$}"  $\alpha < \beta < \gamma \Rightarrow 
\bold c\{\alpha,\gamma\} \le \text{ Max}\{\bold c\{\alpha,\beta\},
\bold c\{\beta,\gamma\}\}$.
\ermn
For $\alpha \in S$ let $c_\alpha:[\kappa^+]^2 \rightarrow \kappa$ be defined
by:

for $\varepsilon < \zeta < \kappa^+$ we let

$$
c_\alpha\{\varepsilon,\zeta\} = \bold c\{\gamma_{\alpha,\varepsilon},
\gamma_{\alpha,\zeta}\}.
$$
\mn
Let ${\Cal C} = \{c_\alpha:\alpha \in S\}$ so $c_\alpha \in
{}^{([\kappa^+]^2)}\kappa$, so $|{\Cal C}| \le 2^{\kappa^+}$.  Let for
$c \in {\Cal C},S_c = \{\alpha \in S:c_\alpha = c\}$, so $\langle S_c:c \in
{\Cal C} \rangle$ is a partition of $S$ to $\le 2^{\kappa^+} < \mu^+$ sets
hence necessarily for some $c \in {\Cal C}$ we have
\mr
\item "{$(*)_4$}"  $S_c \notin I[\lambda]$.
\ermn
We fix $c$.  We define a forcing notion $Q$:
\mr
\item "{$(A)$}"  $p \in Q$ iff $p = (u^p,\xi^p)$ where $u^p \in
[\kappa^+]^{< \kappa}$ and $\xi^p < \kappa$ and Rang$(c \restriction [u^p]^2)
\subseteq \xi^p$
\sn
\item "{$(B)$}"  $Q \models p \le q$ \ub{iff}: $(p,q \in Q$ and)
{\roster
\itemitem{ $(i)$ }  $u^p \subseteq u^q$
\sn
\itemitem{ $(ii)$ }  $\xi^p \le \xi^q$
\sn
\itemitem{ $(iii)$ }  for every 
$\beta \in u^p$ and $\alpha \in (u^q \backslash u^p)
\cap \beta$ we have $\bold c\{\alpha,\beta\} \ge \xi^p$
\endroster}
\item "{$(*)_5$}"  $(a) \quad Q$ is a partial order 
\sn
\item "{${{}}$}"  $(b) \quad Q' = 
\{p \in Q:u^p$ has a maximal element$\}$ is a dense subset of $Q$ \nl
[why?  check Clause (a), as for clause (b), for any 
$p \in Q$ choose $j \in(\sup(u^p)+1,\kappa^+)$ and define
$q = (u^q,\xi^q)$ by $u^q = u^p \cup \{j\}$ and $\xi^q = \xi^p$]
\sn
\item "{$(*)_6$}"  $Q$ satisfies the $\kappa^+$-c.c. \nl
[why?  assume toward contradiction that 
$\langle p_i:i < \kappa^+ \rangle$ are pairwise incompatible.
Without loss of generality $p_i \in Q'$.  Without loss of generality
$\langle u^{p_i}:i < \kappa^+ \rangle$ is a $\Delta$-system with heart $u^*$.
So without loss of generality $\xi^{p_i} = \xi^*$.  
So $C = \{\delta < \kappa^+:u^{p_\delta} \backslash u^*$ is 
disjoint to $\delta$
and $(\forall j < \delta)(u^{p_j} \subseteq \delta)\}$ is a club of
$\kappa^+$.  Let for $\delta \in C,\varepsilon_\delta = \text{ Min}
(u^{p_\delta} \backslash \delta)$ and $\zeta_\delta = \text{ max}
(u^{p_\delta})$ so $\delta \le \varepsilon_\delta \le \zeta_\delta$.  Now
assume $\alpha < \beta$ are from $C$, and $p_\alpha,p_\beta$ is incompatible.
Why is $q = (u^{p_\alpha} \cup u^{p_\beta},\zeta)$ not a common upper bound
where we let $\zeta = \sup(\{\xi^*\} \cup \text{ Rang}(\bold c \restriction 
[u^{p_\alpha}\cup u^{p_\beta}]^2))+1$?  
As $q \in Q$ and as $u^{p_\alpha} \cap \alpha = 
u^{p_\beta} \cap \beta,u^{p_\alpha} \subseteq
\beta$ and $\xi^* = \xi^{p_\alpha} = \xi^{p_\beta}$ so $p_\alpha \le q$, hence
necessarily $\neg p_\beta \le q$ so clause (iii) of (B) fails, i.e. for some
$\gamma_1 \in u^{p_\alpha} \backslash \alpha$ and $\gamma_2 \in u^{p_\beta}
\backslash \beta$ we have $\bold c\{\gamma_1,\gamma_2\} < \xi^{p_\beta} = 
\xi^*$.  
But $\varepsilon_\alpha \le \gamma_1$ and
$\varepsilon_\alpha < \gamma_1 \Rightarrow \bold c
\{\varepsilon_\alpha,\gamma_1\}
< \xi^{p_\alpha} = \xi^*$ and $\gamma_2 \le \zeta_\beta$ and $\gamma_2 <
\zeta_\beta \Rightarrow \bold c
\{\gamma_2,\zeta_\beta\} < \xi^{p_\beta} = \xi^*$.
Hence by $\boxtimes$ necessarily $\bold c
\{\varepsilon_\alpha,\zeta_\beta\} < \xi^*$. \nl
So for $\delta \in S_c,\langle \gamma_{\delta,\varepsilon_i}:i \in C \rangle$
is strictly increasing hence with limit $\delta$ and for each $i \in C,
\gamma_{\delta,\zeta_i}$ is above $\{\gamma_{\delta,\varepsilon_j}:j < i\}$
but $< \delta$ and

$$
j < i \Rightarrow \bold c\{\gamma_{\delta,\varepsilon_j},
\gamma_{\delta,\zeta_i}\} < \xi^* \Rightarrow \gamma_{\delta,
\varepsilon_j} \in a^{\gamma_{\delta,\zeta_i}}_{\xi^*}.
$$
\mn
By \cite{Sh:108} it follows that $S \in I[\lambda]$
(or directly, for every $\gamma < \lambda,|\{\langle \gamma_{\delta,
\varepsilon_j}:j \in C \cap i^* \rangle:\delta \in S,i^* \in C_\zeta,
\gamma_{\delta,\zeta^*_i} = \gamma\}| < \lambda$ as for each $i < \kappa^+$
(and $\gamma$) we have $\le |a^\gamma_{\xi^*}|^{|i^*|} \le
(\lambda_{\xi^*})^{|i|} \le \mu$ possibilities); contradiction.  So $Q$
satisfies the $\kappa^+$-c.c.]
\ermn
Now clearly for every $i < \kappa^+$ there is $p_i \in Q'$ such that $i \in
u^{p_i}$, hence (by $(*)_6$), for some $i(*) < \kappa^+$ we have
$p_{i(*)} \Vdash_Q ``{\underset\tilde {}\to W_1} = \{i:p_i \in G 
\text{ and\, cf}(i) = \kappa\}$ is stationary in $\kappa^+$".  Let $p_{i(*)}
\in G \subseteq Q,G$ generic over $V$ and $W_1 = {\underset\tilde {}\to W_1}
[G]$.  Let $C = \{ \delta < \kappa^+:(\forall i < \delta)$ sup $(u^{p_i}) < 
\delta\}$, it is a club of $\kappa^+$.  Let $W_2 = C \cap W_1$ and for $i \in
S_2$ let $\varepsilon_i = \text{ Min}(u^{p_i} \backslash i),\zeta_i =
\text{ max}(u^{p_i})$.  Now
\mr
\item "{$(*)_7$}"  if $i \in W_2$ and $\xi < \kappa$, \ub{then}
$\{j \in W_1 \cap i:\bold c
(\varepsilon_j,\varepsilon_i) < \xi\}$ has cardinality $< \kappa$.
\sn
[Why?  By density argument for some $q \in G$ we have $p_i \le q$ and
$\xi^q > \xi$.  Now if $j \in S_1 \cap i \backslash u^q$ then $p_j \in G$
hence for some $q^+ \in G \subseteq Q$ we have $q \le q^+ \and p_j \le q^+$,
so $\varepsilon_j \in u^{q^+} \cap \varepsilon_i$ and as $q \le q^+$ by the
definition of $\le^Q$, necessarily $c(\varepsilon_i,\varepsilon_j) \ge
\xi^q > \xi$, as asserted.]
\ermn
Now define for $\delta \in S_c,A'_\delta = \{\gamma_{\delta,\varepsilon}:
\varepsilon \in W_2\}$.  So $A'_\delta$ is an unbounded subset of $\delta$
of order type $\kappa^+$.
\mr
\item "{$(*)_8$}"  if $\delta_1 \ne \delta_2$ are from $S_c$ then
$A'_{\delta_1} \cap A'_{\delta_2}$ has cardinality $< \kappa$.
\sn
[Why?  Without loss of generality, let $\delta_1 < \delta_2$, let
$\varepsilon(*) \in S_2$ be such that $\delta_1 < \gamma_{\delta_2,
\varepsilon(*)}$.  Assume toward contradiction that $A = A'_{\delta_1} \cap
A'_{\delta_2}$ has cardinality $\ge \kappa$.  Recall (by $(*)_3$) that
$\beta \in A \Rightarrow \bold c\{\beta,\delta_1\} = 0$, letting $\xi^* =
\bold c\{\delta_1,\gamma_{\delta_2,\varepsilon(*)}\}$ we get by $\boxtimes$
that $\beta \in A \Rightarrow \bold c\{\beta,
\gamma_{\delta_1,\varepsilon(*)}\} \le \max\{ \bold c\{\beta,\delta_1\},
\bold c\{\delta_1,\gamma_{\delta_1,\varepsilon(*)}\} = \max\{0,\xi^*\} =
\xi^*$. \nl
So $A^- = \{\varepsilon:\gamma_{\delta_2,\varepsilon} \in A\}$ has cardinality
$\kappa$ and $\varepsilon \in A^- \Rightarrow c\{\varepsilon,\varepsilon(*)\}
\le \xi^*$, contradicting $(*)_7$.]  
\ermn
So we are done.  \hfill$\square_{\scite{g.8}}$
\enddemo
\bigskip

\proclaim{\stag{g.9a} Claim}  Assume
\mr
\item "{$(A)(i)$}"  $\lambda > \kappa > \theta > \sigma \ge \aleph_0$
and $\kappa = \kappa^{< \kappa}$
\sn
\item "{$(B)_1$}"  ${\Cal A} \subseteq [\lambda]^\theta$ and
$A_1 \ne A_2 \in {\Cal A} \Rightarrow |A_1 \cap A_2| < \sigma$.
\ermn
\ub{Then} for some forcing notion $Q$
\mr
\item "{$(a)$}"  $Q$ is a strategically $< \kappa$-complete forcing notion
(hence add no new sequence of length $< \kappa$)
\sn
\item "{$(b)$}"  $Q$ is $\kappa^+$-c.c. forcing notion of cardinality
$\lambda^{< \kappa}$
\sn
\item "{$(c)$}"  in $V^Q$, clauses $(A)(i),(B)_1$ above still hold and $(B)_2$
from \scite{t.2}, i.e.
\sn
\item "{$(B)_2$}"  ${\Cal A}$ is $\kappa$-free
\sn
\item "{$(d)$}"  if $\lambda,\kappa,{\Cal A}$ satisfies clause $(C)$ of
\scite{t.2} in $V$, then this still holds in $V^Q$.
\endroster
\endproclaim
\bigskip

\demo{Proof}  Let ${\Cal A} = \{A_\zeta:\zeta < \lambda^*\}$ with no
repetitions.  \nl
Let $Q$ be the set of $p = (v,v_*) = (v^p,v^p_*)$ such that:
\mr
\item "{$(a)$}"  $v_* \subseteq v \in [\lambda^*]^{< \kappa}$
\sn
\item "{$(b)$}"  there is a list $\langle \zeta(\varepsilon):\varepsilon <
\varepsilon^* \rangle$ of $v_*$ such that for every $\varepsilon <
\varepsilon^*$ we have $A_{\zeta(\varepsilon)} \cap \dbcu_{\xi < \varepsilon}
A_{\zeta(\xi)}$ has cardinality $< \sigma$; we call $\langle \zeta
(\varepsilon):\varepsilon < \varepsilon^* \rangle$ a witness, also the well
ordering on $u^p_*$ it induces is called a witness.
\ermn
The order is defined by

$$
\align
p \le q \text{ iff } &(\alpha) \qquad v^p_* \subseteq v^q_* \text{ and} \\
  &(\beta) \qquad v^p \backslash v^p_* \subseteq v^q \backslash v^q_* \\
  &(\gamma) \qquad \text{ every } \bar \zeta \text{ witnessing } p \in Q
\text{ can be end-extended to } \bar \zeta' \\
  &\qquad \qquad \qquad \text{witnessing } q \in Q.
\endalign
$$
\mn
Define a $P$-name $\underset\tilde {}\to Y = \cup\{v^p_*:p \in
{\underset\tilde {}\to G_Q}\},{\Cal A}' = \{A_\zeta:\zeta \in
\underset\tilde {}\to Y\}$.  Now
\mr
\item "{$(*)_1$}"  $Q$ is a partial order
\sn
\item "{$(*)_2$}"  $|Q| = (\lambda^*)^{< \kappa} \le 
(\lambda^{< \theta})^{< \kappa} = \lambda^{< \kappa}$
\sn
\item "{$(*)_3$}"  any increasingly continuous sequence of members of $Q$ of
length $< \kappa$ has a least upper bound. \nl
Hence
\sn
\item "{$(*)_4$}"  $Q$ is strategically $(< \kappa)$-complete. \nl
\sn
For $p \in Q$ let $u^p = \cup\{A_\zeta:\zeta \in v^p\}$
\sn
\item "{$(*)_5$}"  for $p \in Q$ we have $u^p \in [\lambda]^{< \kappa}$
and $p \le q \Rightarrow u^p \subseteq u^q$.
\ermn
Let $Q' = \{p \in Q:\text{if } \zeta < \lambda^* \text{ and } |A_\zeta \cap
u^p| \ge \sigma$ then $\zeta \in v^p\}$. 
For $p \in Q$ let $v^p_\otimes = \{\zeta < \lambda^*:|A_\zeta \cap u^p|
\ge \sigma\}$, so:
\mr
\item "{$(*)_6$}"  $(a) \quad v^p \subseteq 
v^p_\otimes$ and $p \in Q \Rightarrow |v^p_\otimes| \le \kappa$ and
\sn
\item "{${{}}$}"  $(b) \quad$ if $(\forall \alpha < \kappa)
[|\alpha|^\sigma < \kappa]$ then $p \in Q 
\Rightarrow |v^p_\otimes| < \kappa$, and
\sn
\item "{${{}}$}"  $(c) \quad p \in Q \Rightarrow 
p \in Q' \equiv v^p_\otimes = v^p$ and
\sn
\item "{${{}}$}"   $(d) \quad Q'$ is a dense subset of 
$Q$ if $(\forall \alpha < \kappa)[|\alpha|^\sigma < \kappa]$
\sn
[Why?  E.g. for clause (d), let $p \in Q$ we choose 
by induction on $\varepsilon \le \sigma^+
(< \kappa)$ a condition $p_\varepsilon$ such that: $p_0 = p,
v^{p_\varepsilon}_* = v^p_*,p_\varepsilon$ is increasingly continuous with
$\varepsilon$ and $v^{p_{\varepsilon +1}} = \{\zeta < \lambda^*:\zeta \in
v^{p_\varepsilon}$ or just $|A_\zeta \cap u^p| \ge \sigma\}$.  There are
no problems and $p_{\sigma^+}$ is as required as 
$|A_\zeta \cup u^{p_{\sigma^+}}| \ge 
\sigma \Rightarrow \text{ for some } \varepsilon < \sigma^+,|A_\zeta \cap
u^{p_\varepsilon}| \ge \sigma \Rightarrow \text{ for some } \varepsilon <
\sigma^+,\zeta \in v^{p_{\varepsilon +1}} \subseteq v^{p_{\sigma^+}}$.]
\sn
\item "{$(*)_7$}"  if $p \in Q',\zeta \in \lambda^* \backslash v^p$ or just
$p \in Q,\zeta \in \lambda^* \backslash v^p_\otimes$ then $p' = (v^p \cup 
\{\zeta\},v^p_* \cup \{\zeta\})$ and $p'' = (v^p \cup \{\zeta\},v^p_*)$
are in $Q$ (even $p \in Q' \Rightarrow p' \in Q'$) and are $\ge p$. 
\ermn
We say $p_0,p_1 \in Q$ are isomorphic if otp$(v^{p_0}) = 
\text{ otp}(v^{p_1})$, otp$(u^{p_0}) = \text{ otp}(u^{p_1})$,
and OP$_{v^{p_1},v^{p_0}}$ maps $v^{p_0}_*$ onto $v^{p_1}_*,
OP_{u^{p_1},u^{p_0}}$ maps $u^{p_0}$ onto $u^{p_1}$ and for $\zeta \in
v^{p_0},\alpha \in u^{p_0}$ we have $\alpha \in A_\zeta \Leftrightarrow
\text{ OP}_{u^{p_1},u^{p_0}}(\alpha) \in 
A_{\text{OP}_{v^{p_1},v^{p_0}}(\zeta)}$
\mr
\item "{$(*)_8$}"  $Q$ satisfies the $\kappa^+$-c.c. \nl
[Why?  Let $p_\alpha \in Q$ for $\alpha < \kappa^+$.
Let $v_\alpha = \dbcu_{\beta < \alpha} v^{p_\beta}_\otimes$ and $u_\alpha =
\cup\{A_\zeta:\zeta \in v_\alpha\}$ so $u^{p_\beta} \subseteq u_\alpha$ for
$\beta < \alpha$.  As $v^{p_\alpha} \in [\lambda^*]
^{< \kappa}$, we can find stationary $\{S < \kappa^+:\text{cf}(\delta) =
\kappa\}$ and $v$ such that $\alpha \in S \Rightarrow v^{p_\alpha} \cap
v_\alpha = v$.  Similarly \wilog \, $\alpha \in S \Rightarrow u^{p_\alpha}
\cap u_\alpha = u$.  Without loss of generality for $\alpha,\beta \in S$
the conditions $p_\alpha,p_\beta$ are isomorphic, the isomorphism preserving 
$v$ and $u$.  So $v^{p_\alpha}_* \cap v = v_*$ for some 
$v_* \subseteq v$. Let $<^*_\alpha$ be a well ordering of 
$v^{p_\alpha}_*$ which witnesses $p_\alpha \in Q$, so
\wilog \, $<^*_\alpha \restriction v_* = <^*$.  
Let $\alpha < \beta$ be in $S$ and
define $q = (v^{p_\alpha} \cup v^{p_\beta},v^{p_\alpha}_* \cup
v^{p_\alpha}_*$).  Clearly $v^q_* \subseteq v^q \in [\lambda^*]^{< \kappa}$,
also $\zeta \in v^{p_\alpha} \backslash v^{p_\beta}$ or
$\zeta \in v^{p_\beta} \backslash v^{p_\alpha}$ implies $|A_\zeta \cap u) <
\sigma$ (why if not we can find $\zeta_1 \in v^{p_\alpha} \backslash 
v^{p_\beta},\zeta_2 \in v^{p_\beta} \backslash v^{p_\alpha}$ such that
$\zeta \in \{\zeta_1,\zeta_2\}$ and $A_{\zeta_1} \cap u = A_{\zeta_2} \cap
u$, so $|A_{\zeta_\ell} \cap u| \le |A_{\zeta_1} \cap A_{\zeta_2}| < \sigma$
hence $|A_\zeta \cap u| < \sigma$).  Hence $\zeta \in v^{p_\alpha} \backslash
v^{p_\beta} \Rightarrow |A_\zeta \cap u^{p_\beta}| < \sigma$ (otherwise
$A_\zeta \cap u^{p_\beta} \subseteq u_\beta \cap u^{p_\beta} = u$ hence
$|A_\zeta \cap u| \ge \sigma$ and get a contradiction by the previous
statement) and $\zeta \in v^{p_\beta} \backslash v^{p_\alpha} \Rightarrow
|A_\zeta \cap v^{p_\alpha}| < \sigma$ (similar proof).
Now define a two-place relation $<^*$ on $v^q_*$:
$$
\align
\zeta_1 <^* \zeta_2 \text{ \ub{iff} } &\zeta_1 <^*_\alpha \zeta_2
(\text{so } \zeta_1,\zeta_2 \in v^{p_\alpha}) \\
  &\text{or } \zeta_1 \in v^{p_\alpha}_* \and \zeta_2 \in v^{p_\beta}_*
\backslash v^{p_\alpha}_* \\
  &\text{or } \{\zeta_1,\zeta_2\} \subseteq v^{p_\beta}_* \backslash
v^{p_\alpha}_* \and \zeta_1 <^*_\beta \zeta_2
\endalign
$$
Easily $<^*$ is a well order of $v^q_*$ (as $\zeta \in v^{p_\beta}_*
\backslash v^{p_\alpha}_* \Rightarrow |A_\zeta \cap u^{p_\alpha}| < \sigma$,
and it is a witness).   So $q \in Q$.  
Does $p_\alpha \le q$?  Clauses $(\alpha),(\beta)$ are very straight
and for clause $(\gamma)$, as $p_\alpha,p_\beta$ are isomorphic for any given
witness $<^1$, a well ordering of $v^{p_\alpha}_*$, we can find $<^2$, a
witness for $p_\beta$ which is a well ordering of $v^{p_\beta}_*$, and is
conjugate to $<^1$; now use $<^1,<^2$ as we use $<^*_\alpha,<^*_\beta$
above.  So really $p_\alpha \le q$.  Similarly $p_\beta \le q$.]
\sn
\item "{$(*)_9$}"  $\Vdash_Q ``{\underset\tilde {}\to {\Cal A}'} = \{A_\zeta:
\zeta \in \cup\{v^p_*:p \in {\underset\tilde {}\to G_Q}\}\}$ is 
$(< \kappa)$-free". \nl
[Why?  Read the definitions of $Q$ and of being $(< \kappa)$-free, 
remembering that forcing with $Q$ add no new sets of ordinals 
$< \kappa$ as it is strategically $(< \kappa)$-complete.]
\sn
\item "{$(*)_{10}$}"  if $p,q \in Q$ are compatible, \ub{then} they have an
upper bound $r \in Q$ such that $v^r = v^p \cup v^q$
\sn
\item "{$(*)_{11}$}"  if ${\Cal A}$ satisfies clause (C) of \scite{t.2}
then ${\underset\tilde {}\to {\Cal A}'}$ satisfies this in $V^Q$.
\ermn
[Why?  Assume $p^* \in Q,p^* \Vdash_Q ``\underset\tilde {}\to F:\lambda
\rightarrow [\lambda]^{\le \kappa}$ is a counterexample".  Without loss of
generality $p^* \in Q'$.  As $Q$ satisfies the $\kappa^+$-c.c. and as 
increasing the $\underset\tilde {}\to F(\alpha)$ is O.K., \wilog \, each
$\underset\tilde {}\to F(\alpha)$ is from $V$ and
$\underset\tilde {}\to F = F$.  As we can increase each $F(\alpha)$,
\wilog \, $\zeta \in v^{p^*}_\otimes \Rightarrow A_\zeta \subseteq 
\dbca_\alpha
F(\alpha)$.  As $V,{\Cal A}$ satisfies clause (C) there are $\zeta$ and $A
\in [A_\zeta]^\theta$ which is $F$-free, by the prvious sentence $\zeta \notin
v^{p^*}_\otimes$.  Define $q = (v^q,v^q_*),v^q = v^{p^*} \cup \{\zeta\},
v^q_* = v^{p^*}_* \cup \{\zeta\}$.  It is easy to prove $p^* \le q \in Q$,
the point being $|A_\zeta \cap \{A_\xi:\xi \in v^{p^*}_*\}| < \sigma$ which
holds as $\zeta \notin v^{p^*}_\otimes$, and $q$ forces that $A \in
[A_\zeta]^\theta$ is as required concerning $F$.  An alternative (older)
proof is for each $\alpha$ let
$\langle p_{\alpha,i}:i < \kappa \rangle$ be a maximal antichain above
$p^*$ of $Q$ such that $p_{\alpha,i} \Vdash_Q ``\underset\tilde {}\to F
(\alpha) = a_{\alpha,i}"$.
\mn
Define $F:\lambda \rightarrow [\lambda]^{\le \kappa}$ in $V$ by
$F(\alpha) = \dbcu_{i < \kappa} a_{\alpha,i} \cup \dbcu_{i < \kappa}
u^{p_{\alpha,i}}$.  So as in $V,{\Cal A}$ satisfies clause (C), there is
$\zeta \in {\Cal A}$ such that $A_\zeta$ is $F$-free or just some
$A \in [A_\zeta]^\theta$ is $F$-free.  Let $A = 
\{\gamma_\varepsilon :\varepsilon  < \theta\}$.  We can now choose by
induction $\varepsilon < \theta,(p_\varepsilon,j_\varepsilon)$ such that:
\mr
\item "{$(a)$}"  $p_\varepsilon \in Q$ is increasingly continuous
\sn
\item "{$(b)$}"  $j_\varepsilon < \kappa$
\sn
\item "{$(c)$}"  $p_0 = p^*$
\sn
\item "{$(d)$}"  $p_{\varepsilon +1}$ is an upper bound of $p_\varepsilon,
p_{\gamma_\varepsilon,j_\varepsilon}$
\sn
\item "{$(e)$}"  $v^{p_{\varepsilon +1}} = v^{p_\varepsilon} \cup
v^{p_{\gamma_\varepsilon,j_\varepsilon}}$
\sn
\item "{$(f)$}"  $v^{p_\varepsilon} = \cup
\{v^{p_{\gamma_\xi,j_\xi}}:\xi < \varepsilon\} \cup
v^{p^*}$.
\ermn
For $\varepsilon = 0$ use clause (c), for $\varepsilon$ limit use $(*)_3$
and clause (a) for successor $\varepsilon +1$ let $j_\varepsilon = 
\text{ Min}\{j:p_\varepsilon,p_{\gamma_\varepsilon,j}$ are compatible$\}$ 
(well defined as
$\langle p_{\gamma_\varepsilon,j}:j < \kappa \rangle$ is a maximal antichain
above $p^*$), and define $p_{\varepsilon +1}$ by $(*)_{10}$.
\sn
Let $p_\theta = \dbcu_{\varepsilon < \theta} p_\varepsilon$.
Lastly, define $q$:

$$
v^q = v^{p_\theta} \cup \{\zeta\}
$$

$$
v^q_* = v^{p_\theta}_* \cup \{\zeta\}
$$
\mn
It suffices to prove: $p^* \le q \in Q$ and $\varepsilon < \theta
\Rightarrow p_{\gamma_\varepsilon,j_\varepsilon} \le q$ as then $q$ forces
$A_\zeta$ is as required.  The non-trivial part is showing
\mr
\item "{$\boxtimes$}"  if $<^*_0$ is a witness to $p_{\gamma_\varepsilon,
j_\varepsilon}$ then some well ordering $<^*_1$ of $v^q_*$ which end extends
it is a witness to $q$.
\ermn 
First we can find $<^*_1$, a well ordering of $v^{p_\theta}_*$ end extending
$<^*_0$ and which is a witness to $p_\theta$.  We now define $<^*$, a well 
ordering of $v^q_*:<^* \restriction v^{p_\theta}_* = <^*_1$, and by $<^*,
\zeta$ is
just above all elements of $v^{p_{\gamma_\varepsilon,j_\varepsilon}}$ and
below all elements of $v^{p_\theta}_* \backslash
v^{p_{\gamma_\varepsilon,j_\varepsilon}}$.  Now $<^*$ is as required (note
that we have not proved $p_\theta \le q$!).
\hfill$\square_{\scite{g.9a}}$
\enddemo
\bigskip

\demo{\stag{g.10} Observation}  Assume that $\kappa = 
\kappa^{< \kappa} < \lambda$ and $S \subseteq \lambda$ stationary.  
\ub{Then} for some $\kappa^+$-c.c., 
strategically $\kappa$-complete forcing notion $Q$ of cardinality
$\lambda^{< \kappa}$, we have $\Vdash_Q ``S$ is the union of $\le \kappa$ 
sets each not reflecting any $\delta$ of cofinality $\le \kappa$.
\enddemo
\bigskip

\demo{Proof}  Straightforward.
\enddemo
\newpage

\head {\S3 Equi-consistency} \endhead  \resetall 
\bn
The following theorem clarifies the consistency strength of the problem
to a large extent.  We can hardly expect more with no inner models for
supercompacts.  Concentrating on ${}^\omega 2$ is for historical reason;
we can replace $\aleph_0$ by $\mu$.  Also, using the same claims we can
replace $\lambda > \beth_2$ by other restrictions.  Note \scite{3.2}
continues \cite[\S3]{Sh:460}, \cite{HJSh:249}.  The claims will give more, 
naturally.  However: \nl
\ub{\stag{3.0} Problem}:  What occurs if we demand GCH?
\bigskip

\proclaim{\stag{3.1} Theorem}  The following are equi-consistent (with ZFC +
$\kappa = \text{\rm cf\/}(\kappa) > 2^{\aleph_0}$
(in fact we use only forcing which preserves the cardinals $\le 
(2^{\aleph_0})^+$ and do not change the value of $2^{\aleph_0}$, in fact
the composition of $\kappa$-complete and c.c.c. of cardinality
$\le 2^{\aleph_0}$ ones; so we can
add $2^{\aleph_0} = \aleph_1$ or $2^{\aleph_0} = \aleph_2$ or
$2^{\aleph_0} = \aleph_{\omega^3 + \omega +3}$ or whatever)
\mr
\item "{$(A)[{}^\omega 2] = (A)_{(\omega_2)}$}"  there is a compact Hausdorff 
space $X$ such that $X \rightarrow_w ({}^\omega 2)^1_2$ but no subspace 
with $\le 2^{< \kappa}$ points has this property (on $\rightarrow_\omega$
see \scite{t.1}(2) and ${}^\omega 2$ is the Cantor discontinuum)
\sn
\item "{$(A)^+$}"  like $(A)_{({}^\omega 2)}$ replacing ${}^\omega 2$ by
``for any Hausdorff space $Y^*$ with $\le 2^{\aleph_0}$ points" and demand
$X$ has a clopen basis only if $Y$ has
\sn
\item "{$(B)[{}^\omega 2] = (B)_{(\omega_2)}$}"  there is a compact 
Hausdorff space $X$ with clopen basis such that $X \rightarrow ({}^\omega 2)^1
_{< \text{ cf}(2^{\aleph_0})}$ but no subspace with $\le 2^{< \kappa}$
points has this property
\sn
\item "{$(B)^+$}"  like $(B)[{}^\omega 2]$ replacing ${}^\omega 2$ by
``for any Hausdorff space with $\le 2^{\aleph_0}$ points" and demand $X$ has
a clopen basis only if $Y$ has
\sn
\item "{$(C)$}"  there are $\lambda,S,\bar f$ such that
{\roster
\itemitem{ $(a)$ }  $S \subseteq \lambda$ is stationary, $\lambda > 
2^{< \kappa}$
is regular
\sn
\itemitem{ $(b)$ }  $\bar f = \langle f_\delta:\delta \in S \rangle$
\sn
\itemitem{ $(c)$ }  $f_\delta$ is a one-to-one function from ${}^\omega 2$
to $\delta$
\sn
\itemitem{ $(d)$ }  if $\delta_1 \ne \delta_2$ then $\{\eta \in {}^\omega 2:
f_{\delta_1}(\eta) = f_{\delta_2}(\eta)\}$ has scattered closure (in the
topological space ${}^\omega 2$)
\endroster}
\item "{$(D)$}"  there are $\lambda,S,\bar A$ such that
{\roster
\itemitem{ $(a)$ }  $S \subseteq \lambda$ is stationary, $\lambda > 
2^{< \kappa}$ 
is regular
\sn
\itemitem{ $(b)$ }  $\bar A = \langle A_\delta:\delta \in S \rangle$
\sn
\itemitem{ $(c)$ }  $A_\delta$ is a subset of $\delta$ of cardinality
$2^{\aleph_0}$
\sn
\itemitem{ $(d)$ }  for $\delta_1 \ne \delta_2$ from $S$ we have
$A_{\delta_2} \cap A_\delta$ is finite
\sn
\itemitem{ $(e)$ }  $\{A_\delta:\delta \in S\}$ is $\kappa$-free 
\sn
\itemitem{ $(f)$ }  if $F:\lambda \rightarrow \alpha^*,\alpha^* < \lambda$
then for some $\delta,F \restriction A_\delta$ is constant.
\endroster}
\ermn
\endproclaim
\bn
Note that we can easily add clauses sandwiched between two existing ones.
\bn
\ub{Question}:  With what can we replace the space ${}^\omega 2$?

We make some definitions and prove some claims before the proof.  The
following definition is used in \scite{3.2}.
\definition{\stag{3.1A} Definition}  
1) For a cardinal $\kappa$ and $I_0,I_1 \subseteq
\{(a,b):a,b \subseteq \kappa \text{ are disjoint}\}$ and cardinal $\theta$
we say that a cardinal $\lambda$ is $(I_0,I_1,\theta)$-approximate or
$(\kappa,I_0,I_1,\theta)$-approximate \ub{if} we 
can find $\bar{\Cal P} = \langle {\Cal P}_\alpha:\alpha \in C \rangle$ such 
that
\mr
\widestnumber\item{$(iii)$}
\item "{$(i)$}"  $C$ a club of $\lambda$
\sn
\item "{$(ii)$}"  ${\Cal P}_\alpha \subseteq [\alpha]^{< \theta}$ for
$\alpha \in C$ and $|{\Cal P}_\alpha| \le \text{ Min}(C \backslash (\alpha +
1))$
\sn
\item "{$(iii)$}"  for any 1-to-1 function $f$ from $\kappa$ to $\lambda$,
for some $\alpha \in C$ one of the following holds
{\roster
\itemitem{ $(a)$ }  for some $c \in {\Cal P}_\alpha$ and $(a,b) \in I_1$ we
have $(\forall i \in a)(f(i) \in c)$ and 
$(\forall i \in b)[f(i) \ge \alpha]$
\sn
\itemitem{ $(b)$ }  for some $(a,b) \in I_0$ we have
$$
(\alpha) \qquad (\forall i < \kappa)(f(i) < \alpha \rightarrow i \in a)
$$

$$
(\beta) \qquad (\forall i < \kappa)
[i \in b \rightarrow \alpha \le f(i) < \text{ Min}(C \backslash (\alpha +1))]
$$
\endroster}
\ermn
2) If $c \ell$ is a function from ${\Cal P}(\kappa)$ to ${\Cal P}
(\kappa)$ and $K \subseteq {\Cal P}(\kappa)$ and

$$
I_1 = \{(a,b):a \subseteq \kappa,b \in K \text{ and } b \subseteq
c \ell(a)\}
$$

$$
I_0 = \{(a,b):a \subseteq \kappa,b \in K \text{ and } a \cap b =
\emptyset\}
$$
\mn
\ub{then} we may say $\lambda$ is 
$(K,c \ell,\theta)$-approximate or $(\kappa,K,c \ell,\theta)$-approximate
instead of $\lambda$ is $(I_0,I_1,\theta)$-aproximate. \nl
3) We may replace $\kappa$ by another set of this kind call the domain of
the tuple (understood from $I_0,I_1$).  We may write this set before 
$I_0$ for clarification. \nl
4) We may replace $(I_0,I_1,\theta)$ by $(\bold I,\theta)$ if $\bold I$
is a set of pairs $(I_0,I_1)$ such that $\langle {\Cal P}_\alpha:\alpha
\in C \rangle$ satisfies the requirement above all the triples
$(I_0,I_1,\theta)$ such that $(I_0,I_1) \in \bold I$ (not
necessarily all pairs have the same domain $A$). \nl
Similarly, $\bold K$ stands for a set of tuples $(\kappa,K,c \ell,\theta)$
or in short $(\kappa,K,c \ell)$ when $\theta$ is understood from the context
or even $(K,c \ell)$ as in part (2).  
(We may even vary $\theta$).
\bn
Concerning \scite{3.f3} below \nl
\ub{\stag{3.1B} Examples}:  1) Let $\bold C$ 
be a Cantor set (say ${}^\omega 2$), \nl
$c \ell^{\bold C}$ is the (topological) closure operation on subsets of
$\bold C$ \nl
$K^{\bold C} = \{A \subseteq \bold C:A \text{ is closed perfect uncountable}
\}$. \nl
2) Let $\Bbb R$ be the real line, $c \ell^{\Bbb R}$ be the (topological)
closure operation on subsets of $\Bbb R$ and $K^{\Bbb R} = \{A \subseteq
\Bbb R:A \text{ is closed perfect uncountable, bounded (from below and
above)}\}$.
\enddefinition
\bigskip

\proclaim{\stag{3.2} Lemma}  Assume
\mr
\item "{$(a)$}"  $\lambda > \chi \ge \kappa \ge \theta$ and $\sigma$ are
infinite cardinals,
\sn
\item "{$(b)$}"  $c \ell$ is a partial function from $[\lambda]^{< \theta}$
to $K \subseteq [\lambda]^{\le \kappa}$
\sn
\item "{$(c)$}"  $\bold K$ is a set of triples $(\kappa,K^*,c \ell^*)$ 
with $K^* \subseteq {\Cal P}(\kappa),c \ell^*$ a function from 
$[\kappa]^{< \theta}$ to ${\Cal P}(\kappa)$ as in
Definition \scite{3.1A}(2) above (for $\theta$)
\sn
\item "{$(d)$}"  if $b \in K$, then for some $(\kappa,K^*,c \ell^*) \in 
\bold I$ and one to one function $f$ from $\kappa$ into $b$, we have:

$$
b' \in K^* \Rightarrow \{f(\alpha):\alpha \in b'\} \in K
$$

$$
a',b' \subseteq \kappa \and c \ell^*(a') = b' \Rightarrow c \ell\{f(\alpha):
\alpha \in a'\} \supseteq \{f(\alpha):\alpha \in b'\}
$$
\sn
\item "{$(e)$}"  for every $A \in [\lambda]^{\le \chi}$
we can find a $[K,\sigma]$-colouring $\bold c$ of $A$, which is defined
for any $A \subseteq \lambda$ as saying that: $\bold c$ is a function 
from $A$ to $\sigma$ such that $a \in K \and a \subseteq A \Rightarrow 
\text{ Rang} (\bold c \restriction a) = \sigma$
\sn
\item "{$(f)$}"  for every $\mu$, if 
$\chi < \mu \le \lambda$ then $\mu$ is $(\bold K,\theta)$-approximate.
\ermn
\ub{Then} there is $[K,\sigma]$-colouring $\bold c$ of $\lambda$.
\endproclaim
\bigskip

\demo{Proof}  See after the proof of \scite{3.f10} below.  (The reader may
prefer to read first \S4 up to the proof of \scite{3.2}, \scite{3.5} first).
\enddemo
\bn
\ub{\stag{3.3} Conclusion}:  1) Assume
\mr
\item "{$(a)$}"  every cardinal $\mu,2^{\aleph_0} < \mu \le \lambda$ is
$(\bold C,K^{\bold C},c \ell^{\bold C},\aleph_1)$-approximate
\sn
\item "{$(b)$}"  $X$ is a Hausdorff topological space.
\ermn
Then $X \nrightarrow$ [Cantor set]$^1_{2^{\aleph_0}}$. \nl
2) We can replace in part (1), $\bold C$ by $\Bbb R$. 
\bigskip

\demo{Proof}  By \scite{3.2} (and \scite{3.1B}).\hfill$\square_{\scite{3.3}}$
\enddemo
\bigskip

\proclaim{\stag{3.4} Claim}  The forcing notions in \scite{t.2} and 
in \scite{g.8} satisfies e.g. the condition $*^{\sigma^+}_{\kappa^+}$; 
see below Definition \scite{3.cc}(1A).
\endproclaim
\bigskip

\definition{\stag{3.cc} Definition}  1) Let $D$ be a normal filter on
$\mu^+$ to which $\{\delta < \mu^+:\text{cf}(\delta) = \mu\}$ belongs.
A forcing notion $Q$ satisfies $*^\epsilon_D$ where $\epsilon$ is a 
limit ordinal $< \mu$, if player I has a winning strategy in the 
following game $*^\epsilon_D[Q]$ defined as follows:
\newline
\noindent
\underbar{Playing}: the play finishes after  $\epsilon$ moves.   

In the $\zeta$-th move:   

Player I --- if  $\zeta \ne 0$  he chooses  
$\langle q^\zeta_i:i < \mu^+ \rangle$ such that $q^\zeta_i \in Q$
\roster
\item "{{}}"  $\quad \quad$
and $(\forall \xi < \zeta)(\forall^D i < \mu^+)
p^\xi_i \le q^\zeta_i$ and he chooses a 
\item "{{}}"   $\quad \quad$
function $f_\zeta:\mu^+ \rightarrow \mu^+$ such that for a club of 
$i < \mu^+,f_\zeta(i) < i$;  
\item "{{}}"  $\quad \quad$
if $\zeta = 0$ let  $q^\zeta_i = \emptyset_Q$, $f_\zeta =$ is identically
zero.
\endroster   
\medskip

Player II --- he chooses  $\langle p^\zeta_i:i < \mu^+ \rangle$  
such that $(\forall^D i) q^\zeta_i \le p^\zeta_i$ and $p^\zeta_i \in 
Q$.
\medskip

\noindent
\underbar{The Outcome}:  Player I wins provided that for some $E \in D$: 
if \newline
$\mu < i < j < \mu^+,i,j \in E$,
$cf(i) = cf(j) = \mu$ and $\dsize \bigwedge_{\xi < \epsilon} f_\xi
(i) = f_\xi(j)$ \underbar{then} the set  
$\{p^\zeta_i:\zeta < \epsilon \} \cup \{p^\zeta_j:\zeta <
\epsilon \}$  has an upper bound in  $Q$. \newline
1A) If $D$ is $\{A \subseteq \mu^+:\text{for some club } E \text{ of }
\mu^+ \text{ we have } i \in E \and \text{ cf}(i) = \mu \Rightarrow i \in A\}$
we may write $\mu$ instead of $D$ (in $*^\varepsilon_D$ and in the related
notions defined below and above). \nl
2) We may allow the strategy to be non-deterministic, e.g. choose not
$f_\zeta$ just $f_\zeta/D_{\mu^+}$. \nl
3) We say a forcing notion $Q$ is $\varepsilon$-strategically complete if for
the following game, $\bigotimes^\varepsilon_Q$ player I has a winning
strategy. \nl
In the $\zeta$-th move: \nl
Player I - if $\zeta \ne 0$ he chooses $q_\zeta \in Q$ such that
$(\forall \xi < \zeta) p_\xi \le q_\zeta$ if $\zeta = 0$ let $q_\zeta =
\emptyset_Q$. \nl
Player II - he chooses $p_\zeta \in Q$ such that $q_\zeta \le p_\zeta$.
\mn
\ub{The Outcome}:  In the end Player I wins provided that he always has a
legal move.
\enddefinition
\bigskip

\proclaim{\stag{3.ccd} Lemma}  The property ``$Q$ is 
$(< \mu)$-strategically complete and has $*^\varepsilon_\mu$" is preserved
by $(< \kappa)$-support iteration (see \cite{Sh:546}).  FILL?
\endproclaim
\bigskip

\demo{Proof}  Straight; in each coordinate we preserve that the sequence of
conditions is increasingly continuous and on each stationary $S \subseteq
\{\delta < \kappa^+:\text{cf}(\delta) = \kappa\}$ on which the pressing down
function is constant the conditions form a $\Delta$-system.
\hfill$\square_{\scite{3.ccd}}$
\enddemo
\bn
We can also consider
\definition{\stag{3.4a} Definition}  1) We say $X^* \rightarrow [Y^*]^n_\theta$
\ub{if} $X^*,Y^*$ are topological spaces and for every $h:[X^*]^n \rightarrow
\theta$ there is a closed subspace $Y$ of $X^*$ homeomorphic to $Y^*$ such
that for some $\alpha < \theta,\alpha \notin \text{ Rang}(h \restriction
[Y]^n)$ is not $\theta$. \nl
2) If we omit the ``closed" we shall write $\rightarrow_w$ instead of
$\rightarrow$ and $\nrightarrow,\nrightarrow_w$ denote the negations.
[FILL? \scite{3.1}, \scite{3.2}.]
\enddefinition
\bigskip

\proclaim{\stag{3.5} Claim}  1) Assume $X$ is a Hausdorff space with 
$\lambda$ points and $D$ is a filter on ${}^\omega 2$ containing the
co-countable subsets of ${}^\omega 2$.  Assume further $X \nrightarrow
[{}^\omega 2]^1_\theta$ and $\chi \ge 2^{\aleph_0}$ but no subspace $X^*$ of
$X$ with $< \chi$ points satisfy $X^* \nrightarrow ({}^\omega 2)^1
_{2^{\aleph_0}}$ and $\chi = \chi^{\aleph_0}$.  \ub{Then}
\mr
\item "{$(*)$}"  we can find a regular $\kappa \in (\chi,\lambda]$, a
stationary $S \subseteq \kappa$ and $\bar f = \langle f_\alpha:\alpha \in S
\rangle$ such that:
{\roster
\itemitem{ $(i)$ }  Dom$(f_\alpha) \in D^+$
\sn
\itemitem{ $(ii)$ }  $f_\alpha$ is one-to-one and is a homeomorphism from
${}^\omega 2 \restriction \text{ Dom}(f_\alpha)$ onto $X \restriction
\text{ Rang}(f_\alpha)$
\sn
\itemitem{ $(iii)$ }  if $\alpha \ne \beta$ are from $S$, \ub{then}
$\{\eta \in \text{ Dom}(f_\alpha):f_\alpha(\eta) \in \text{ Rang}(f_\beta)\}$
has scattered closure in ${}^\omega 2$
\sn
\itemitem{ $(iv)$ }  for a club of $\delta \in S$ we have 
Rang$(f_\alpha) \subseteq \dbcu_{\beta \in \alpha \cap S} 
\text{ Rang}(f_\beta)$.
\endroster}
\ermn
2) Similarly for $\nrightarrow_w$. \nl
We shall prove it later (afer the proof of \scite{3.f10}.
\endproclaim
\bn
\ub{\stag{3.6} Observation}:  There is a c.c.c. forcing notion $Q$ of
cardinality $2^{\aleph_0}$ such that:

$$
\align
\Vdash_Q ``&\text{there is }h:{}^\omega 2 \rightarrow \omega 
\text{ such that}: \\
  &(\alpha) \quad \text{ if } C (\in V) \text{ is closed scattered then
each} \\
  &\qquad \quad C \cap h^{-1}\{n\} \text{ is finite, and} \\
  &(\beta) \quad  \text{if } A \subseteq ({}^\omega 2) 
\text{ is uncountable and from }  V \\
  &\qquad \text{then } |A \cap h^{-1}\{n\}| = |A| \text{ for each } n".
\endalign
$$
\bigskip

\demo{Proof}  Let $p \in Q$ be $(f^p,{\Cal C}^p)$ where $f^p$ is a finite
function from ${}^\omega 2$ to $\omega$ and ${\Cal C}^p$ is a finite family
of closed scattered subsets of ${}^\omega 2$.
\sn
Order is: \nl
$p \le q$ iff $f^p \subseteq f^q,{\Cal C}^p \subseteq {\Cal C}^q$ 
and $C \in {\Cal C}^p \and \eta \in C \cap \text{ Dom}(f^p) \and \nu \in 
\text{ Dom}(f^q) \backslash \text{Dom}(f^p) \and \eta \ne
\nu \rightarrow f^q(\eta) \ne f^q(\nu)$.
\sn
Clearly
\mr
\item "{$(*)_1$}"  $Q$ is a forcing notion of cardinality $2^{\aleph_0}$
\sn
\item "{$(*)_2$}"  $Q$ satisfies the c.c.c. \nl
[why?  let $p_\alpha \in Q$ for $\alpha < \omega_1$, let Dom$(f_\alpha) =
\{\eta_{\alpha,\ell}:\ell < \ell_\alpha\},{\Cal C}^{p_\alpha} = \{
C_{\alpha,k}:k < k_\alpha\}$ and let $m_\alpha = \text{ Min}\{m:\langle
\eta_{\alpha,\ell} \restriction m:\ell < \ell_\alpha \rangle$ is with no
repetitions$\}$.  Without loss of generality $m_\alpha = m(*),\ell_\alpha
= \ell(*),k_\alpha = k(*),\eta_{\alpha,\ell} \restriction m(*) = \nu_\ell$.
By $\Delta$-system lemma \wilog \, for some $\ell(**) \le \ell(*)$ we have:
\mn
$(\alpha) \qquad \ell < \ell(**) \Rightarrow 
\langle \eta_{\alpha,\ell}:\alpha < \omega_1 \rangle 
\text{ is with no repetitions}$
\sn
$(\beta) \qquad \ell \in [\ell(**),\ell(*)) 
\Rightarrow \eta_{\alpha,\ell} = \eta_\alpha$
\sn
$(\gamma) \qquad \{\eta_{\alpha,\ell}:\alpha < \omega_1,\ell < \ell(**)\}
\text{ is with no repetitions}$. 
\mn
Now as each $C_{\alpha,k}$ is closed and scattered it is necessarily countable
so \wilog

$$
\alpha < \beta < \omega_1 \and \ell < \ell(**) \Rightarrow
\eta_{\beta,\ell} \notin \dbcu_{k < k(*)} C_{\alpha,k}.
$$
\mn
We now choose by induction on $\ell < \ell(**)$ sets $A_\ell,B_\ell \in
[\omega_1]^{\aleph_1}$, decreasing with $n$ such that

$$
\alpha \in A_{\ell +1} \and \beta \in B_{\ell +1} \and \alpha < \beta
\rightarrow \eta_{\alpha,\ell} \notin \dbcu_{k < k(*)} C_{\beta,k}.]
$$
\sn
This is straight: let $A_\gamma = \omega_1 = B_0$; if $A_\ell,B_\ell$ are 
given, then for some $\alpha^*_\ell
\in A_\ell$ the set $\{\eta_{\alpha,\ell}:\alpha \in A_\ell \backslash
\alpha^*_\ell\}$ is $\aleph_1$-dense in itself, i.e. $(\forall \alpha \in
A_\ell \backslash \alpha^*_\ell)(\forall n < \omega)(\exists^{\aleph_1} 
\beta \in A_\ell)(\eta_{\beta,\ell} \restriction n =
\eta_{\alpha,\ell} \restriction n)$, let $T_\ell = \{\eta_{\alpha,\ell}
\restriction n:\alpha \in A_\ell \backslash \alpha^*_\ell\}$.  
So for each $\beta \in B_\ell$ for some $\nu^\ell_\beta \in T_\ell$ we have
$(\forall \rho \in \dbcu_k C_{\beta,k})(\neg \nu^\ell_\beta \triangleleft
\rho)$ so for some $\nu_\ell \in T_\ell$ we have $B_{\ell +1} =:
\{\beta \in B_\ell:\nu^\ell_\beta = \nu_\ell\}$ is uncountable and let
$A_{\ell +1} = \{\alpha \in A_\ell:\nu_\ell \triangleleft 
\eta_{\alpha,\ell}\}$. \nl
For $\alpha < \beta,\alpha \in A_{\ell(**)},\beta \in
B_{\ell(**)}$, we have $p_\alpha,p_\beta$ are compatible.]
\sn
\item "{$(*)_3$}"  if $A \subseteq {}^\omega 2$ is uncountable and 
$n < \omega$ then ${\Cal I}_{A,n} =:\{p:\text{for some } \eta \in A,
f^p(\eta) = n\}$ is dense open \nl
[why?  as $C^* = \cup\{C:C \in {\Cal C}^p\}$ is closed and scattered hence
countable clearly for some $\eta \in A$ we have $\eta \notin C^*$ so 
$q = (f^p \cup \{(\eta,n)\},{\Cal C}^p)$ satisfies $p \le q \in Q \cap 
{\Cal I}_{A,n}$.]
\sn
\item "{$(*)_4$}"  for each $\eta \in {}^\omega 2$ the set \nl
${\Cal I}_\eta = \{p:\eta \in \text{ Dom}(f^p)\}$ is dense open \nl
[why?  for $p \in Q$ let $n = \sup(\text{Rang}(f^p))+1$ and letting
$q = (f^p \cup\{(\eta,n)\},{\Cal C}^p)$ satisfies $p \le q \in 
{\Cal I}_\eta$.]
\sn
\item "{$(*)_5$}"  for each closed scattered $C$, the set ${\Cal I}_C =
\{p:C \in {\Cal C}^p\}$ is dense open \nl
[why? immediate as $p \in Q \Rightarrow p \le (f^p,{\Cal C}^p \cup \{C\}) \in
Q$.]
\ermn
Let $\underset\tilde {}\to f = \cup\{f^p:p \in \underset\tilde {}\to G\}$
\mr
\item "{$(*)_6$}"  $\underset\tilde {}\to f$ is a function from
$({}^\omega 2)^v$ to $\omega$ and for each closed scattered $C \in V,
f \restriction C$ is one to one except on a finite set \nl
[why?  easy].
\ermn
Together we are done.  \hfill$\square_{\scite{3.6}}$
\enddemo
\bigskip

\demo{Proof of \scite{3.1}}
\sn
\ub{$(B)^+ \Rightarrow (B)[{}^\omega 2]$}

Trivial (special case).
\mn
\ub{$(A)^+ \Rightarrow (A)[{}^\omega 2]$}

Trivial (a special case).
\mn
\ub{$(B)^+ \Rightarrow (A)^+$}

Trivial (stronger demands).
\mn
\ub{$(B)[{}^\omega 2] \Rightarrow (A)[{}^\omega 2]$}

Trivial (stronger demands).
\mn
\ub{$(A)[{}^\omega 2] \Rightarrow (C)$}

By \scite{3.5} for the filter $D$, using the monotonicity of
$X \rightarrow ({}^\omega 2)^1_\theta$ in 
$\theta = \{2^\omega \backslash X:X \subseteq {}^\omega 2$: and the 
closure of $X$ in ${}^\omega 2$ is scattered (equivalently countable)$\}$.
\mn
\ub{$(C) \Rightarrow (D)$}

Forcing by Levy$(\kappa,2^{< \kappa})$ change nothing so \wilog \,
$\kappa = \kappa^{< \kappa}$.
Let $\mu,S,\bar f = \langle f_\alpha:\alpha \in S \rangle$ as there.  Next
let $Q$ be the forcing notion from \scite{3.6} so we get the conclusion of
\scite{3.6} for the ideal $[A]^{< \aleph_0}$ where $A \subseteq
({}^ \omega 2)^V$ has cardinality $2^{\aleph_0}$.  So letting
$A_\alpha = \text{ Rang}(f_\alpha \restriction A)$ we get: $A_\alpha
\subseteq \alpha,|A_\alpha| = 2^{\aleph_0}$ and for $\alpha \ne \beta$ from
$S,A_\alpha \cap A_\beta$ is finite.  So clauses (a)-(d) of (D) holds.  Then
we force by Levy$(\lambda,2^{< \lambda})$ nothing changes but we get
$\diamondsuit_S$.  By \scite{g.10} \wilog \, $S$ does not reflect in ordinal
$\delta$ if cf$(\delta) \le \kappa$.  So $(*)$ of \scite{g.6} holds hence
by \scite{g.6A} we get the $\kappa$-freeness (clause (e) of
$(D) \equiv (B)_2$ of \scite{t.2}) and clause (f) (= clause (C) of
\scite{t.2}). \nl
[Question:  \scite{g.9a} not needed.]
\mn
\ub{$(D) \Rightarrow (B)^+$}

We do it by forcing but for the proof any $\kappa$ 
such that $\aleph_1 \le \text{ cf}(\kappa) = \kappa,2^{< \kappa} < 
\lambda$ can serve. \nl
If $\kappa > 2^{\aleph_0}$ (as in the main case) and we restrict ourselves
to spaces $Y^*$ with a basis of cardinality $< \kappa$, \ub{then} we can use
a product instead of iteration for $P_{i(*)}/P_2$ below and proof is easier.
By forcing by Levy$(\lambda,2^{< \lambda})$ (see \scite{3.7}) \wilog \, 
$\diamondsuit_S$ for the $S$ of (D), this will be
preserved for any forcing notion $P$ if $P$ has density $\le \lambda$, which
is our case.  
Now we use iterated forcing $\langle P_j,Q_i:j \le i(*),i < i(*) \rangle$
with $(< \kappa)$-support, each satisfying the $*^{\sigma^+}_{\kappa^+}$
version of $\kappa^+$-c.c. (see \scite{3.4}).  Now let $Q_0$ be
as Levy$(\kappa,2^{< \kappa})$ and each
${\underset\tilde {}\to Q_{1+i}}$ be as in \scite{t.2} for some
${\underset\tilde {}\to Y^*_{1+i}}$ (a $P_i$-name of a topological space as in
\scite{t.2}) and it forces an example ${\underset\tilde {}\to X^*_{1+i}}$.  
With suitable bookkeeping (if $\kappa > 2^{\aleph_0}$ is
easier) we finish as those iterations preserve ``$i < \kappa$-strategic
completeness hence no new set of ordinals of
cardinality $< \kappa$ and (the strong version of) $\kappa^+$-c.c." is
preserved.  \nl
Still we have to prove that the example $X^*_{1+i}$ we force to satisfy
``$X^*_{1+i} \rightarrow (Y^*_{1+i})^1_\sigma$" has this property not only in
$V^{P_{1+i+1}}$ but also in $V^{P_{i(*)}}$.  For this we repeat the relevant
part of the proof of \scite{t.2} together with the preservation of
$(*)^{\sigma^+}_{\kappa^+}$.
\hfill$\square_{\scite{3.1}}$
\enddemo
\newpage

\head {\S4 Helping equiconsistency} \endhead  \resetall 
\bigskip

\definition{\stag{3.f1} Definition}  1) Let ${\Cal Y}$ denote a set of pairs
of the form $(I,J)$ where $I \subseteq J$ are ideals over a common set
called Dom$(I,J) = \text{ Dom}((I,J))$.  Let $\kappa({\Cal Y}) = \sup
\{|\text{Dom}(I,J)|:(I,J) \in {\Cal Y}\}$.  We call ${\Cal Y}$ standard if
for each $(I,J) \in {\Cal Y}$, the set Dom$(I,J)$ is a cardinal. \nl
2) NFr$_1(\lambda,{\Cal Y})$ if for some $\lambda^* > \lambda$ we have
NFr$_1(\lambda^*,\lambda,{\Cal Y})$ which means $\lambda \ge |{\Cal Y}| + 
\kappa({\Cal Y})$ and there are $\langle {\Cal F}_{(I,J)}:(I,J) \in {\Cal Y}
\rangle$ exemplifying it which means:
\mr
\item "{$(a)$}"  ${\Cal F}_{(I,J)} \subseteq \{f:f \text{ a function, Dom}
(f) \in J^+\}$
\sn
\item "{$(b)$}"  if $f \ne g \in {\Cal F}_{(I,J)}$ then \nl
$\{x:x \in \text{ Dom}(f) \cap \text{ Dom}(g) \text{ but } f(x) \ne g(x)\}$
belong to $I$ 
\sn
\item "{$(c)$}"  $\lambda \ge |\cup\{\text{Rang}(f):f \in {\Cal F}_{(I,J)}$
and $(I,J) \in {\Cal Y}\}|$
\sn
\item "{$(d)$}"  $\lambda < \lambda^* = \sum\{|{\Cal F}_{(I,J)}|:
(I,J) \in {\Cal Y}\}$.
\ermn
2) NFr$_2(\lambda,{\Cal Y})$ if $\lambda$ is regular $> |{\Cal Y}| + \theta
(Y)$ and there are $(I,J) \in {\Cal Y}$ and $\langle f_\delta:\delta \in S
\rangle$ such that
\mr
\item "{$(a)$}"  $S \subseteq \lambda$ is stationary
\sn
\item "{$(b)$}"  Dom$(f_\delta) \in J^+$, Rang$(f_\delta) \subseteq \delta$
\sn
\item "{$(c)$}"  $\delta_1 \ne \delta_2 \Rightarrow \{x:x \in \text{ Dom}
(f_{\delta_1}) \cap \text{ Dom}(f_{\delta_2})$ and $f_{\delta_1}(x) =
f_{\delta_2}(x)\} \in I$.
\ermn
3) We omit $N$ from NFr in parts (1) and (2) for the negation.
\enddefinition
\bn
\ub{\stag{3.f2} Fact}:  1) NFr$_1(\lambda,{\Cal Y})$ is preserved by 
increasing ${\Cal Y}$ to ${\Cal Y}'$ when $|{\Cal Y}'| + \kappa({\Cal Y}')
< \lambda$.  Also NFr$_2(\lambda,{\Cal Y})$ is preserved by increasing
${\Cal Y}$.  Similarly if
NFr$_1(\lambda^*,\lambda,\mu,{\Cal Y}),\lambda^* \ge \lambda^*_1 > \lambda_1
\ge \lambda,\lambda_1 \ge |{\Cal Y}_1| + \theta({\Cal Y}_1)$ and ${\Cal Y}_1
\supseteq Y$ then NFr$_1(\lambda^*_1,\lambda_1,{\Cal Y}_1)$. \nl
2)  NFr$_1(\lambda,{\Cal Y})$ is equivalent to NFr$_1(\lambda^+,\lambda,
{\Cal Y})$. \nl
3) If $\lambda^*$ is regular or at least cf$(\lambda^*) > |{\Cal Y}|$ and
NFr$_1(\lambda^*,\lambda,{\Cal Y})$, \ub{then} there is 
$(I,J) \in {\Cal Y}$ such that NFr$_1(\lambda^*,\lambda,\mu,\{(I,j)\})$. \nl
4) NFr$_1(\lambda,{\Cal Y})$ implies NFr$_2(\lambda^+,{\Cal Y})$.
\bigskip

\demo{Proof}  Check.
\enddemo
\bigskip

\proclaim{\stag{3.7} Claim}  Asume NFr$_2(\lambda,\{(I,J)\})$ and let
$\bar A = \langle A_\delta:\delta \in S \rangle$ exemplifies it and let
$\mu = (\text{Dom}(I,J))$ and assume $\mu^{++} < \lambda$. \nl
1) If $\diamondsuit_S$ then we can find $\langle A'_\delta:\delta \in S
\cap E \rangle$ exemplifying NFr$_2(\lambda,\{(I,J)\}),E$ a club of $\lambda$
such that
\mr
\item "{$(*)$}"  if $\tau^+ < \lambda,F:\lambda \rightarrow [\lambda]
^{\le \tau}$ \ub{then} for some $\delta \in S \cap E$ the set $A_\delta$ is
$F$-free (i.e. $\alpha \ne \beta \in A_\delta \Rightarrow \beta \notin
F(\alpha)$. 
\ermn
2) The forcing of adding a Cohen subset of $\lambda$ (i.e. $({}^{\lambda >}
2,\triangleleft))$ preserve ``$\bar A$ exemplifies NFr$_2(\lambda,\{(I,J)\})$"
(as it preserves ``$S$ is stationary"), add no bounded subsets to $\lambda$
and forces $\diamondsuit_S$.
\endproclaim
\bigskip

\demo{Proof}  1) As in \scite{g.5}. \nl
2) Straight.  \hfill$\square_{\scite{3.7}}$
\enddemo
\bn
\centerline {$* \qquad * \qquad *$}
\bn
Now we give sufficient conditions for the existence of colouring.
\proclaim{\stag{3.f3} Claim}  Assume:
\mr
\item "{$(a)$}"  ${\Cal Y}$ is as in Definition \scite{3.f1}(1)
\sn
\item "{$(b)$}"  $\lambda > \mu \ge |{\Cal Y}| + \kappa({\Cal Y})$
\sn
\item "{$(c)$}"  for no regular $\kappa \in (\mu,\lambda]$ do we have
NFr$_2(\kappa,{\Cal Y})$
\sn
\item "{$(d)$}"  $c \ell$ is a function from $[\lambda]^{\le \mu}$ to
$[\lambda]^{\le \mu}$
\sn
\item "{$(e)$}"  for $A,B \in [\lambda]^{\le \mu}$ we have $A \subseteq
c \ell(A) = c \ell(c \ell(A))$ and $A \subseteq B \Rightarrow c \ell(A)
\subseteq c \ell(B)$
\sn
\item "{$(f)$}"  ${\Cal P} \subseteq [\lambda]^{\le \mu}$ satisfies: \nl
for every $A \in {\Cal P}$ there are a pair 
$(I,J) \in {\Cal Y}$, a set ${\Cal U} \in J^+$
and one to one $f:{\Cal U} \rightarrow A$ such that:
{\roster
\itemitem{ $(\alpha)$ }  if \,${\Cal U}' \subseteq {\Cal U} \and
{\Cal U}' \in I^+$ \ub{then} for some $A' \in {\Cal P}$ we have $A' \subseteq
A \cap c \ell(\{f(i):i \in {\Cal U}'\})$
\sn
\itemitem{ $(\beta)$ }  there are ${\Cal U}'_\alpha \subseteq {\Cal U},
{\Cal U}'_\alpha \in I^+$ for $\alpha < \alpha^*$ for some $\alpha^* \le \mu$
such that for any ${\Cal U}' \subseteq {\Cal U},{\Cal U}' \in I^+$ for some
$\alpha < \alpha^*$ we have ${\Cal U}'_\alpha \subseteq {\Cal U}'$ or at
least $c \ell(\{f(i):i \in {\Cal U}'_\alpha\}) \subseteq c \ell \{f(i):i \in
{\Cal U}'\}$
\sn
\itemitem{ $(\gamma)$ }  $c \ell(A) = A$.
\endroster}
\ermn
\ub{Then}
\mr
\item "{$Dec(\lambda,{\Cal P},\mu,{\Cal Y})$}"  for every $\chi > \lambda$
and $x \in {\Cal H}(\chi)$ there is $\langle M_i:i < \lambda \rangle$
such that:
{\roster
\itemitem{ $(i)$ }  $M_\alpha \prec ({\Cal H}(\chi),\in,<^*_\chi)$
\sn
\itemitem{ $(ii)$ }  $\mu \cup \{{\Cal Y},\lambda,\mu,x\} \subseteq M_\alpha$
and $\|M_\alpha\| = \mu$
\sn
\itemitem{ $(iii)$ }  $\dbcu_{\alpha < \lambda} M_\alpha$ includes
$\lambda$
\sn
\itemitem{ $(iv)$ }  assume $A \in {\Cal P}$ and define 
$\alpha(A) = \text{ Min}\{\alpha \le \lambda:\text{if } \alpha < \lambda 
\text{ then for some }
(I,J) \in {\Cal Y} \text{ and } {\Cal U} \in J^+$ and $f:{\Cal U} \rightarrow
A$ which is one-to-one, we have $\{i \in {\Cal U}:f(i) \in \dbcu_{\beta \le
\alpha} M_\beta\} \in J^+\}\}$ and letting $(I,J),{\Cal U},f$ be witnesses to
$\alpha(A) = \alpha$ we have: \nl
$\{i \in {\Cal U}:f(i) \in M_\alpha \backslash \dbcu_{\beta < \alpha}
M_\beta\} \in J^+$ \nl
moreover for some $X \in M_\alpha$ of cardinality $\le \mu$ (so $X \subseteq
M_\alpha$) we have $\{i \in {\Cal U}:f(i) \in X \backslash \dbcu_{\beta <
\alpha} M_\beta\} \in J^+$
\sn
\itemitem{ $(v)$ }  for any pregiven $\sigma = \text{ cf}(\sigma) \le \mu$
we can demand $M_\alpha = \dbcu_{\varepsilon < \sigma} M_{\alpha,\varepsilon},
M_{\alpha,\varepsilon}$ increasing with $\varepsilon$ and $\langle
M_{\alpha,\zeta}:\zeta \le \varepsilon \rangle \in M_{\alpha,\varepsilon +1}$
and $\langle M_\beta:\beta < \alpha \rangle \in M_{\alpha,\varepsilon}$.
\endroster}
\endroster
\endproclaim
\bn
Now below we shall prove Claim \scite{3.f3} follows from the 
following variant (we change (d), (e), (f)).
\proclaim{\stag{3.f3a} Claim}  Assume
\mr
\item "{$(a)'$}"  ${\Cal Y}$ is as in Definition \scite{3.f1}(1)
\sn
\item "{$(b)'$}"  $|X| = \lambda > \mu \ge |{\Cal Y}| + \kappa({\Cal Y})$
\sn
\item "{$(c)'$}"  for no regular $\kappa \in (\mu,\lambda]$ do we have
NFr$_2(\kappa,{\Cal Y})$
\sn
\item "{$(d)'$}"  $\bar{\Cal F} = \langle {\Cal F}_t:t \in T \rangle,T$ is a
partial order; we consider the ${\Cal F}_t$'s as index sets such that
$t \ne s \Rightarrow {\Cal F}_s \cap {\Cal F}_t = \emptyset$
\sn
\item "{$(e)'$}"  each member $f \in \dbcu_{t \in T} {\Cal F}_t$ is a
function such that for some $(I,J) = (I_f,J_f) \in {\Cal Y}$ we have
Dom$(f) \in J^+$, Rang$(f) \subseteq X$
\sn
\item "{$(f)'$}"  if $t \in T$ and $f \in {\Cal F}_t$, \ub{then} there 
is a subset $T[f]$ of $T_{<f>}$ of cardinality $\le \mu$ which is a cover
which means $(\forall s \in T_{<f>})(\exists t \in T[f])[s \le_T t]$ where
$$
\align
T_{<f>} =: \{r \in T:&\text{for some } g \in {\Cal F}_r \text{ we have }
(I_g,J_g) = (I_f,J_f) \text{ and} \\
  &\{i:i \in \text{ Dom}(f),i \in \text{ Dom}(g) \text{ and } f(i) = 
  g(i)\} \in I^+_g = I^+_f\}.
\endalign
$$
\ermn
\ub{THEN}
\mr
\item "{$Dec(\lambda,\bar{\Cal F},\mu,{\Cal Y})$}"  for every $\chi > \lambda$
and $x \in {\Cal H}(\chi)$ there is $\langle M_i:i < \lambda \rangle$
such that:
{\roster
\itemitem{ $(i)$ }  $M_\alpha \prec ({\Cal H}(\chi),\in,<^*_\chi)$
\sn
\itemitem{ $(ii)$ }  $\mu \cup \{{\Cal Y},\lambda,\mu,x\} \subseteq M_\alpha$
and $\|M_\alpha\| = \mu$
\sn
\itemitem{ $(iii)$ }  $\dbcu_{\alpha < \lambda} M_\alpha$ includes
$\lambda$
\sn
\itemitem{ $(iv)$ }  if $s \in T$, \ub{then} for some $t,s \le_T t \in T$
and for some $\alpha < \lambda$ and $g \in {\Cal F}_t$ we have \nl

$\qquad (\alpha) \qquad \{i \in \text{ Dom}(g):g(i) \in \dbcu_{\beta <
\alpha} M_\beta\} \in J_g$ \nl

$\qquad (\beta) \qquad t,g \in M_\alpha$ hence Rang$(g) \subseteq M_\alpha$
\sn
\itemitem{ $(v)$ }  for any pregiven $\sigma = \text{ cf}(\sigma) \le \mu$ we
can demand $M_\alpha = \dbcu_{\varepsilon < \sigma} M_{\alpha,\varepsilon}$
where $\langle M_{\alpha,\varepsilon}:\varepsilon < \sigma \rangle$ is
increasing, $\mu \cup \{Y,\lambda,\mu,M,*\} \subseteq M_{\alpha,\varepsilon},
M_\alpha = \dbcu_{\varepsilon < \sigma} M_{\alpha,\varepsilon},
\langle M_{\alpha,\zeta}:
\zeta \le \varepsilon \rangle \in M_{\alpha,\varepsilon +1}$
and $\langle M_\beta:\beta < \alpha \rangle \in M_{\alpha,\varepsilon}$.
\endroster}
\endroster
\endproclaim
\bn
Before proving \scite{3.f3a} we deduce \scite{3.f3} and prepare the ground.
\demo{Proof of \scite{3.f3} from \scite{3.f3a}}  Straight: just let
$T = {\Cal P}$ and for $A \in T$ we let 

$$
\align
{\Cal F}_A = \biggl\{f:&\text{for some } (I,J) \in {\Cal Y},\text{Dom}(f)
\in J^+, \\
  &f \text{ is one to one into } A \text{ and} \\
  &{\Cal U}' \subseteq \text{ Dom}(f) \and {\Cal U}' \in I^+ \Rightarrow
c \ell(\text{Rang}(f \restriction {\Cal U}')) \in {\Cal P} \biggr\}.
\endalign
$$
\mn
We define the partial order $\le_T$ on $T$ by: $A_1 \le_T A_2$ iff
$A_2 \subseteq A_1$.
We have to check that the assumptions in \scite{3.f3a} holds, now clauses
$(a)', (b)', (c)'$ are the same as $(a), (b), (c)$ of \scite{3.f3}, and
clauses $(d)', (e)'$ are obvious.  As for clause $(f)'$ we shall use 
clause $(f)$ and the definition of ${\Cal F}_A$. \nl
[Why?  Let $t \in T,f \in {\Cal F}_t$ for $t = A \in {\Cal P}$, let
$(I_A,J_A) \in {\Cal Y},{\Cal U} \in J^+$ and $f^*$ be as in clause (f) of
\scite{3.f3} (with $f^*$ here standing for $f$ there), and let $\langle
{\Cal U}'_\alpha:\alpha < \alpha^* \rangle$ be as in subclause $(\beta)$
there.  For each $\alpha < \alpha^*$ choose $A'_\alpha$ as in subclause
$(\alpha)$ of clause (f) of \scite{3.f3} for ${\Cal U}' = {\Cal U}'_\alpha$.
Let us choose $T[f] =: \{A'_\alpha:\alpha < \alpha^*\}$, so $T[f] \in
[{\Cal P}]^{\le \mu} = [T]^{\le \mu}$.  Let us check that $T[f]$ is as
required; being a cover: let $r \in T_{\langle f \rangle}$, i.e. let $A' =r$
and (by the definition of $T_{\langle f \rangle}$), there is $g \in
{\Cal F}_r$ such that ${\Cal U}' = \{i:i \in \text{ Dom}(f)$ and $i \in
\text{ Dom}(g)$ and $f(i) = g(i)\} \in I^+_f$, so $r = A'' \in {\Cal P}$ and
${\Cal U}' \in I$ so for some $\alpha < \alpha^*$ we have $c \ell\{f(i):i \in
{\Cal U}'_\alpha\} \subseteq c \ell\{f(i):i \in {\Cal U}'\}$ so by the choice
of $A'_\alpha$ we have $A'_\alpha \subseteq c \ell\{f(i):i \in {\Cal U}'\}$
and $A'_\alpha \in T[f]$, let $s = A'_\alpha$, so $r \in T[f]$ is enough.
\nl
Now $s \in T_{\langle f \rangle}$ and (by the choice of $s = A'_\alpha$)
clearly $A'_\alpha \subseteq A$ which means $r \le s$.  Also $g' = f
\restriction \{i \in \text{ Dom}(f):f(i) \in A'_\alpha\}$ belongs to
${\Cal F}_{A'_\alpha}$ (as $f \in {\Cal F}_A$) and so $g'$ witness $s \in
T[f]$, proving $T[f]$ covers.  Lastly, $T[f]$ has cardinality $\le |\alpha^*|
\le \mu$.] \nl

Lastly let $\chi$ be large enough and $x \in 
{\Cal H}(\chi)$.  So by \scite{3.f3a} there is $\langle M_i:i < \lambda 
\rangle$ for our $\langle {\Cal F}_A:A \in T \rangle,x,\chi$ as required 
there.  It is enough 
to show that $\langle M_i:i < \lambda \rangle$ is as required in the 
conclusion of \scite{3.f3}.  Now clauses (i), (ii), (iii) of the 
conclusion of \scite{3.f3}
are just like clauses (i), (ii), (iii) and (v) of the conclusion of 
\scite{3.f3a}, so
we should check only clause (iv).  So assume $A \in {\Cal P}$ and $\alpha(A)$
is as defined there.  By clause (iv) of the conclusion of \scite{3.f3a} 
applied to $s=A$ there are $t,\alpha,g$ as there, i.e. $s \le_T t,\alpha <
\lambda,g \in {\Cal F}_t$ and $\{i \in \text{ Dom}(g):g(i) \in \dbcu_{\beta
< \alpha} M_\beta\} \in J_g$ and $t,g \in M_\alpha$.  So $t \in P,t \subseteq
A$, Dom$(g) \subseteq t \subseteq s = A$, 
Dom$(g) \subseteq M_\alpha$ and Rang$(g) \in M_\alpha$.  So 
Dom$(g) \in J^+_g$ and $\alpha$ is as required. 
\hfill$\square_{\scite{3.f3a}}$
\enddemo
\bigskip

\proclaim{\stag{3.f3b} Claim}  1) In \scite{3.f3} we can conclude
$(\alpha)_\theta \Rightarrow (\beta)_\theta$ when
\mr
\item "{$(\alpha)_\sigma$}"  if ${\Cal P}' \subseteq {\Cal P}$ has
cardinality $\le \mu$, \ub{then} we can find $h:\cup\{A:A \in {\Cal P}'\}$
to $\theta$ such that $A \in {\Cal P}' \Rightarrow \theta = \text{ Rang}
(h \restriction A)$
\sn
\item "{$(\beta)_\sigma$}"  we can find $h:\lambda \rightarrow \theta$ such
that $A \in {\Cal P} \Rightarrow \theta = \text{ Rang}(h \restriction A)$.
\ermn
2) In \scite{3.f3a} we can conclude $(\alpha)_\theta \Rightarrow
(\beta)_\theta$ when
\mr
\item "{$(\alpha)_\theta$}"  if $T' \subseteq T,|T'| \le \mu$ and $G$ is
a function with domain $\cup\{{\Cal F}_t:(\exists s \in T')(s \le_T t)\}$
such that $G(f) \in J_f$, \ub{then} we can find a function $h$ and 
$\langle (t_s,f_s):s \in T' \rangle$ such that 
$s \le_T t_s,f_s \in {\Cal F}_{t_s}$ and 
$s \in T_s \Rightarrow \theta = \{(h(f_s(i)):i \in \text{Dom}(f_s)
\backslash G(f_s)\}$
\sn
\item "{$(\beta)_\theta$}"  we can find a function $h:\lambda \rightarrow
\theta$ as in $(\alpha)_\theta$ for $T' = T$.
\endroster
\endproclaim
\bigskip

\demo{Proof}  1) Let $\{A^*_{\alpha,\zeta}:\zeta < \zeta_\alpha \le \mu\}$
list $\{A \in {\Cal P}:\alpha(A) = \alpha\}$ and let $(I^\alpha_\zeta,
J^\alpha_\zeta),{\Cal U}^\alpha_\zeta,f^\alpha_\zeta$ witness $\alpha(A_i)
= \alpha$.  Let $A'_{\alpha,\zeta} = \{f^\alpha_\zeta(i):i \in {\Cal U}_\zeta$
and $f_\zeta(i) \in M_\alpha \backslash \dbcu_{\beta < \alpha} M_\beta\}$.
Clearly $A'_{\alpha,\zeta} \in {\Cal P}$ and we apply 
clause $(\alpha)_\theta$ to
${\Cal P}_\zeta = \{A'_{\alpha,\zeta}:\zeta < \zeta_\alpha\}$ getting
$h_\alpha:\dbcu_{\zeta < \zeta_\alpha}A'_{\alpha,\zeta} \rightarrow
\theta$ so \wilog \, $h_\alpha:\lambda \cap M_\alpha \backslash
\dbcu_{\beta < \alpha} M_\beta \rightarrow \theta$.  Now $h = \dbcu_{\alpha <
\lambda} h_\alpha$ is as required. \nl
2) Similar.  \hfill$\square_{\scite{3.f3b}}$
\enddemo
\bn
The following is close to \cite[\S3]{Sh:161} (or see \cite[\S3]{Sh:523} or
\cite{EM}).
\definition{\stag{3.f4} Definition}  1) We say $\Gamma = (S,\bar \lambda)$ is
a full $(\lambda,\mu)$-set if:
\mr
\widestnumber\item{$(f)(\alpha)$}
\item "{$(a)$}"  $S$ is a set of finite sequences of ordinals
\sn
\item "{$(b)$}"  $S$ is closed under initial segments
\sn
\item "{$(c)$}"  $\bar \lambda = \langle \lambda_\eta:\eta \in S \rangle,
\lambda_{<>} = \lambda$
\sn
\item "{$(d)$}"  for each $\eta \in S,\{\alpha:\eta \char 94 \langle \alpha
\rangle \in S\}$ is empty or the regular cf$(\lambda_\eta)$
\sn
\item "{$(e)$}"  $\lambda_\eta > \mu$ \ub{iff} $\lambda_\eta \ne \mu$ 
\ub{iff} $(\exists \alpha)(\eta \char 94 \langle \alpha \rangle \in S)$
\ub{iff} $\eta \in S \backslash S^{mx}$
\sn
\item "{$(f)$}"  if 
$\eta \in S,\lambda_\eta \ne \mu$ \ub{then} for every ordinal 
$\alpha$ we have $\alpha < \text{ cf}(\lambda_\eta) \Leftrightarrow \eta
\char 94 \langle \alpha \rangle \in S$
\sn
\item "{$(g)$}" $(\alpha) \quad$  if 
$\lambda_\eta > \mu$ is a successor cardinal then
$\alpha < \lambda_\eta \Rightarrow \lambda^+_{\eta \char 94 \langle 
\alpha \rangle} = \lambda_\eta$
\sn
\item "{${{}}$}" $(\beta) \quad$ if $\lambda_\eta > \mu$ is a 
limit cardinal then $\langle \lambda_{\eta \char 94 \langle \alpha \rangle}:
\alpha < \text{ cf}(\lambda_\eta) \rangle$ \nl

$\qquad$ is strictly increasing with limit $\lambda_\eta$.
\ermn
2) Let $S^{mx} = \{\eta \in S:\lambda_\eta = \mu\}$.
\enddefinition
\bn
\ub{\stag{3.f5} Observation/Definition}:  If 
$\Gamma = (S,\bar \lambda)$ is a full $(\lambda,\mu)$-set, \ub{then} 
from $S$ we can reconstruct $\bar \lambda$ 
hence $\Gamma$, so we may say ``$S$ is a full $(\lambda,\mu)$-set, $\bar 
\lambda = \bar \lambda^{[S]}$.
\bn
\ub{\stag{3.f6} Fact/Definition}:  1) If $S$ is a full $(\lambda,\mu)$-set and
$\eta \in S$ let $S^{<\eta>} = \{\nu:\eta \char 94 \nu \in S\}$, is a full
$(\lambda_\eta,\mu)$-set. \nl
2) If for $\alpha < \text{ cf}(\lambda),S_\alpha$ is a full $(\lambda_\alpha,
\mu)$ set and $(\forall \alpha < \text{ cf}(\lambda))(\lambda_\alpha =
\lambda_0 \and \lambda = \lambda^+_0)$ or 
$\langle \lambda_\alpha:\alpha < \text{ cf}(\lambda) \rangle$
is strictly increasing with limit $\lambda,\lambda_0 \ge \mu$, \ub{then} 
$S = \{<>\} \cup
\dbcu_{\alpha < \text{ cf}(\lambda)} \{ \langle \alpha \rangle \char 94
\eta:\eta \in S_\alpha\}$ is a full $(\lambda,\mu)$-set. \nl
3) For a full $(\lambda,\mu)$-set $S$ and $\eta \in S$, if 
$\lambda_\eta > \mu$ let $\eta^+ = \langle \eta(\ell):\ell < k \rangle 
\char 94 \langle \eta(k)+1 \rangle$ if $\ell g(\eta) = k+1,<>^+$ will
be used though not defined.
\bigskip

\demo{Proof}  Straightforward.
\enddemo
\bigskip

\definition{\stag{3.f7} Definition}  1) We define by induction on $\lambda$ the
following.  For a set $X$ of cardinality $\lambda,\chi$ large enough and
$x \in {\Cal H}(\chi)$ we say $\bar N$ is a $\mu$-decomposition of $X$ for
${\Cal H}(\chi),x$ (or $(\lambda,\mu)$-decomposition) \ub{if} for some full
$(\lambda,\mu)$-set $S$ it is an $S$-decomposition of $X$ inside ${\Cal H}
(\chi)$, which means:
\mr
\item "{$(a)$}"  $\bar N = \langle (N_\eta,N^+_\eta):\eta \in S \rangle$
\sn
\item "{$(b)$}"  $N_\eta \prec N^+_\eta \prec ({\Cal H}(\chi),\in,<^*)$
except that $N_\eta$ is an empty set if Rang$(\eta) \subseteq \{0\}$
\sn
\item "{$(c)$}"  $\{X,x\} \in N^+_\eta$ and $\ell < \ell g(\eta) \Rightarrow
N_{\eta \restriction \ell},N^+_{\eta \restriction \ell} \in N^+_\eta$
\sn
\item "{$(d)$}"  $\|N^+_\eta\| = \lambda_{\eta^+} = \|(N^+_\eta \backslash
N_\eta) \cap X\|$ and $\lambda_{\eta^+} \subseteq N^+_\eta$
\sn
\item "{$(e)$}"  if $\lambda_{<>} > \mu$, \ub{then} 
$\langle N_{< \alpha >}:\alpha < \text{ cf}(\lambda_{<>})
\rangle$ is $\prec$-increasingly continuous with union containing 
$N^+_{<>}$
\sn
\item "{$(f)$}"  $N^+_{<\alpha>} = N_{<\alpha +1 >}$
\sn
\item "{$(g)$}"  for each $\alpha < \text{ cf}(\lambda_{<>}(S))$ the sequence
$\langle (N_{<\alpha> \char 94 \eta},N^+_{< \alpha> \char 94 \eta}):\eta \in
S^{<\alpha>} \rangle$ is a $(\lambda_\eta,\mu)$-decomposition of 
$X \cap N^+_{< \alpha>}$ for ${\Cal H}(\chi),
\langle x,N_\alpha,N^+_\alpha \rangle$.
\ermn
2) We say $\bar N$ is a $(\lambda,\mu,\sigma)$-decomposition of $X$ for
${\Cal H}(\chi),x$ if $\sigma = \text{ cf}(\sigma) \le \mu$ and in addition
\mr
\item "{$(h)$}"  for each $\eta \in S \backslash S^{\text{max}}$ the
sequence $\langle N_{\eta,\varepsilon}:\varepsilon \le \sigma \rangle$ 
is increasing continuous,
$N_{\eta,0} = N_\eta,\langle N_{\eta,\zeta}:\zeta \le \varepsilon
\rangle \in N_{\eta,\varepsilon +1}$ and the objects we demand $\in N^+_\eta$
belongs to $N_{\eta,1}$ (in clauses (c) and (h)).
\endroster
\enddefinition
\bigskip

\definition{\stag{3.f9} Definition}  1) Let $X,\lambda,\mu,{\Cal Y},
\bar{\Cal F}$ be as in \scite{3.f3a} so $T = \text{ Dom}(\bar{\Cal F})$.  \nl
We say $\bar N$ is a full $\mu$-decomposition of $X$ for $x,\chi$ is good
for $(X,{\Cal Y},\bar{\Cal F})$ if:
\mr
\item "{$(a)$}"  $\bar N$ is a full $(\lambda,\mu)$-decomposition of $X$ for
${\Cal H}(\chi),\langle x,X,\lambda,\mu,{\Cal Y},{\Cal P} \rangle$; let \nl
$\bar N =
\langle N_\eta:\eta \in S \rangle$ and $\bar \lambda = \bar \lambda^{[S]}$
\sn
\item "{$(b)$}"  if $s \in T$, \ub{then} for some $t \in T,s \le_T t$ and
for some $\eta \in S^{mx}$ (i.e. $\lambda_\eta = \mu$) there is 
$f \in {\Cal F}_t$ and so $(I_f,J_f) \in {\Cal Y},{\Cal U}_f \in J^+_f,
f:{\Cal U}_f \rightarrow \text{ Rang}(g)$ such that:
{\roster
\itemitem{ $(*)_1$ }  $\{i \in {\Cal U}_f:
f(i) \in \cup\{N_\nu:\nu \le_{\ell x} \eta \text{ and }
\nu \in S^{mx}\}\}$ belongs to $J_f$
\sn
\itemitem{ $(*)_2$ }  $\{i \in {\Cal U}_f:f(i) \in N^+_\eta \backslash \cup 
\{N_\nu:\nu <_{\ell x} \eta \text{ and } \nu \in S^{mx}\}\}$ 
belongs to $J^+_f$
\sn
\itemitem{ $(*)_3$ }  $t,f$ belong to $N^+_\eta$.
\endroster}
\ermn
2) We define a full $(\lambda,\mu,\sigma)$-decomposition similarly.
\enddefinition
\bigskip

\proclaim{\stag{3.f10} Claim}  Under the assumption of \scite{3.f3a}, for
$x \in {\Cal H}(\chi),\sigma = \text{ cf}(\sigma) \le \mu$ and $\chi$ 
large enough there is a $(\lambda,\mu,\sigma)$-decomposition of 
$X$ for $\chi,x$ good for $(X,{\Cal Y},\bar{\Cal F})$.
\endproclaim
\bigskip

\demo{Proof}  By induction on $\lambda = |X|$.
\enddemo
\bn
\ub{Case 1}:  $\lambda = \mu$.

Trivial.
\bn
\ub{Case 2}:  $\lambda = \text{ cf}(\lambda) > \mu$.

Choose $\langle N_\alpha:\alpha < \text{ cf}(\lambda) \rangle$ such that
$\{x,X,\bar{\Cal F},\mu,\lambda\} \in N_0,N_\alpha \prec ({\Cal H}(\chi),
\in$, \nl
$<^*_\chi),N_\alpha$ is $\prec$-increasingly continuous, $\langle N_\beta:
\beta \le \alpha \rangle \in N_{\alpha +1}$, each $N_\alpha$ has cardinality
$< \lambda$ and $N_\alpha \cap \lambda$ is an
initial segment.  For $t \in T$ let $\alpha(t) = \text{ Min}\{\alpha:
\text{for some } f \in \dbcu_{s \ge t} {\Cal F}_s$ and $(I_f,J_f) \in 
{\Cal Y}$ (as in \scite{3.f3a} clause $(e)'$) we have $\{i:i \in
\text{ Dom}(f)$ and $f(i) \in N_\alpha\} \in J^+\}$. \nl
Let $S = \{\beta < \lambda:\text{ for some } t \in T \text{ we have }
\beta = \alpha(t)\} \subseteq \lambda$.  For each $\beta \in S$ choose
$t_\beta \in T$ and $s_\beta,t_\beta \le_T s_\beta$ such that 
$\beta = \alpha(t_\beta)$ and $f_\beta \in {\Cal F}_{s_\beta}$ witness 
this.  Let ${\Cal U}_\beta = \text{ Dom}(f_\beta)$ and let 
$(I_\beta,J_\beta) = (I_{f_\beta},J_{f_\beta})$.  Now \wilog \, 
$f_\beta \in N_{\beta +1}$ (hence $s_\beta,I_\beta,J_\beta \in 
N_{\beta +1}$ (as all the requirements on
$f_\beta$ have parameters in $N_{\beta +1}$).  
First assume toward contradiction that $S$ is stationary. 
Now as ${\Cal Y} \in N_0,|{\Cal Y}| < \lambda$ clearly ${\Cal Y} \subseteq
N_0$ hence for some $y \in {\Cal Y}$ the set 
$S_y = \{\beta \in S:(I_\beta,J_\beta) = y\}$ is
stationary.  Let $y = (I^*,J^*)$ and $S'_y = \{\beta \in S_y:N_\beta \cap
\lambda = \beta\}$, clearly it is stationary.  It suffices to show that
$\langle f_\delta:\delta \in S'_y \rangle$ exemplifies
NFr$_2(\lambda,{\Cal Y})$ contradicting assumption $(c)'$ from \scite{3.f3a}.
If not, for some
$\delta_1 < \delta_2$ in $S_\beta$ we have 
$B =: \{i:i \in \text{ Dom}(f_{\delta_1}),
i \in \text{ Dom}(f_{\delta_2})$ and $f_{\delta_1}(i) = f_{\delta_2}(i)\}
\in I^+$, hence $t_{\delta_2} \in T_{\langle f_{\delta_1} \rangle}$ (see
\scite{3.f3a}, clause $(f)'$) hence by an assumption there is $t'_{\delta_2}$
such that $t_{\delta_2} \le_T 
t'_{\delta_2} \in T[f_{\delta_1}]$.  But $\bar{\Cal F},f_{\delta_1}$ belong
to $N_{\delta_1+1} \prec N_{\delta_2}$ hence $T_{\langle f_{\delta_1} \rangle}
\in N_{\delta_1 +1}$ but $T[f_{\delta_1}]$ has cardinality $\le \mu$ (see
clause $(f)'$ of \scite{3.f3a}) hence $T_{[f_{\delta_1}]} \subseteq
N_{\delta_1 +1}$ but $t'_{\delta_2} \in T[f_{\delta_1}]$ so $t'_{\delta_2} \in
N_{\delta_1 +1}$ hence ${\Cal F}_{t_{\delta_2}} \in N_{\delta_1 +1}$ hence
(see \scite{3.f3a}, $(d)'$) we have 
$F_{t_{\delta_2}} \subseteq N_{\delta_1 +1}$ hence
there is $f' \in {\Cal F}_{t'_{\delta_2}} \cap N_{\delta_1 +1}$ hence
Rang$(f') \subseteq N_{\delta_1 +1}$ contradicting the demand 
$\alpha(t_{\delta_2}) = \delta_2$.  So $S$ is not stationary. 
\sn
Let $E$ be a club of $\lambda$ disjoint to $S$ and we can find $\bar N' =
\langle N'_\alpha:\alpha < \lambda \rangle$ like $\langle N_\alpha:\alpha <
\lambda \rangle$ such that $E,\bar N \in N_0$ so for $\bar N',S = \emptyset$.
Now for each $\alpha$ we use the induction hypothesis on $X_\alpha = X \cap 
N'_{\alpha +1} \backslash N'_\alpha$ and $\langle \bar{\Cal F}^{\langle \alpha
\rangle}_t:t \in T^{\langle \alpha \rangle} \rangle$ where 
$T^{\langle \alpha \rangle} = T \cap N'_{\alpha +1} \backslash N'_\alpha$
and ${\Cal F}^{\langle \alpha \rangle}
_t = \{f \restriction {\Cal U}:{\Cal U}$ is $\{i \in \text{ Dom}(f):f(i) 
\in X_\alpha\}$ and $f \in {\Cal F}_t \cap N'_{\alpha +1}\}$.
\bn
\ub{Case 3}:  $\lambda$ singular $> \mu$.

Let $\lambda = \dsize \sum_{i < \text{ cf}(\lambda)} \lambda_i,
\langle \lambda_i:i < \text{ cf}(\lambda) \rangle$ increasingly continuous,
$\lambda_0 > \mu^+$.  We choose by induction on $\zeta < \mu^+,\langle
N^\zeta_i:i < \text{ cf}(\lambda) \rangle$ such that:
\mr
\item "{$(a)$}"  $N^\zeta_i$ is $\prec$-increasing in $i$
\sn
\item "{$(b)$}"  $\langle \lambda_i:i < \text{ cf}(\lambda) \rangle,X,
\lambda,\mu,\bar{\Cal F}$ all belong to $N^\zeta_0$
\sn
\item "{$(c)$}"  $\lambda_i \subseteq N^\zeta_i$ and $\|N^\zeta_i\| = 
\lambda_i$
\sn
\item "{$(d)$}"  for each $i, \langle N^\zeta_i:\zeta \le \mu^+ \rangle$ is
$\prec$-increasingly continuous
\sn
\item "{$(e)$}"  $\left< \langle N^\varepsilon_i:i < \text{ cf}(\lambda)
\rangle:\varepsilon \le \zeta \right> \in N^{\zeta +1}_i$.
\ermn
For each $i < \lambda$ and $\zeta < \mu^+$ and $(I,J) \in {\Cal Y}$ let
${\Cal F}^{\zeta,i}_{(I,J)}$ be a maximal family of functions $f \in \{f
\restriction {\Cal U}:{\Cal U} \in J^+_f,f \in \dbcu_{t \in T} {\Cal F}_t,
{\Cal U} \subseteq \text{ Dom}(f)\}$, Rang$(f) \subseteq X \cap N^\zeta_i$ 
and $f \ne g \in {\Cal F}^{\zeta,i}_{(I,J)} \Rightarrow 
\{i:i \in \text{ Dom}(f),i \in
\text{ Dom}(g)$ and $f(i) \ne g(i)\} \in I$.  Without loss of generality
${\Cal F}^{\zeta,i}_{(I,J)} \in N^\zeta_{i+1}$ and by \scite{3.f2}(4) and
assumption \scite{3.f3a} clause $(c)'$ we know
$|{\Cal F}^{\zeta,i}_{(I,J)}| \le \lambda_i$, so a list of it of length
$\le \lambda_i$ belongs to $N^{\zeta +1}_i$ hence ${\Cal F}^{\zeta,i}
_{(I,J)} \subseteq N^{\zeta +1}_i$.  So if $t \in {\Cal T}$ and we define
$\alpha(t)$ as in Case 2 for $\langle N^{\mu^+}_\alpha:\alpha \le
\text{ cf}(\mu) \rangle$, we get that $\alpha(t)$ is necessarily nonlimit.
\hfill$\square_{\scite{3.f10}}$
\bigskip

\demo{Proof of \scite{3.f3a}}

Just by \scite{3.f10} above and \scite{3.f3b} 
(reading Definition \scite{3.f9}).
\enddemo
\bigskip

\demo{Proof of Lemma 3.2}

Just by \scite{3.f10} above and \scite{3.f3b}.
\enddemo
\bigskip

\demo{Proof of \scite{3.5}}  We use \scite{3.f3} above. \nl
1) Without loss of generality let $\lambda$ be 
the set of points of $X,\mu = \chi,I = \{A \subseteq {}^\omega 2:
\text{the closure of } A \text{ is countable}\},J$ the following ideal on
${}^\omega 2$

$$
\{{\Cal U} \subseteq {}^\omega 2:|{\Cal U}| < 2^{\aleph_0}\}
$$
\mn
and

$$
{\Cal Y} = \{(I,J)\}.
$$
\mn
So the conclusion $(*)$ of \scite{3.5} just means ``for some regular
$\kappa \in (\mu,\lambda]$ we have NFr$_1(\kappa,{\Cal Y})$" and toward
contradiction assume it fails.  Clearly $\chi \ge \mu$.  Also without loss
of generality the set of points of $X$ is $\lambda$, let $c \ell:
[\lambda]^{\le \mu} \rightarrow [\lambda]^{\le \mu}$ be

$$
\align
c \ell(A) = \bigl\{ \alpha:&\alpha \in A \text{ or for some countable }
B \subseteq A,\alpha \text{ belongs} \\
  &\text{to the closure of } B \text{ in the topological space } X
\text{ and} \\
  &c \ell(B) \text{ has cardinality } \le 2^{\aleph_0} \bigr\}.
\endalign
$$
\mn
Let us consider the assumptions of \scite{3.f3}.  Now clause (a) is by the
explicit choice of ${\Cal Y}$ above, also (b).  Clause (c) is the assumption
toward contradiction above, clause (d) (on $c \ell$) holds as clearly
$A \in [\lambda]^{\le \mu}$ implies $c \ell(A) = \cup \{c \ell(B):B \in
[A]^{\le \aleph_0}\}$ and $[A]^{\le \aleph_0}$ has cardinality
$\le \mu^{\aleph_0} = \chi^{\aleph_0} = \chi$ and for each countable $B$,
contribute at most $2^{\aleph_0}$ points.  Clause (e) holds by the properties
of closure.  Lastly, for clause (f) we define

$$
\align
{\Cal P} = \{A:&A \subseteq \lambda \text{ is a closed subset of } A, \\
  &\text{has cardinality continuum and } X \restriction A \text{ is
homeomorphic to } {}^\omega 2 \}.
\endalign
$$
\mn
So as all the assumptions of \scite{3.f3} holds so we can apply \scite{3.f3b}.
There for $\theta = 2^{\aleph_0}$, if we can apply \scite{3.f3b}(1)
we get $X \nrightarrow ({}^\omega 2)^1_{2^{\aleph_0}}$.  
But $(\alpha)_\theta$ of \scite{3.f3b} is immediate
for any closed subspace $Y$ of $X$ homeomorphic to
${}^\omega 2$, we have $|Y \cap \dbcu_{\beta < \alpha} M_\beta| < 
2^{\aleph_0}$ and for some $A \in M_\alpha,|A| = 2^{\aleph_0},Y \cap 
c \ell(A)$ contains a set $Y'$ homeomorphic to ${}^\omega 2$, and this has
$2^{\aleph_0}$ pairwise disjoint subspaces which $\in M_\alpha$ so at least
one is disjoint to $\dbcu_{\beta < \alpha} M_\beta$ so we are done by the
choice of $h_\alpha$.]  \nl
2) Similar, we just should be more accurate about closure; note that the
topological closure of a countable set may have cardinality bigger than
$2^{\aleph_0}$. For $A \subseteq
X$ let $c \ell(A) = c \ell(A,X) = \cup \{\text{Rang}(f):f$ a one to one
mapping from $\Bbb R$ to $X$ which is a homeomorphism and such that
$Y_f = \{x \in \Bbb R:f(x) \in A\}$ is dense$\}$.  But for any such $f_1,f_2$,
if some $Y \subseteq Y_{f_1} \cap Y_{f_2}$ is countable dense and
$[x \in Y \Rightarrow f_1(y) = f_2(y)]$ then $f_1 = f_2$, so the proof is
similar.   \hfill$\square_{\scite{3.5}}$ 
\enddemo
\bigskip

\remark{\stag{3.10} Concluding Remark}  1) Of course, we may replace in
\scite{3.1} the space ${}^\omega 2$ by many others, e.g. $\Bbb R$, or any
Hausdorff $Y^*$ with $2^{\aleph_0}$ points such that for any uncountable
$A \subseteq Y^*$, for some countable $B \subseteq A,|c \ell_{Y^*}(B)| =
2^{\aleph_0}$ moreover if $Z \subseteq Y^*,|Z| < 2^{\aleph_0}$ for some
uncountable $B' \subseteq c \ell_{Y^*}(B)$ we have $c \ell_{Y^*}(B')$ is
disjoint to $Z$.

We can also add variants with $\rightarrow_w$ replacing $\rightarrow$.  As
long as the space has $\le 2^{\aleph_0}$ points, the only place we should be
concerned is the proof of \scite{3.5}, we reconsider the choice of $c \ell$
in the proof.  In all cases for an embedding $f$ from $Y \subseteq Y^*$ to 
$X$, let $c \ell(\text{Rang}(f)) = \{x \in X:
\text{for some } y \in Y^*,f \cup \{ \langle y,x \rangle\}$ is an embedding
of $Y^* \restriction ({\Cal U} \cup \{x\})$ to $X \restriction
((\text{Rang}(f)) \cup \{y\})\}$ and $f^+ = f \cup \{ \langle y,x \rangle:
x,y \text{ as above}\}$.
The point is that for this choice of $c \ell$, if $Y_1 \subseteq Y_2
\subseteq Y^*,Y_2 \subseteq c \ell_{Y^*}(X_1)$ if $f$ embeds $Y_2$ into $X$
with Rang$(f)$ not necessarily close, then $(f \restriction X_1)^+$ is a
function from some $Y_3 \subseteq Y^*$ into $X$ extending $f$. \nl
2) We may like to add to \scite{3.1} the case with continuum many colours
that is let $(B_m)_{< \mu}[{}^\omega 2]$ and $(B_m)^+_{< \mu}$ be defined
like $(B)[{}^\omega 2],(B)^+$, replacing $)^1_{< \text{ cf}(2^{\aleph_0})}$
by $)^1_{< \mu}$ and we add $(B_m)_{< \beth^+_2}[{}^\omega 2],
(B)^+_{< \beth^+_2}$ to the list of equivalent statements.  Similarly for
$(A)$.  More is proved $X \rightarrow ({}^\omega 2)^1_{< \lambda}$ where
$X$ has $\lambda$ points (or we get $\lambda$ when we ask for compact $X$).
The main point is adopting \scite{t.2} (and \scite{t.4}).
\sn
For this we add also $(C_m)_{\beth_2,\beth_2,\aleph_2}$ where for 
$\kappa \ge \theta \ge \sigma$ we let \nl

$(C_m)_{\kappa,\theta,\sigma} \quad$  there are $\lambda,S,\bar f$ 
such that
\mr
\item"{$(a)$}"  $S \subseteq \lambda$ is stationary 
$> \kappa^+,\kappa > \theta \ge \sigma$
\sn
\item"{$(b)$}"  $\bar f = \langle f_\delta:\delta \in S \rangle$
\sn
\item"{$(c)$}"  Dom$(f_\delta) = \theta$, each $f_\delta(i)$ is a subset
of $\delta \backslash i$ of cardinality $\le \kappa$ and $\langle \min
(f_\delta(i)):i < \theta \rangle$ is increasing with limit $\delta$ (can
ask $i < j < \theta \Rightarrow f (\sup(f_\delta(i)) < \min(f_\delta(j))$
\sn
\item"{$(d)$}"  if $\delta_1 < \delta_2$ are in $S$ then $\{i < \theta:
f_{\delta_2}(i) \cap \dbcu_{j < \theta} f_{\delta_1}(j) \ne \emptyset\}$ has
cardinality $< \sigma$
\sn
\item"{$(e)$}"  if $F_\ell:\lambda \rightarrow [\lambda]^{\le \kappa}$ for
$\ell = 0,1$ and $F_0(\alpha) \in [\lambda \backslash \alpha]
^{\le \kappa}$, \ub{then} for some $\delta \in S$ we have:
{\roster
\itemitem{ $(\alpha)$ }  $f_\delta$ is $(F_0,F_1)$-free which means: \nl
for $i \ne j < \theta$, the set $F_1(f_\delta(i))$ is disjoint to 
$F_0(f_\delta(j))$
\sn
\itemitem{ $(\beta)$ }  there are $\langle \alpha_i:i < \theta \rangle$
such that $f_\delta(i) = F_0(\alpha_i)$ and sup$[\dbcu_{j<i} f_\delta(i)]
< \alpha_i$. \nl
Similarly for (D).  Why is this O.K.?  See below, noting that we get more.
\endroster}
\ermn
3) As before, $(B_m)^+ \Rightarrow (B_m)[{}^\omega 2] \Rightarrow (A_m)
[{}^\omega 2]$ and $(B_m)^+ \Rightarrow (A_m)^+ \Rightarrow (A_m)
[{}^\omega 2]$, also easily $(C) \Rightarrow (C)^+_{\beth_2,\beth_2,\aleph_2};
(B_m)^+ \Rightarrow (B)^+,(A_m)^+ \Rightarrow (A)^+,(B_m)[{}^\omega 2]
\Rightarrow (B)[{}^\omega 2]$ and $(A_m)[{}^\omega 2] \Rightarrow (A)
[{}^\omega 2]$
\mr
\item "{$(f)$}"  if $(F_0,F_1)$ is a pair of functions with domain $\lambda$
and $F_0(i) \in [\lambda \backslash i]^{\le \kappa}$
\ermn
3A) The forcing in \scite{g.9a}, with the role of $A_\zeta$ being replaced by
$\dbcu_{i < \theta} f_\zeta(i)$ and $A^p_\zeta \subseteq \dbcu_{i < \theta}
f_\delta(i)$ such that $i < \theta \Rightarrow |A^p_\zeta \cap f_\delta(i)|
\le 1$ works. \nl
4) Also
\mr
\item "{$\boxtimes_4$}"  $(D)_{\beth_2,\beth_2,\aleph_2}$
implies the consistency of $(B_m)^+_{< \beth^+_2}$.
\ermn
As before \wilog \, for some $\kappa = \kappa^{< \kappa} \ge \theta = 
2^{\aleph_0},\sigma$ are such that $(C)_{\kappa,\theta,\sigma}$ hold.  Now we
just need to repeat the proof of \scite{t.2}.
The asymmetry in clause (d) does not hurt as if $\delta_2 \ne
\delta_2,A^{p_1}_{\delta_1},A^{p_2}_{\delta_2}$ are well defined, then it
follows that $|A^{p_1}_{\delta_1} \cap A^{p_2}_{\delta_2}| < \sigma$.
\sn
In the crucial point we let $p^* \Vdash ``\underset\tilde {}\to c:\lambda
\rightarrow \underset\tilde {}\to \mu$ for some $\underset\tilde {}\to \mu
< \lambda"$.  Really less is enough: let $p^* \Vdash ``
\underset\tilde {}\to Z \subseteq \lambda$ is unbounded" and we shall find
$q$ and $\delta \in S$ such that $p^* \le q \in P$ and $q \Vdash ``
{\underset\tilde {}\to X^*} \restriction A^p_\delta$ is a copy of the space
$Y$ (e.g. ${}^\omega 2$) and $A^p_\delta \subseteq Y"$.  How?  We define
\nl
$F_0(\alpha) = \{\beta:\beta \in [\alpha,\lambda) \text{ and }
p^* \nVdash \beta \ne \text{ Min}(\underset\tilde {}\to Z \backslash \alpha)
\}$. \nl
$F_1(\alpha) = \cup\{u^{p_{\alpha,i}}:i < \kappa\}$ where $\langle
p_{\alpha,i}:i < \kappa \rangle$ is a maximal antichain above $p^*$ such that
$p_{\alpha,i}$ forces $\alpha \in \underset\tilde {}\to Z$ or forces
$\alpha \notin \underset\tilde {}\to Z$. \nl
Now we repeat the proof of \scite{t.2}, but instead deciding the colour we
decide the right member of $\underset\tilde {}\to Z$. \nl
5) Lastly, we get $(C)^+_{\kappa,\theta,\nu}$ from $(C)_{\kappa,\theta,\nu}$.
So assume $\lambda > \kappa^+,\kappa > \theta \ge \sigma$ and $\langle
A_\delta:\delta \in S \rangle$ are as in $(C)$ and as before (by forcing)
\wilog \, $\diamondsuit_S$.  Now we can actually prove $(C)_{\kappa,\theta,
\sigma}$ for $\lambda$.  So we prove
\mr
\item "{$\boxtimes_5$}" if
{\roster
\itemitem{ $(\alpha)$ }  $\lambda > \kappa^+,\kappa > \theta \ge \sigma,
\kappa^\sigma < \lambda$
\sn
\itemitem{ $(\beta)$ }  $J$ an ideal on $\theta$ such that $(\forall A \in
J^+)(\exists a \in J^+)(a \subseteq A)$
\sn
\itemitem{ $(\gamma)$ }  $S \subseteq \lambda$ is stationary, $\bar f = 
\langle f_\delta:\delta \in S \rangle,f_\delta:\theta \rightarrow \theta$ 
increasing, $\delta_1 < \delta_2 \Rightarrow \{i < \theta:f_{\delta_1}(i) = 
f_{\delta_2}(i)\} \in J^+$ 
\sn
\itemitem{ $(\delta)$ }  $\diamondsuit_S$.
\endroster}
\ermn
Then $(C)_{\kappa,\theta,\sigma}$ as witnessed by $\lambda$. \nl
So let $\langle (F^\delta_0,F^\delta_1):\delta \in S \rangle$ be such that
$F^\delta_\ell:\delta \rightarrow [\delta]^{< \kappa}$ for $\ell=0,1$ be
such that: if $F_\ell:\lambda \rightarrow [\lambda]^{\le \kappa}$ for
$\ell=0,1$ then $S_{(F_0,F_1)} = \{\delta \in S:F_0 \restriction \delta = 
F^\delta_0 \text{ and } F_1 \restriction \delta = F^\delta_1\}$ is stationary.
We now choose by induction on $\delta \in S$ a function $f_\delta$ such 
that:
\mr
\item "{$(a)$}"   if there is a function $f$ with domain $\theta$ satisfying
the conditions below then $f_\delta$ is such a function, otherwise $f_\delta$
is constantly $\emptyset$
{\roster
\itemitem{ $(\alpha)$ }  $f(i) \in [\delta]^{\le \kappa} \backslash \{
\emptyset\}$
\sn
\itemitem{ $(\beta)$ }  $i < j \Rightarrow \sup(f_\delta(i)) < 
\min(f_\delta(j))$
\sn
\itemitem{ $(\gamma)$ }  for each $i < \theta$ for some $\alpha_i < \delta$
we have $F^\delta_0(\alpha_i) = f_\delta(i)$ and \nl
sup$(\dbcu_{j<i} f(j)] < \alpha_i \le \min f(i))$
\sn
\itemitem{ $(\delta)$ }  $\langle \min(f(i)):i < \theta \rangle$ converge
to $\delta$
\sn
\itemitem{ $(\varepsilon)$ }  for $i \ne j < \theta$ the set
$F^\delta_1(f(i))$ and $F^\delta_6(f(j))$ are disjoint
\sn
\itemitem{ $(\zeta)$ }  if $\delta_1 \in \delta \cap S$ then \nl
$\{i < \delta:f(i) \cap \dbcu_{j < \theta} f_{\delta_1}(j) \ne \emptyset\}$
has cardinality $< \sigma$.
\endroster}
\ermn
Let $S^- = \{\delta \in S:f_\delta \text{ is not constantly } \emptyset\}$
and we suffice to prove that $\bar f = \langle f_\delta:\delta \in S^-
\rangle$ is as required.  Most clauses hold by the definition and we should
check clause (e), so let $F_0,F_1$ be as there.  Let $S_{F_0,F_1} =
\{\delta \in S:F_0 \restriction \delta = F^\delta_0 \text{ and } F_1
\restriction \delta = F^\delta_1\}$, so this set is stationary. \nl
For every $\alpha \in S^* = \{\delta < \lambda:\text{cf}(\delta) =
\kappa^+\}$ let $g(\alpha) = \sup(\alpha \cap F_1(\alpha)) < \alpha$ so 
$g$ is constantly $\alpha(*)$ on some stationary
$S^{**} \subseteq S$.
\sn
$E_0 = \{\delta < \lambda:\text{otp}(S^{**} \cap \delta) = \delta$ and
$\alpha < \delta \Rightarrow \sup(F_0(\alpha)) < \delta$ and $\alpha <
\delta \Rightarrow \sup (F_1(\alpha)) < \delta\}$.
\sn
Let $E^*_1 = \{\delta < \lambda:\text{otp}(E_0 \cap \delta) = \delta\}$ and
for $\delta \in E_1 \cap S_{F_0,F_1}$ let $A'_\delta = \{\alpha \in E_0:
\text{otp}(\alpha \cap E_0) \in A_\delta\}$, so $A_\delta \subseteq \delta
= \sup(A_\delta)$, otp$(A_\delta) = \theta$ and $\delta_1 \ne \delta_2 \in
E_1 \cap S_{F_0,F_1} \Rightarrow |A_{\delta_1} \cap A_{\delta_2}| < \sigma$.
\sn
Let $A_\delta = \{\alpha'_{\delta,i}:i < \theta\}$ increasingly and let
$\alpha_{\delta,i} = \text{ Min}(S^{**} \backslash (\alpha'_{\delta,i} +1))$
so $\alpha_{\delta,i} < \alpha''_{\delta,i+1}$ (even
$\alpha_{\delta,i} < \text{ Min}(E_1 \backslash (\alpha'_{\delta,i} +1)$ and
choose $f'_\delta$ a function with domain $\theta$ by

$$
f'_\delta(i) = F_0(\alpha_{\delta,i}) = F^\delta_0(\alpha'_i)
$$
\mn
(the last equality as $F_\ell \restriction \delta = F^\delta_\ell$ as
$\delta \in S_{F_0,F_1}$). \nl
Clearly $f'_\delta(i) = F_0(\alpha_i) \subseteq \text{ Min}(E_1 \backslash
(\alpha'_\ell +1))$ and

$$
\gamma \in f'_\delta(i) \Rightarrow F(\gamma) \subseteq \text{ Min}
(E_1 \backslash (\alpha'_{\delta,i} +1)) \le \alpha'_{\delta,i+1} <
\alpha'_{\delta,i+1}
$$ 

$$
\gamma \in f'_\delta(i) \Rightarrow F(\gamma) \cap \alpha_{\delta,i} 
\subseteq \alpha(*) < \alpha_0
$$
\mn
Now $f'_\delta$ satisfies almost all the requirements on $f_\delta$ and if
$f'_\delta = f_\delta$ for stationarily many $\delta \in E_1 \cap 
S_{F_0,F_1}$ we are done.  Let $W = \{\delta \in E_1 \cap S_{F_0,F_1}:
f'_\delta \ne f_\delta\}$, we shall prove that $W$ is not stationary - this
is more than enough.
\sn
So for $\delta \in W$ necessarily for some $h(\delta) \in \delta \cap S$
we have

$$
w_\delta = \{i < \theta:f'_\delta(i) \cap \dbcu_{j < \theta} f_{(\delta)}(j)
\ne 0\}
$$
\mn
has cardinality $\ge \sigma$, so by Fodor's lemma for some $\delta(*)$ we
have $W_1 = \{\delta \in W:h(\delta) = \delta(*)\}$ is stationary. \nl
Similarly as $\theta^\sigma < \lambda = \text{ cf}(\lambda)$ for some
$w^* \in [\theta]^\sigma,w_2 = \{\delta \in w_1:w^* \subseteq w_\delta\}$ is
stationary.  As ${}^\sigma[\dbcu_{j < \theta} f_{\delta(*)}(j)]^\sigma$ has
cardinality $\kappa^\sigma$ which is $< \lambda$ \wilog \, for some
$h^*:w^* \rightarrow \dbcu_{j < \delta}f_{\delta(*)}(j)$ the set

$$
W_3 = \{\delta \in W_2:(\forall i \in w^*)(h^*(i) \in f'_\delta(i) \cap
\dbcu_{j < \theta} f_{\delta(*)}(j))\}
$$
\mn
is stationary.  So if $\delta_1 < \delta_2$ are in $w_3$ the set
$\{i < \theta:f'_{\delta_1}(i) = f'_{\delta_2}(i)\}$ include $w^*$.
But $f'_{\delta_1}(i) = f'_{\delta_2}(i)$ implies that $\alpha^{\delta_1}_i=
\alpha^{\delta_2}_i$, hence $A_{\delta_1} \cap A_{\delta_2}$ has cardinality
$\ge \sigma$ continuously. \nl
6) $W$ has a $\boxtimes$ clause $(\delta)$, we add: Rang$(f_\delta)$ is
bound in $\delta$? \nl
This is equivalent to: for some fixed $\mu < \lambda,(\forall \delta)
(\text{Rang}(f_\delta) \subseteq \mu)$.  Repeating the proof and replacing
club of $C \in [\mu]^\mu$ we get clause $(C)_{\kappa,\theta,\sigma}$
witnessing $\lambda$ with Rang$(f_\delta) \subseteq \mu$.  We then get
versions of the $(A)$'s and $(B)$'s with $\mu$ points. \nl
(Note one special point: we should rephrase the ``weak $\Delta$-system
argument, by using it on a tree with two levels. \nl
7) Note that by part (5) we get a stronger version of the topological
statements: for any $\lambda$ 
(or $\mu$ in (6)) points there is a close copy of
${}^\omega 2$ (or the space $Y$) included in it.  Of course, if we like the
space to be compact this refers only to any set of $\lambda$ (or $\mu$)
points among the original ones.  Note the Boolean Algebra of clopen sets
(when $Y$ has such a basis) satisfies the c.c.c. (remember in the cases only
$u^p_{\zeta,2i} \cap u^p_{\zeta,2i+1} = \emptyset$ is demanded, the Boolean
Algebra is free) so we cannot control the set of ultrafilters (= points), but
if we allow more disjointness demand we may, but we have not considered it.
\endremark
\bigskip

\proclaim{\stag{3.11} Claim}  Assume $\mu = \mu^{< \mu}$.  There is a
$\mu$-complete $\mu^+$-c.c. forcing notion $Q$ such that

$$
\align
\Vdash_Q ``&\text{ there is a function } h:{}^\mu \mu \rightarrow \mu
\text{ such that} \\
  &(\alpha) \quad \text{ if } C \in V \text{ is a closed subset of } 
{}^\mu \mu \text{ of cardinality } \le \mu \\
  &\qquad \text{ then } \alpha < \mu \Rightarrow |C \cap h^{-1}\{\alpha\}|
< \mu \\
  &(\beta) \quad \text{ if } A \in V \text{ is a subset of } {}^\mu \mu 
\text{ of cardinality } > \mu \\
  &\qquad \text{ then } \alpha < \mu \Rightarrow |A \cap h^{-1}\{\alpha\}|
= |A|".
\endalign
$$
\endproclaim
\bigskip

\demo{Proof}  As in the proof of \scite{3.6}, it suffices to prove:
\mr
\item "{$(*)$}"  assume $i^*,j^* < \mu$ and $\eta_{\alpha,i} \in {}^\mu \mu$
for $\alpha < \mu^+,i < i^*$ is with no repetitions and $C_{\alpha,j}
\subseteq {}^\mu \mu$ is closed with $\le \mu$ points for $\alpha < \mu^+,
j < j^*$.  Find $\alpha < \beta$ such that $i < i^* \and j < j^* \Rightarrow
\eta_{\alpha,i} \notin C_{\beta,j}$. \nl
\sn
Why $(*)$ holds?  Assume not.  First choose $\delta^* < \mu^+$ such that:
\sn
\item "{$(**)$}"  if $\beta < \mu^+$ and $\zeta < \mu$ then for some
$\alpha < \delta^*$ we have $i < i^* \Rightarrow \eta_{\alpha,i} \restriction
\zeta = \eta_{\beta,i} \restriction \zeta$.
\ermn
We can find $\beta$ such that $\delta^* < \beta < \mu^+$ and
$\{\eta_{\beta,i}:i < i^*\}$ is disjoint to $\dbcu_{j < j^*} C_{\delta^*,j}$.
$\beta$ exists as $|\dbcu_{j<j^*} C_{\delta^*,j}| \le \mu$.
Let $\zeta^* < \mu$ be large enough such that 
$i < i^* \and j < j^* \Rightarrow \neg
(\exists \nu)(\eta_{\beta,i} \restriction \zeta \triangleleft \nu \in
C_{\delta^*,j})$.  Lastly, choose $\alpha < \delta^*$ such that $i < i^*
\Rightarrow \eta_{\alpha,i} \restriction \zeta = \eta_{\beta,i} \restriction
\zeta$.  Now $\langle \alpha,\delta^* \rangle$ can serve as $(\alpha,\beta)$
above.  \hfill$\square_{\scite{3.11}}$
\enddemo
\newpage

\head {Appendix: similar proofs} \endhead  \resetall 
\bn
\ub{\stag{p.1} Proof of \scite{t.4}} \nl
\demo{Proof}  We write the proof for part (1) and indicate the changes for
part (2).  Without loss of generality
\mr
\item "{$\bigotimes_1$}"  $(\forall \alpha < \beta < \lambda)
(\forall B \in [\lambda]^{< \lambda})
(\exists^{\kappa^+}A \in {\Cal A})[\{\alpha,\beta\} \subseteq A \and A \cap
B \subseteq \{\alpha,\beta\}]$. \nl
\mn
[Why?  As we can use $\bigl\{ \{2 \alpha:\alpha \in A\}:A \in {\Cal A}
\bigr\}$, \wilog \, $\bigcup\{A:A \in {\Cal A}\} = \{2\alpha:\alpha <
\lambda\}$ and choose $A_{\alpha,\beta,\gamma} \in [\lambda]^\theta$ for
$\alpha < \beta < \gamma < \lambda$ such that $\{\alpha,\beta\} \subseteq
A_{\alpha,\beta,\gamma}$ and $\langle A_{\alpha,\beta,\gamma} \backslash
\{\alpha,\beta\}:\alpha < \beta < \gamma < \lambda \rangle$ are pairwise
disjoint subsets of $\{2 \alpha +1:\alpha < \lambda\}$ each of cardinality
$\theta$ and replace ${\Cal A}$ by ${\Cal A}^* =: {\Cal A} \cup \{
A_{\alpha,\beta,\gamma}:\alpha < \beta < \gamma < \lambda\}$.
Now clause (A), (D), (E) are not affected.
Clearly clause $(B)_1$ holds (i.e. ${\Cal A}^* \subseteq [\lambda]^\theta$
and $A \ne B \in {\Cal A}^* \Rightarrow |A \cap B| < \sigma$).  Also clause
(C) is inherited by any extension of the original ${\Cal A}$.  Lastly for
clause $(B)_2$, if ${\Cal A}' \subseteq {\Cal A}^*,|{\Cal A}'| < \kappa$,
let $\langle A_\zeta:\zeta < \zeta^* \rangle$ be a list of ${\Cal A}' \cap
{\Cal A}$ as guaranteed by $(B)_2$ and let $\langle A_\zeta:\zeta \in
[\zeta^*,\zeta^* + |{\Cal A}' \backslash {\Cal A}|) \rangle$ list with no
repetitions ${\Cal A}' \backslash {\Cal A}$, now check.]
\sn
\item "{$\bigotimes_2$}"   ${\Cal B}$ is a basis of $Y^*$ of cardinality
$\theta^*$, and for part (2), ${\Cal B}$ is as there. \nl
\mn
[Why?  Straight.]
\ermn
Let ${\Cal A} = \{A_\zeta:\zeta < \lambda^*\}$ and ${\Cal B} = \{b_i:i <
\theta^*\}$. \nl
We define a forcing notion $P$:
\sn
$p \in P$ has the form $p = (u,u_*,v,v_*,\bar w) = (u^p,u^p_*,v^p,v^p_*,
\bar w^p)$ such that:
\mr
\item "{$(\alpha)$}"  $u_* \subseteq u \in [\lambda]^{< \kappa}$
\sn
\item "{$(\beta)$}"  $v_* \subseteq v \in [\lambda^*]^{< \kappa}$
\sn
\item "{$(\gamma)$}"  $\bar w = \bar w^p = 
\langle w_{\zeta,i}:\zeta \in v_* \text{ and } i < \theta^* \rangle =
\langle w^p_{\zeta,i}:\zeta \in v_*,i < \theta^* \rangle$
\sn
\item "{$(\delta)$}"  $w_{\zeta,i} \subseteq u_*$ and \nl
$b_i \cap b_j = \emptyset \Rightarrow w_{\zeta,i} \cap w_{\zeta,j} =
\emptyset$ this is toward being Hausdorff
\sn
\item "{$(\varepsilon)$}"  $\zeta \in v_* \Rightarrow A_\zeta \subseteq u$
\sn
\item "{$(\zeta)$}"  letting $A^p_\zeta = \cup\{w_{\zeta,i}:i < \theta^*\}
\cap A_\zeta$ for $\zeta \in v^p_*$ it has cardinality $\theta$ and for
simplicity even order type $\theta$ and letting
$\langle \gamma^p_{\zeta,j}:j < \theta \rangle$ list its members with no
repetitions we have \nl
$w^p_{\zeta,i} \cap A^p_\zeta = \{\gamma^p_{\zeta,j}:j < \theta \text{ and }
j \in b_i\}$
\sn
\item "{$(\eta)$}"  if $\zeta \in v^p_*,i < \theta^*$ and $\xi \in v^p_*$
\ub{then} the set ${\Cal U}^p_{\zeta,\xi,i}$ is an open subset (for part 
(2), clopen subset) of the space $Y^*$ where ${\Cal U}^p_{\zeta,\xi,i} =:
\{j < \theta:\gamma^p_{\xi,j} \in w^p_{\zeta,i}\}$.
\ermn
$\bigoplus \qquad$  \ub{convention} if $\zeta \in \lambda^* \backslash v^p_*$ 
we stipulate $w^p_{\zeta,i} = \emptyset$.
\mn
The \ub{order} is: $p \le q$ iff $u^p \subseteq u^q,u^p_* =u^q_* \cap u^p,
v^p \subseteq v^q,v^p_* = v^q_* \cap v^p$ and $\zeta \in v^p_* \Rightarrow
w^p_{\zeta,i} = w^q_{\zeta,i} \cap u^p$.
\mn
Clearly
\mr
\item "{$(*)_0$}"  $P$ is a partial order. 
\ermn
What is the desired space in $V^P$?  We define a $P$-name
${\underset\tilde {}\to X^*}$ as follows: \nl
set of points $\bigcup\{u^p_*:p \in {\underset\tilde {}\to G_P}\}$ \nl
The topology is defined by the following basis: \nl
$\{ \dbca_{\ell < n} {\underset\tilde {}\to {\Cal U}_{\zeta_\ell,i_\ell}}:
n < \omega,\zeta_\ell < \lambda^*,i_\ell < \theta^*\}$ where \nl
${\underset\tilde {}\to {\Cal U}_{\zeta,i}}
[{\underset\tilde {}\to G_P}] = \cup\{w^p_{\zeta,i}:p \in 
{\underset\tilde {}\to G_P},\zeta \in v^p_*\}$ \nl
(for part (2), also their compliments and even their Boolean combinations)
\mr
\item "{$(*)_1$}"  for $\alpha < \lambda$ and $p \in P$ will have
$p \Vdash ``\alpha \in \underset\tilde {}\to X'$ iff $\alpha \in u^p_*$ and
$p \Vdash ``\alpha \notin \underset\tilde {}\to X'$ iff $\alpha \in
u^p_\alpha \backslash u^p_*$
\sn
\item "{$(*)_2$}"  $P$ is $\kappa$-complete, in fact if $\langle
p_\varepsilon:\varepsilon < \delta \rangle$ is increasing in $P$ and $\delta
< \kappa$ then $p = \dbcu_{\varepsilon < \delta} p_\varepsilon$ is an upper
bound where $u^p = \dbcu_{\varepsilon < \delta} u^{p_\varepsilon},
u^p_* = \dbcu_{\varepsilon < \delta} u^{p_\varepsilon}_*,
v^p = \dbcu_{\varepsilon < \delta} v^{p_\varepsilon},v^p_* = 
\dbcu_{\varepsilon < \delta} v^{p_\varepsilon}_*$
and $w^p_{\zeta,i} = \cup\{w^{p_\varepsilon}_{\zeta,i}:\zeta \in
v^{p_\varepsilon}_*,\varepsilon < \delta\}$ \nl
[why?  straight]
\sn
\item "{$(*)_3$}"  $P' = \{p \in P:\text{if } \zeta < \lambda^* \text{ and }
|A_\zeta \cap u^p| \ge \sigma$ then $\zeta \in v^p\}$ is a dense subset of $P$
\nl
[why?  for any $p \in P$ we define by induction on $\varepsilon \le
\sigma^+:p_\varepsilon \in P$, increasingly continuous with $\varepsilon$.
Let $p_0 = p$, if $p_\varepsilon$ is defined, we define $p_{\varepsilon +1}$
by

$$
v^{p_{\varepsilon+1}} = \{\zeta < \lambda^*:\zeta \in v^{p_\varepsilon}
\text{ or } |A_\zeta \cap u^{p_\varepsilon}| \ge \sigma\}
$$

$$
v^{p_{\varepsilon +1}}_* = v^{p_\varepsilon}_*
$$

$$
u^{p_{\varepsilon +1}} = u^{p_\varepsilon} \cup \bigcup\{A_\zeta:
\zeta \in v^{p_{\varepsilon +1}}\}
$$

$$
u^{p_{\varepsilon +1}}_* = u^{p_\varepsilon}_* (= u^p_*)
$$

$$
w^{p_{\varepsilon +1}}_{\zeta,i} \text{ is: } w^{p_\varepsilon}_{\zeta,i}
\text{ if } \zeta \in v^{p_\varepsilon}_*,i < \theta^*
$$
(and there are no other cases). \nl
By assumption $(A)(ii)$, the set $v^{p_{\varepsilon +1}}$ has cardinality
$< \kappa$, so $p_{\varepsilon +1}$ belongs to $P$.
\ermn
Clearly $p_\varepsilon \le p_{\varepsilon +1} \in P$.
\mn
Now for $\varepsilon$ limit let $p_\varepsilon = \dbcu_{\xi < \xi} p_\xi$.
So we can carry the definition.  Now $p_{\sigma^+} = \dbcu_{\varepsilon <
\sigma} p_\varepsilon$ is as required because if $A_\zeta \in {\Cal A},
|A_\zeta \cap u^{p_{\sigma^+}}| \ge \sigma$ then for some $\varepsilon <
\sigma^+,|A_\zeta \cap u^{p_\varepsilon}| \ge \sigma$ hence $\zeta \in
v^{p_{\varepsilon +1}}$ hence $A_\zeta \subseteq u^{p_{\varepsilon +1}}
\subseteq u^{p_{\sigma^+}}$. \nl
Note that we use here $\sigma^+ < \kappa$.]
\mr
\item "{$(*)_4$}"  $P$ satisfies the $\kappa^+$-c.c. \nl
[Why?  Also easy.  Let $p_j \in P$ for $j < \kappa^+$, \wilog \,
$p_j \in P'$ for $j < \kappa^+$.  Now by the $\Delta$-system lemma for some
unbounded $S \subseteq \kappa^+$ and $v^\otimes \in [\lambda^*]^{< \kappa},
u^\otimes \in [\lambda]^{< \kappa}$ we have: \nl
$j \in S \Rightarrow v^\otimes \subseteq v^{p_j} \and u^\otimes \subseteq
u^{p_j} \text{ and } \langle v^{p_j} \backslash v^\otimes:j \in S \rangle$
are pairwise disjoint and $\langle u^{p_j} \backslash u^\otimes:j \in S
\rangle$ are pairwise disjoint.  Without loss of generality otp$(v^{p_j})$,
otp$(u^{p_j})$ are constant for $j \in S$ and any two $p_i,p_j$ are 
isomorphic over $v^\otimes,u^\otimes$ (if not clear see \scite{t.4}). \nl
Now for $j_1,j_2 \in S,p_{j_1},p_{j_2}$ are compatible because of the
following $(*)_5$]
\sn
\item "{$(*)_5$}"   assume $p^1,p^2 \in P$ satisfies
{\roster
\itemitem{ $(i)$ }  $v^{p^1}_* \cap (v^{p^2} \backslash v^{p^2}_*) 
= \emptyset$ and 
$u^{p^1}_* \cap (u^{p^2} \backslash u^{p^2}_*) = \emptyset$
\sn
\itemitem { $(ii)$ }  $v^{p^2}_* \cap (v^{p^1} \backslash v^{p^1}_*) = 
\emptyset$ and $u^{p^2}_* \cap (u^{p^1} \backslash u^{p^1}_*) =
\emptyset$
\sn
\itemitem{ $(iii)$ }  if $\zeta \in v^{p^1}_* \cap v^{p^2}_*$ then
$A^{p^1}_\zeta = A^{p^2}_\zeta$ and \nl
$i < \theta^* \Rightarrow w^{p^1}_{\zeta,i} \cap (u^{p^1} \cap 
u^{p^2}) = w^{p^2}_{\zeta,i} \cap (u^{p^1} \cap u^{p^2})$
\sn
\itemitem{ $(iv)_1$ }  if $\zeta \in v^{p^1}_* \backslash v^{p^2}_*$ then
$|A_\zeta \cap u^{p^2}| < \sigma$ or just $|A^{p^1}_\zeta \cap u^{p^2}| <
\sigma$
\sn
\itemitem{ $(iv)_2$ }  similarly \footnote{note that if $p^1,p^2 \in P'$, then
clauses $(iv)_1,(iv)_2$ holds automatically.} for 
$\zeta \in v^{p^2}_* \backslash v^{p^1}_*$ 
\endroster}
\ermn
\ub{then} there is $q \in P$ such that: 
\mr
\item "{$(a)$}"  $v^q = v^{p^1} \cup v^{p^2}$ 
\sn
\item "{$(b)$}"  $v^q_* = v^{p^1}_* \cup v^{p^2}_*$ 
\sn
\item "{$(c)$}"  $u^q = u^{p^1} \cup u^{p^2}$ 
\sn
\item "{$(d)$}"  $u^q_* = u^{p^1}_* \cup u^{p^2}_*$
\sn
\item "{$(e)$}"  $p^1 \le q,p^2 \le q$.
\ermn
[Why?   To define the condition $q$ we just have to define
$w^q_{\zeta,i}$ (for $\zeta \in v^q_* = v^{p^1}_* \cup v^{p^2}_*$ 
and $i < \theta^*$).  If $\zeta \in v^{p^1}_* \cap v^{p^2}_*$ we let 
$w^q_{\zeta,i} = w^{p^1}_{\zeta,i} \cup w^{p^2}_{\zeta,i}$ for $i < \theta^*$.
\nl
Now for $\ell = 1,2$, let $v^{p^\ell}_* \backslash v^{p^{3 - \ell}}_*$
be listed as $\langle \Upsilon(\varepsilon,\ell):\varepsilon < 
\varepsilon^*_\ell \rangle$ with no repetitions such that 
$B^\ell_\varepsilon =: A^{p^\ell}_{\Upsilon(\varepsilon,\ell)} \cap 
(\dbcu_{\xi < \varepsilon} A^{p^\ell}_{\Upsilon(\xi,\ell)} \cup 
u^{p^{3 - \ell}})$ is of cardinality $< \sigma$. \nl
\sn
[Why possible?  By the assumption $(B)_2$ and clause $(iv)$ above.] \nl
Now for each
$\zeta \in v^{p^{3 - \ell}}_* \backslash v^{p^\ell}_*$ we choose by induction
on $\varepsilon < \varepsilon^*_\ell$ the sequence 
$\langle w^{\ell,\varepsilon}_{\zeta,i}:i < \theta^* \rangle$ such that
\nl
1)  $w^{\ell,\varepsilon}_{\zeta,i} \subseteq u^{p^{3 - \ell}} \cup 
\dbcu_{\xi < \varepsilon} A^{p^\ell}_{\Upsilon(\xi,\ell)}$.
\nl
2)  $w^{\ell,\varepsilon}_{\zeta,i}$ is increasingly continuous with
$\varepsilon$.
\nl
3) $w^{\ell,0}_{\zeta,i} = w^{p^{3 - \ell}}_{\zeta,i}$. \nl
4) $\varepsilon' < \varepsilon \Rightarrow w^{\ell,\varepsilon}_{\zeta,i}
\cap (u^{p^{3 - \ell}} \cup \dbcu_{\xi < \varepsilon_1} 
A^{p^\ell}_{\Upsilon(\xi,\ell)}) = w^{\ell,\varepsilon'}_{\zeta,i}$.
\nl
5) if $i < j < \theta^*$ and $b_i \cap b_j = \emptyset$ (hence
$w^{p^\ell}_{\zeta,i} \cap w^{p^\ell}_{\zeta,j} = \emptyset)$ then 
$w^{\ell,\varepsilon}_{\zeta,i} \cap w^{\ell,\varepsilon}
_{\zeta,j} = \emptyset$.
\nl
6) $\{j < \theta:\gamma^{p^\ell}_{\Upsilon(\varepsilon,\ell),j} \in
w^{\ell,\varepsilon +1}_{\zeta,i}\}$ is an open set in $Y^*$
(for part (2): clopen).]
\sn
For $\varepsilon = 0$ use clause (3) and for limit $\varepsilon$ take unions
(see clause (2)).  Suppose we have defined for $\varepsilon$ and let us
define for $\varepsilon +1$.  By an assumption above $B^\ell_\varepsilon$ has
cardinality $< \sigma$ and so $Z^\ell_\varepsilon = \{j < \theta:
\gamma^{p^\ell}_{\Upsilon(\varepsilon,\ell),j} \in B^\ell_\varepsilon\}$
is a subset of $\theta$ of cardinality $< \sigma$.  Hence, by assumption (E),
we can find a sequence $\langle t_j(\varepsilon,\ell):j \in Z^\ell
_\varepsilon \rangle$ such that: $t_j(\varepsilon,\ell) < \theta^*$ and
$j \in b_{t_j(\varepsilon,\ell)}$ for $j \in Z^\ell_\varepsilon$ and
$\langle b_{t_j(\varepsilon,\ell)}:j \in Z^\ell_\varepsilon \rangle$ is a
sequence of pairwise disjoint subsets of $Y^*$.  \nl
Lastly, we let

$$
\align
w^{\ell,\varepsilon +1}_{\zeta,i} = w^{\ell,\varepsilon}_{\zeta,i}
\cup \bigl\{ \gamma^{p^\ell}_{\Upsilon(\varepsilon,\ell),s}:
&\text{ for some } j \in Z^\ell_\varepsilon \text{ we have}: \\
  &\,\gamma^{p^\ell}_{\Upsilon(\varepsilon,\ell),t_j(\varepsilon,\ell)}
\in w^{\ell,\varepsilon}_{\zeta,i} \text{ and} \\
  &\,s \in b_{t_j(\varepsilon,\ell)} \bigr\}.
\endalign
$$
\mn
Clearly this is O.K. and we are done.  Remember that the union of
$< \sigma$ set from ${\Cal B}$ is clopen for part (2).]
\mr
\item "{$(*)_6$}"  in $(*)_5$ if in addition for $\ell =1,2$ we have
$Z_\ell \subseteq u^{p^\ell} \backslash u^{p^{3 - \ell}}$ such that
$(\forall \zeta \in v^{p^\ell}_*)[|A^{p^\ell}_\zeta \cap Z_\ell| < \sigma]$
\ub{then} we may add to the conclusion \nl
$\ell \in \{1,2\},\zeta \in v^{p^{3 - \ell}}_* \backslash v^{p^\ell}_*,
i < \theta^* \Rightarrow w^q_{\zeta,i} \cap Z_\ell = \emptyset$. \nl
More generally if $g_\ell:(v^{p^{3 - \ell}}_* \backslash v^{p^\ell}_*) \times
\theta^* \times Z_\ell \rightarrow \{0,1\}$ we can add \nl
$\ell \in \{1,2\},\zeta \in v^{p^{3 - \ell}}_* \backslash v^{p^\ell}_*,
i < \theta^*,\gamma \in Z_\ell \Rightarrow [\gamma \in w^q_{\zeta,i} 
\leftrightarrow g_\ell(\zeta,i,\gamma) = 1]$. \nl
[Why?  When for $\zeta \in v^{p^{3 - \ell}}_* \backslash v^{p^\ell}_*$, we
define $\langle w^{\ell,\varepsilon}_{\zeta,i}:i < \theta^* \rangle$ by
induction on $\varepsilon$ we add \nl
(7)  $i < \theta^*,\gamma \in Z_\ell \cap (u^{p^{3 - \ell}} \cup \dbcu
_{\xi < \varepsilon} A^{p^\ell}_{\Upsilon(\xi,\ell)})$ implies $\gamma \in
w^{\ell,\varepsilon}_{\zeta,i} \leftrightarrow g_\ell(\zeta,i,\gamma)=1$.  In
the proof when we use clause (E), instead of using $B^\ell_\varepsilon = 
A^{p^\ell}_{\zeta(\varepsilon,\ell)} \cap (\dbcu_{\xi < \varepsilon} 
A^{p^\ell}_{\zeta(\xi,\ell)} \cup u^{p^{3 - \ell}})$ we use 
$B^\ell_\varepsilon = A^{p^\ell}_{\zeta(\varepsilon,
\ell)} \cap (\dbcu_{\xi < \varepsilon} A^{p^\ell}_{\zeta(\xi,\ell)} \cup
u^{p^{3 - \ell}} \cup Z_\ell)$ which still has cardinality $< \sigma$.]
\mn
Now we come to the main point
\sn
\item "{$(*)_7$}"  in $V^P$, if $i(*) < \text{ cf}(\theta)$ and
$X^* = \dbcu_{i < i(*)} X_i$ then some closed $Y \subseteq X^*$ 
is homeomorphic to $Y^*$.
\ermn
[Why?  Toward contradiction assume $p^* \in P$ and $p^* \Vdash_P ``\langle
{\underset\tilde {}\to X_i}:i < i(*) \rangle$ is a counterexample to
$(*)_7$". \nl
Without loss of generality 
$p^* \Vdash_P ``\langle {\underset\tilde {}\to X_i}:i < i(*) \rangle$ is a
partition of $X^*$, i.e. of $\lambda$". \nl
For each $\alpha < \lambda$ let $\langle(p_{\alpha,j},i_{\alpha,j}):j <
\kappa \rangle$ be such that:
\mr
\widestnumber\item{$(iii)$}
\item "{$(i)$}"  $\langle p_{\alpha,j}:j < \kappa \rangle$ is a maximal
antichain of $P$ above $p^*$
\sn
\item "{$(ii)$}"  $p_{\alpha,j} \Vdash_P ``\alpha \in X_{i_{\alpha,j}}"$,
so $i_{\alpha,j} < i(*)$
\sn
\item "{$(iii)$}"  $p^* \le p_{\alpha,j}$.
\ermn
Now choose a function $F$, Dom$(F) = \lambda$ as follows:

$$
F(\alpha) \text{ is } \bigcup\{u^{p_{\alpha,j}}:j < \kappa\}.
$$
\mn
So we can find $\zeta(*) < \lambda^*$ and $A \subseteq A_{\zeta(*)}$ of
order type $\theta$ such that: if $\alpha \ne \beta$ are from $A$ then 
$\alpha \notin F(\beta)$.  Let $A = \{\beta_\varepsilon:\varepsilon < 
\theta\}$ with no repetitions.  Now we 
shall choose by induction on $\varepsilon \le \theta,p_\varepsilon,
g_\varepsilon$ and if
$\varepsilon < \theta$ also $j_\varepsilon < \kappa$ such that:
\roster
\item "{$(a)$}"  $p_\varepsilon \in P$ and $u^{p_\varepsilon} = 
u^{p^*} \cup \dbcu_{\varepsilon(1) < \varepsilon} 
u^{p_{\beta_{\varepsilon(1)},j_{\varepsilon(1)}}}$ \nl
$v^{p_\varepsilon} = v^{p^*} \cup \dbcu_{\varepsilon(1) < \varepsilon}
v^{p_{\beta_{\varepsilon(1)},j_{\varepsilon(1)}}}$ \nl
$v^{p_\varepsilon}_* = v^{p^*}_* \cup \dbcu_{\varepsilon(1) < \varepsilon}
v^{p_{\beta_{\varepsilon(1)},j_{\varepsilon(1)}}}_*$ \nl
$w^{p_\varepsilon}_{\zeta,i} = w^{p^*}_{\zeta,i} \cup \dbcu_{\varepsilon(1) <
\varepsilon} w^{p_{\beta_{\varepsilon(1)},j_{\varepsilon(1)}}}$ (remember the
convention $\bigoplus$) \nl
(so $p_0 = p^*$)
\sn
\item "{$(b)$}"  $j_\varepsilon = 
\text{ Min}\{j < \kappa:p_{\beta_\varepsilon,j}$ is compatible with 
$p_\varepsilon\}$
\sn
\item "{$(c)$}"  $g_\varepsilon$ is a function, increasing with
$\varepsilon$, from $v^{p_\varepsilon}_*
\times \theta^*$ into the family  of open subsets of $Y^*$ (for part (2),
clopen)
\sn
\item "{$(d)$}"  if $b_{i_1} \cap b_{i_2} = \emptyset$ then $g_\varepsilon
(\zeta,i_1) \cap g_\varepsilon(\zeta,i_2) = \emptyset$ (if defined)
\sn
\item "{$(e)$}"  letting $\Upsilon_\varepsilon = \text{ otp}\{\xi <
\varepsilon:i_{\beta_\varepsilon,j_\varepsilon} = i_{\beta_\xi,j_\xi}\}$ we
have for $\zeta \in v^{p_\varepsilon}_*$:
$$
\beta_\varepsilon \in w^{p_{\varepsilon +1}}_{\zeta,i} \Leftrightarrow
\Upsilon_\varepsilon \in g_\varepsilon(\zeta,i)
$$
\sn
\item "{$(f)$}"  $p_\varepsilon$ is increasing continuous.
\ermn
No problem to carry the definition.  As for $\varepsilon$ successor, for 
this $(*)_6$ was prepared.  In limit $\varepsilon$ take union.  In all 
cases $j_\varepsilon$ is well defined by clause $(i)$ above.  Let 
$i^* < i(*)$ be minimal such that the set $Z = \{\varepsilon < \theta:
i_{\beta_\varepsilon,j_\varepsilon} = i^*\}$ has cardinality $\theta$.  
Note: $\zeta(*) \notin
v^{p_{\beta_\varepsilon,j}}$ as $A \cap F(\beta_\varepsilon)$ is a
singleton so $|A \cap u^{p_{\beta_\varepsilon,j}}| \le 1$ and 
$p_{\beta_\varepsilon,j} \in P'$.  Now we define $p$:

$$
u^p = u^{p_\theta}
$$

$$
v^p = v^{p_\theta} \cup \{\zeta(*)\}
$$

$$
v^p_* = v^{p_\theta}_* \cup \{\zeta(*)\}
$$
\mn
$w^p_{\zeta,i}$ is
\mr
\item "{$(\alpha)$}"  $w^{p_\theta}_{\zeta,i}$ \ub{if} $\zeta \in 
v^{p_\theta}$
\sn
\item "{$(\beta)$}"  $\{\beta_\varepsilon:\varepsilon \in Z \text{ and otp}
(Z \cap \varepsilon) \in b_i\}$ \ub{if} $\zeta = \zeta(*)$ so
\sn
\item "{$(\gamma)$}"  $A^p_{\zeta(*)} = \{\beta_\varepsilon:\varepsilon \in
Z\}$ and $\gamma^p_{\zeta(*),\varepsilon}$ is the $\varepsilon$-th member
of $A^p_{\zeta(*)}$.
\ermn
We can easily check that $p \in P$ and $p^* \le 
p_{\beta_\varepsilon,j_\varepsilon} \le p \in P$ (but we do not ask
$p_\varepsilon \le p$).
Clearly $p$ forces that $\{\beta_\varepsilon:\varepsilon \in Z\}$ is
included in one $X_i$. \nl
Let $g:\theta \rightarrow \lambda$ be $g(\xi) = \beta_\varepsilon$ when
$\xi < \theta,\varepsilon \in Z$, otp$(Z \cap \varepsilon) = \xi$.
Now $p \ge p^*$ and we are done by $(*)_8$ below.]
\mr
\item "{$(*)_8$}"  if $p \in P$ and $\zeta \in v^p_*$ then \nl
$p \Vdash$ ``the mapping $j \mapsto \gamma^p_{\zeta,j}$ for $j < \theta$ is
a homeomorphism from $Y^*$ onto the closed subspace
$\underset\tilde {}\to X \restriction \{\gamma^p_{\zeta,j}:j < \theta\}$ of
$\underset\tilde {}\to X$" \nl
[Why?  Let $p \in G,G \subseteq P$ is generic over $V$.
{\roster
\itemitem{ $(\alpha)$ } If $b \in {\Cal B}$, then for 
some open set ${\Cal U}$ of $\underset\tilde {}\to X$ (clopen for
part (2)) we have
$$
{\Cal U} \cap \{\gamma^p_{\zeta,j}:j < \theta\} = \{\gamma^p_{\zeta,j}:j \in
b\}
$$
[Why?  As $b = b_i$ for some $i < i(*)$ and $p$ forces that \nl
$w_{\zeta,i} \cap \{\gamma^p_{\zeta,j}:j < \theta\} = \{\gamma^p_{\zeta,j}:
j \in b_i\}$.]
\sn
\itemitem{ $(\beta)$ }  If $b$ is an open set for $Y^*$, then for some open
subset ${\Cal U}$ of $\underset\tilde {}\to X$ we have
$$
{\Cal U} \cap \{\gamma^p_{\zeta,j}:j < \theta\} = \{\gamma^p_{\zeta,j}:
j \in b\}
$$
[Why?  As $b = \dbcu_{i \in Z} b_i$ for some $Z \subseteq \theta^*$ 
and apply clause $(\alpha)$]
\sn
\itemitem{ $(\gamma)$ }  if ${\Cal U}$ is an open subset of
$\underset\tilde {}\to X$ and $\gamma^p_{\zeta,j(*)} \in {\Cal U}$
(and $\zeta \in u^p_*$), then for some $i(*) < \theta^*$ we have
$$
\gamma^p_{\zeta,j(*)} \in w^p_{\zeta,i(*)} \cap
\{\gamma^p_{\zeta,j}:j < \theta\} \subseteq
{\underset\tilde {}\to {\Cal U}_{\zeta,i(*)}} \cap
\{\gamma^p_{\zeta,j}:j < \theta\} \subseteq {\Cal U}.
$$
[Why?  By the definition of the topology $\underset\tilde {}\to X$ we can
find $n < \omega,\xi_\ell < \lambda^*$ and $i_\ell < \theta^*$ such that
$\gamma^p_{\zeta,j(*)} \in \dbca_{\ell < n} 
{\underset\tilde {}\to {\Cal U}_{\xi_\ell,i_\ell}}[G]$.  
We can find $q \in P$ such that $p \le q$ and $\xi_\ell \in v^q_*$ for
$\ell < n$.  For each $\ell < n$,
by clause $(\eta)$ in the definition of $P$ we have
${\Cal U}^q_{\zeta,\xi_\ell,j_\ell}$ is an open set
for $Y^*$, and necessarily $j(*) \in {\Cal U}^q_{\zeta,\xi_\ell,j_\ell}$.  Let
$i(*)$ be such that $j(*) \in b_{i(*)} \subseteq \dbca_{\ell < n}
{\Cal U}^q_{\zeta,\xi_\ell,j_\ell}$ hence
$\gamma^p_{\zeta,j(*)} \in {\underset\tilde {}\to {\Cal U}_{\zeta,i(*)}}[G]
\cap \{\gamma^p_{\zeta,j}:j < \theta\} \subseteq \dbca_{\ell < n}
{\underset\tilde {}\to {\Cal U}_{\xi_\ell,j_\ell}}[G] \subseteq {\Cal U}$ as
required.  So $i(*)$ is as required.]
\sn
\itemitem{ $(\delta)$ }  $\{\gamma^p_{\zeta,j}:j < \theta\}$ is a closed
subset of $\underset\tilde {}\to X$ \nl
[Why?   Let $\beta \in \lambda \backslash \{\gamma^p_{\zeta,j}:j < \theta\}$
and let $p \le q \in P$; it suffices to find $q^+,q \le q^+ \in P$ and 
$\xi \in v^{q^+}_*$ and $i < \theta^*$ such that 
$\beta \in w^{q^+}_{\xi,i}$ and
$w^{q^+}_{\xi,i} \cap \{\gamma^p_{\zeta,j}:j < \theta\} = \emptyset$.
Without loss of generality $\beta \in u^q_*$. \nl
We can find a set $u \subseteq u^q_*$ such that $\beta \in u,A^q_\zeta \cap
u = \emptyset$ and $\zeta' \in v^q_* \Rightarrow \{j < \theta:
\gamma^q_{\zeta',j} \in u\}$ is a clopen subset of $Y$, (just as in the 
proof of $(*)_5$).
We can find by $\otimes_1$ $\xi \in \lambda^* \backslash v^q$ such that
$\{\emptyset\} = A_\xi \cap u^q_*$ (why? apply $\otimes_1$ with $\alpha <
\beta \in \lambda \backslash u^q$ and $B = u^q$) and let
$\gamma_{\varepsilon,i} \in A_\xi$ for $i < \theta$ be increasing.  We define
$q^+$ as follows.
\endroster}
\ermn

$$
v^{q^+} = v^q \cup \{\xi\}
$$

$$
v^{q^+}_* = v^q_* \cup \{\xi\}
$$

$$
u^{q^+} = u^q \cup A_\xi
$$

$$
u^{q^+}_* = u^q_* \cup \{\gamma_{\xi,j}:j < \theta\}
$$

$w^{q^+}_{\zeta,i} \text{ is } w^q_{\zeta,i} \text{ if } \zeta \in v^q_*
\text{ and is } \{\gamma_{\xi,j}:j \in b_i\} \cup u$ if $\zeta = \xi \and 
0 \in b_i \text{ and is } \{\gamma_{\xi,j}:j \in b_i\}$ if $\zeta = \xi 
\and 0 \notin b_i$.
\mn
Lastly, we would like to know that 
$\underset\tilde {}\to X$ is a Hausdorff space.  We prove more
\mr
\item "{$(*)_9$}"  In $V^P$ if $u_1 \subseteq u_2 \in [\lambda]^{< \sigma}$
then for some $\zeta,i$ we have
$$
w_{\zeta,i} \cap u_2 \cap \underset\tilde {}\to X = u_1
\cap \underset\tilde {}\to X
$$
[Why?  Let $p_0 \in P$ force ${\underset\tilde {}\to u_1} \subseteq
{\underset\tilde {}\to u_2}$ form a counterexample, as 
$P$ is $\kappa$-complete
some $p_1 \ge p_0$ forces ${\underset\tilde {}\to u_1} = u_2,
{\underset\tilde {}\to u_2} = u_2$ and $p_1 \in P'$.  Necessarily $u_2
\subseteq u^{p_1}_*$. \nl
Let $\zeta(*) \in \lambda^* \backslash v^{p_1}$ be such that $A_{\zeta(*)}
\cap u^{p_1} = \emptyset$ (as in the proof of $(*)_8$).  Let
$\gamma_{\zeta(*),j} \in A_{\zeta(*)}$, for $j < \theta$ be increasing.  Let
$u \subseteq u^{p_1}_*$ be such that $u \cap u_2 = u_1$ and $\zeta' \in
v^{p_1}_* \Rightarrow \{j < \theta:\gamma^{p_1}_{\zeta',j} \in u\}$ is clopen
in $Y$ (exists as in the proof of $(*)_5$) and define $q \in P$:

$$
u^q = u^{p_1} \cup u_2
$$

$$
u^q_* = u^{p_1}_* \cup (u_2 \backslash u^{p_1})
$$

$$
v^q = v^p \cup \{\zeta(*)\}
$$

$$
v^q_* = v^{p_1} \cup \{\zeta(*)\}
$$

$w^q_{\zeta,i} \text{ is: } w^{p_1}_{\zeta,i} \text{ if } \zeta \in v^q,
\text{ is } \{\gamma_{\zeta(*),j}:j \in b_i\} \cup u \text{ if } 
\zeta = \zeta(*) \and 0 \in b_i \text{ and is } \{\gamma_{\zeta(*),j}:j \in
b_i\} \text{ if } \zeta = \zeta(*) \and 0 \notin b_i$.]
\sn
Together all is done.   \hfill$\square_{\scite{t.2}}$
\endroster
\enddemo
\bn
\ub{\stag{p.2} Proof of \scite{t.5}(2)}

Saharon: after \scite{p.1}.
\bn
\ub{\stag{p.3} Proof of \scite{3.7}} \nl
\demo{Proof}  Let $F:\lambda \rightarrow [\lambda]^{\le \kappa}$ be given.
Choose by induction on $\zeta \le \lambda$ a set $U_\zeta \subseteq \lambda$
and $g_\zeta:U_\zeta \rightarrow \kappa^+$, both increasingly continuous
with $\zeta$ such that:
\mr
\item "{$(*)(i)$}"  if $\alpha \in U_\zeta$ then $F(\alpha) \subseteq
U_\zeta$ and
\sn
\item "{$(ii)$}"  if $\alpha \in U_\zeta$ then $F(\alpha) \backslash
\{\alpha\} \subseteq \{\beta \in U_\zeta:g_\zeta(\beta) \ne g_\zeta
(\alpha)\}$.
\ermn
For $\zeta = 0$ let $U_\zeta = \emptyset = g_\zeta$, for $\zeta$ limit
take unions.  If $U_\zeta = \lambda,U_{\zeta + 1} = U_\zeta,g_{\zeta +1} =
g_\zeta$, otherwise let $\alpha_\zeta = \text{ Min}\{\lambda \backslash
U_\zeta\}$ and let $W_\zeta \in [\lambda]^{\le \kappa}$ be such that
$\alpha_\zeta \in W_\zeta$ and 
$(\forall \alpha \in W_\zeta)[F(\alpha) \subseteq
W_\zeta]$.  Let $\varepsilon_\zeta = 
\sup\{g_\zeta(\beta):\beta \in U_\zeta \cap
W_\zeta\}$ so $\varepsilon_\zeta < \kappa^+$ and let 
$U_{\zeta +1} = U_\zeta \cup W_\zeta,g_{\zeta +1}$ extends $g_\zeta$
such that $g_{\zeta +1} \restriction (W_\zeta \backslash U_\zeta)$ is one
to one with range $[\varepsilon,\varepsilon + \kappa)$. \nl
Now applying $(C)^+$ to the partition which $\dbcu_\zeta g_\zeta$ defines,
we get some $A \in {\Cal A}$ on which $\dbcu_\zeta g_\zeta$ is constant
so by $(*)(ii)$ we are done.  \hfill$\square_{\scite{g.5}}$
\enddemo
\bn
\ub{\stag{p.4} Proof of \scite{3.10}(2)}

Saharon.
\newpage

\head{Private Appendix} \endhead  \resetall 
\bn
\ub{moved from p.2} \nl
\ub{Problem}:  1) Assume ``there is a supercompact cardinal" (as in
\cite{Sh:108}, \cite{HJSh:249}) (for (A), (B), (C) \ub{or} ``there is a
2-huge cardinal" (as in \cite{HJSh:249} for (B), (C)
\mr
\item "{$(A)$}"  can we get the consistency of the assumption of 2.2?
(not just of 2.5)?
\sn
\item "{$(B)$}"  can we get the consistency required for 2.5 or even
2.6 for $X' \rightarrow (Y^*)^1_\theta,\theta > \kappa$? \nl
For this it suffices:
{\roster
\itemitem{ $(*)$ }  if $h:\lambda \rightarrow \kappa^+$ then for some
$A \in {\Cal A}$ and $\zeta < \kappa^+$ we have \nl
$|A \cap h^{-1}\{\zeta\}| = \theta$.
\endroster}
\ermn
2) Can we get the examples with GCH?
\bigskip

\remark{Remark}  Seems easier to get $\theta > \kappa$, if say the space
has size $\aleph_1$ and we have MA,Ded not ?
\endremark
\bn
\centerline{$* \qquad * \qquad *$}
\bn
\ub{Moved from pgs.9-10} \nl
If we look at spaces with clopen basis, it seems easier to force such spaces.
To enable the reader to read one proof the similar parts of the proofs are
repeated.
\proclaim{\stag{c.4} Theorem}  Assume
\mr
\item "{$(A)$}"  $\lambda > \kappa > \theta > \sigma \ge \aleph_0$ and
$\kappa = \kappa^{< \kappa},\kappa > \theta^*$
\sn
\item "{$(B)$}"  ${\Cal A} \subseteq [\lambda]^\theta$ and $A_1 \ne A_2 \in
{\Cal A} \Rightarrow |A_1 \cap A_2| < \sigma$;
\sn
\item "{$(C)$}"  if $F:\lambda \rightarrow [\lambda]^{\le \kappa}$ then for
some $A \in {\Cal A}$ is $F$-free, i.e. $\alpha \ne \beta \in A \Rightarrow
\alpha \in F(\beta)$
\sn
\item "{$(D)$}"  $Y^*$ is a $T_3$ topological space with set of points
$\theta$ and with a clopen basis ${\Cal B} = \{b_i:i < \theta^*\}$
\sn
\item "{$(E)$}"  if $Y \subseteq Y^*$ has cardinality $< \sigma$ \ub{then}
$Y$ is closed.
\ermn
\ub{THEN} for some $\kappa$-complete $\kappa^+$-c.c. forcing notion $P$ in
$V^P$ there is $X^*$ such that:
\mr
\item "{$(a)$}"  $X^*$ is a $T_3$ topological space even with a clopen
basis of size $|{\Cal A}| + \theta^*$ with $\lambda$ points (so its
compatification is also compact and the other properties remain, also the
Boolean Algebra of clopen sets is a free Boolean Algebra so the space is
${}^\kappa 2$)
\sn
\item "{$(b)$}"  $X^* \rightarrow (Y^*)^1_{< \text{ cf}(\theta)}$, i.e. if
$X^* = \dbcu_{i < i(*)} X_i,i(*) < \text{ cf}(\theta)$, \ub{then} some closed
subspace $Y$ of $Y^*$ homeomorphic to $Y^*$, is included in some single $X_i$.
\endroster
\endproclaim
\bn
\centerline{$* \qquad * \qquad *$}
\bn
\ub{Moved from pgs.11-14} \nl
\demo{Proof of g.4}  Let ${\Cal A} = \{A_\zeta:\zeta < \lambda^*\}$.

We define $P$:
\sn
$p \in P$ has the form $p = (u,v,v_*,\bar w) = (u^p,v^p,v^P_*,\bar w^p)$
such that:
\mr
\item "{$(\alpha)$}"  $u \in [\lambda]^{< \kappa}$
\sn
\item "{$(\beta)$}"  $v_* \subseteq v \in [\lambda^*]^{< \kappa}$
\sn
\item "{$(\gamma)$}"  $\bar w = \langle w_{\zeta,i}:\zeta \in v \text{ and }
i < \theta^* \rangle$
\sn
\item "{$(\delta)$}"  $w_{\zeta,i} \subseteq u$
\sn
\item "{$(\varepsilon)$}"  $\zeta \in v_* \Rightarrow A_\zeta \subseteq u$
\sn
\item "{$(\zeta)$}"  if $\zeta \in v$ and $i < \theta^*$ and $\xi \in v_*
\backslash \{\zeta\}$ then $w_{\zeta,i} \cap A_\xi \in [A_\xi]^{< \sigma}$.
\ermn
$\bigoplus \quad$ If $\zeta \in \lambda^* \backslash v$ we stipulate
$w_{\zeta,i} = \emptyset$.
\sn
The order is: $p \le q$ iff $u^p \subseteq u^q,v^p \subseteq v^q,v^p_* =
v^q_* \cap v^q$ and

$$
w^p_{\zeta,i} = w^q_{\zeta,i} \cap u^p.
$$
\mn
Clearly
\mr
\item "{$(*)_1$}"  $P$ is a partial order. \nl
Our intention is to have $X^*$ as follows: set of points $\lambda$ \nl
the clopen basis is generated by $\{{\underset\tilde {}\to u_{\zeta,i}}:
\zeta < \lambda^*,i < \theta^*\}$ where \nl
${\underset\tilde {}\to u_{\zeta,i}}
[G_P] = \cup \{w^p_{\zeta,i}:p \in G_P,\zeta \in v^p\}$. \nl
No need to prove Hausdorf as otherwise we just need to identify $\alpha,
\beta < \lambda$ if $(\forall \zeta,i)(\alpha \in u_{\zeta,i} \equiv
\beta \in u_{\zeta,i})$ but we will still do it.
\sn
\item "{$(*)_2$}"  $P$ is $\kappa$-complete, in fact if $\langle
p_\varepsilon:\varepsilon < \delta \rangle$ is increasing in $P$ and
$\delta < \kappa$ a limit ordinal then 
$p = \dbcu_{\varepsilon < \delta} p_\varepsilon$ is an upper
bound where $u^p = \dbcu_{\varepsilon < \delta} u^{p_\varepsilon},v^p = 
\dbcu_{\varepsilon < \delta} v^{p_\varepsilon},v^p_* = 
\dbcu_{\varepsilon < \delta} v^{p_\varepsilon}_*$
and $w^p_{\zeta,i} = \bigcup\{w^{p_\varepsilon}_{\zeta,i}:\varepsilon
\text{ satisfies } \zeta \in v^{p_\varepsilon},\varepsilon < \delta\}$  \nl
[why?  straight]
\sn
\item "{$(*)_3$}"  $P' = \{p:\text{ if } \zeta < \lambda^*,|A_\zeta \cap
u^p| \ge \sigma$ then $\zeta \in v_\zeta\}$ \nl
[why? for any $p \in P$ we define by induction on $\varepsilon < \sigma^+,
p_\varepsilon \in P$ increasingly continuous with $\varepsilon$.  Let
$p_\sigma = p$, if $p_\varepsilon$ is defined let

$$
v^{p_{\varepsilon +1}} = \{\zeta < \lambda^*:\zeta \in v^{p_\varepsilon}
\text{ or } (A_\zeta \cap u^{p_\varepsilon}) \ge \sigma\}
$$

$$
v^{p_{\varepsilon +1}}_* = v^{p_\varepsilon}_*
$$

$$
u^{p_{\varepsilon +1}} = u^{p_\varepsilon} \cup \bigcup\{A_\zeta:\zeta \in 
v_{\varepsilon +1}\}
$$
\mn
$w^{p_{\varepsilon +1}}_{\zeta,i}$ is: $w^{p_\varepsilon}_{\zeta,i}$ if
$\zeta \in v^{p_\varepsilon} \text{ and } \emptyset \text{ if } i \in
v^{p_{\varepsilon +1}} \backslash v^{p_\varepsilon}$. \nl
Clearly $p_\varepsilon \le p_{\varepsilon + 1} \in P$.  Now $p_{\sigma^+} =
\dbcu_{\varepsilon < \sigma^+} p_\varepsilon$ is as required.]
\sn
\item "{$(*)_4$}"  $P$ satisfies the $\kappa^+$-c.c.c. \nl
[why?  also easy.  Let $p_j \in P$ for $j < \kappa^+$, \wilog \, $p_j \in
P'$ for $j < \kappa^+$.  By the $\Delta$-system lemma for some unbounded
$S \subseteq \kappa^+$ and $v^\otimes \in [\lambda^*]^{< \kappa},u^\otimes
\in [\lambda]^{< \kappa}$ we have: $j \in S \Rightarrow v^\otimes \subseteq
v^{p_j} \and u^\otimes \subseteq u^{p_j}$ and $\langle v^{p_j} \backslash
v^\otimes:j \in S \rangle$ are pairwise disjoint and $\langle u^{p_j}
\backslash u^\otimes:j \in S \rangle$ are pairwise disjoint.  Without loss of
generality otp$(v^{p_j})$, otp$(u^{p_j})$ are constant and any two are
isomorphic over $v^\otimes,u^\otimes$. \nl
Now for $j_1,j_2 \in S,p_{j_1},p_j$ are compatible.]
\sn
\item "{$(*)_5$}"  $\{p:\alpha \in u^p$ and $\zeta \in v^p\}$ is dense open
for each $\alpha < \lambda$ \nl
[why?  if $p \in P$ lets us define $q:u^q = u^p \cup \{\alpha\},v^q = v^p
\cup \{\zeta\},v^q_* = v^p_*$ and $w^q_{\zeta,i}$ is $w^p_{\zeta,i}$ is well
defined and empty otherwise.]
\ermn
Now we come to the main point
\mr
\item "{$(*)_6$}"  in $V^P$, if $i(*) < \text{ cf}(\theta),X^* = 
\dbcu_{i < i(*)} X_i$ then some closed $Y \subseteq X^*$ is homeomorphic to
$Y^*$. 
\ermn
Why?  Toward contradiction assume $p^* \in P$ and $p^* \Vdash_P ``\langle
{\underset\tilde {}\to X_i}:i < i(*) \rangle$ is a counterexample to $(*)_5$".
Without loss of generality $p^* \Vdash_P ``\langle 
{\underset\tilde {}\to X_i}:i < i(*) \rangle$ is a partition of $X^*$, i.e.
of $\lambda"$.
\sn
For each $\alpha < \lambda$ let $\langle (p_{\alpha,j},i_{\alpha,j}:j <
\kappa \rangle$ be such that:
\mr
\widestnumber\item{$(iii)$}
\item "{$(i)$}"  $\langle p_{\alpha,j}:j < \kappa \rangle$ is a maximal
antichain of $P$ above $p^*$
\sn
\item "{$(ii)$}"  $p_{\alpha,j} \Vdash_P ``\alpha \in X_{i_{\alpha,j}}"$
\sn
\item "{$(iii)$}"  $p_{\alpha,j} \in P'$ (and $p^* \le p_{\alpha,j})$.
\ermn
Now choose a function $F$, Dom$(F) = \lambda$ and $F(\alpha)$ is $\kappa
\cup \{u^{p_{\alpha,j}}:j < \kappa\}$. \nl
So we can find $\zeta < \lambda^*$ such that:

if $\alpha \ne \beta$ are from $A_{\zeta(*)}$ then $\alpha \notin F(\beta)$.
\nl
Let $A_{\zeta(*)} = \{\beta_\varepsilon:\varepsilon < \theta\}$ with no
repetitions.  Now we shall choose by induction on $\varepsilon \le \theta,
p_\varepsilon$ and if $\varepsilon < \theta$ also $j_\varepsilon < \kappa$
such that:
\mr
\item "{$(a)$}"  $p_\varepsilon \in P$ and $u^{p_\varepsilon} =
u^{p^*} \cup \dbcu_{\varepsilon(1) < \varepsilon} u^{p_{\varepsilon(1)},
j_{\varepsilon(1)}}$

$$
v^{p_\varepsilon} = v^{p^*} \cup \dbcu_{\varepsilon(1) < \varepsilon}
v^{p_{\beta_{\varepsilon(1)},j_{\varepsilon(1)}}}
$$

$$
v^{p_\varepsilon}_* = v^{p^*}_* \cup \dbcu_{\varepsilon(1) < \varepsilon}
v^{p_{\beta_{\varepsilon(1)},j_{\varepsilon(1)}}}_*
$$

$$
w^{p_\varepsilon}_{\zeta,i} = w^{p^*}_{\zeta,i} \cup 
\dbcu_{\varepsilon(1) < \varepsilon} 
w^{p_{\beta_{\varepsilon(1)},j_{\varepsilon(1)}}} 
$$

$$
\text{ (remember the convention } \oplus)
$$
\sn
\item "{$(b)$}"  $j_\varepsilon = \text{ Min}\{j < \kappa:
p_{\beta_\varepsilon,j}$ is compatible with $p_\varepsilon\}$.
\ermn
No problem to carry the definition (in limit $\varepsilon$ take union,
$j_\varepsilon$ is well defined by clause (i) above).  Let $i^* < i(*)$ be
such that $Z = \{\varepsilon < \theta:i_{\beta_\varepsilon,j_\varepsilon} =
i^*\}$ has cardinality $\theta$.  Note: $\zeta(*) \notin 
v^{p_{\beta_\varepsilon,j}}$ as $A_{\zeta(*)} \cap F(\beta_\varepsilon)$ is a
singleton so $|A_\zeta \cap u^{p_{\beta_\varepsilon,j_\varepsilon}}| \le 1,
p_{\beta_\varepsilon,j_\varepsilon} \in P'$).  Now we define $p$:

$$
u^p = u^{p_\theta}
$$

$$
v^p = v^{p_\theta} \cup \{\zeta(*)\}
$$

$$
v^p_* = v^{p_\theta}_* \cup \{\zeta(*)\}
$$
\mn
$w_{\zeta,i}$ is:
\mr
\item "{$(c)$}"  $w^{p_\theta}_{\zeta,i}$ \ub{if} $\zeta \in v^{p_\theta}$
\sn
\item "{$(d)$}"  $\{\beta_\varepsilon:\varepsilon \in Z \text{ and otp}
(Z \cap \varepsilon) \in b_i\}$ \ub{if} $\zeta = \zeta(*)$.
\ermn
Let $g:\theta \rightarrow \lambda$ be $g(\xi) = \beta_\varepsilon$
where $\xi < \theta,\varepsilon \in Z$, otp$(Z \cap \varepsilon) = \xi$.  Now
$p \ge p^*$ and $p$ forces that:
\mr
\item "{$(\alpha)$}"  if $b \in {\Cal B}$ then for some open set $U$ of
$\underset\tilde {}\to X,X \cap \{\beta_\varepsilon:\varepsilon \in Z\} =
\{g(\varepsilon):\varepsilon \in b\}$ \nl
[why?  as $b = b_i$ for some $i$ and $p$ forces
${\underset\tilde {}\to w_{\zeta,i}} \cap \{\beta_\varepsilon:\varepsilon \in
Z\} = \{g(\varepsilon):\varepsilon \in b_i\}$]
\sn
\item "{$(\beta)$}"   if 
$i < \theta^*,{\underset\tilde {}\to u_{\zeta(*),i}}
\cap \{\beta_\varepsilon:\varepsilon \in Z\}$ is of the form above \nl
[why? clear]
\sn
\item "{$(\gamma)$}"  if $i < \theta^*,\zeta \in \lambda^* \backslash
\{\zeta(*)\}$ then $u_{\zeta,i} \cap \{\beta_\varepsilon:\varepsilon \in Z\}$
has cardinality $< \sigma$ hence $g^{-1}(u_{\zeta,i} \cap 
\{\beta_\varepsilon:\varepsilon \in Z)$ is a clopen subset of $Y^*$ \nl
[why?  the first phrase as $\zeta(*) \in v^p_*$ and clause $(\zeta)$ in the
definition of $P$; the second follows]
\sn
\item "{$(\delta)$}"  $\{\beta_\varepsilon:\varepsilon \in Z\}$ is a closed
set in $X$ \nl
[why?  let $\beta \in \lambda \backslash \{\beta_\varepsilon:\varepsilon \in
Z\},p \le q \in P$, choose $\xi \in \lambda^* \backslash v^q$ and define
$q^+$

$$
v^{q^+} = v^q \cup \{\xi\}
$$

$$
v^{q^+}_* = v^q_*
$$

$$
u^{q^+} = u^q
$$

$$
w^{q^+}_{\zeta,i} \text{ is } w^q_{\zeta,i} \text{ if } \zeta \in v^q
\text{ and is } \{\beta\} \text{ if } \zeta = \xi.]
$$
\mn
\item "{$(*)_7$}"  in $V^P$, if $u_1 \subseteq u_2 \in [\lambda]^{< \sigma}$
then for some $\zeta,i$ we have
$$
w_{\zeta,i} \cap u_2 = u_1
$$
[why?  let $p_0 \in P$ force ${\underset\tilde {}\to u_1} \subseteq
{\underset\tilde {}\to u_2}$ form a counterexample, as 
$P$ is $\kappa$-complete some $p_1 \ge p_0$ forces 
${\underset\tilde {}\to u_1} = u_1,
{\underset\tilde {}\to u_2} = u_2$ and $p_1 = P'$. \nl
Let $\zeta(*) \in \lambda^* \backslash v^{p_1}$ and define $q \in P$:

$$
u^q = u^{p_1} \cup u_2
$$

$$
v^q = v^p \cup \{\zeta\}
$$

$$
v^q_* = v^p_*
$$

$$
w^q_{\zeta,i} \text{ is: } w^p_{\zeta,i} \text{ if } \zeta \in v^q,u_1
\text{ if } \zeta = \zeta(*).
$$
\mn
Now check. \nl
Together all is done.] \hfill$\square_{\sciteu{g.4}}$\sciteuphantom{g.4}
\ermn
\enddemo
\bigskip

\remark{g.8 Concluding Remark}  1) As in \scite{t.5}(1) we can allow
$\kappa = \theta^*$.
\endremark
\bn
\ub{22/12/97} \nl
1) in \scite{3.1} without $\beth_2 = (2^{\aleph_0})^+$: more involved forcing:
${\Cal U}^p \in [\lambda]^{2^{\aleph_0}}$ but we give only countable
information (or $< 2^{\aleph_0}$?) \nl
2) To get GCH?  Try local forcing?  (decision by u.t.?) 
\nl
3) start with \nl
\mn
4) \ub{31/12/97}  \nl
\ub{Question}:  in \scite{3.f3} try without (d), (e), (f) for every
$f:{\Cal U} \rightarrow X,{\Cal U} \in J^+$?
\mn
\ub{Question}:  Replace ${\Cal F}$ by a family of functions?
\sn
We will prove a more detailed result.  The analysis below is somewhat closed
to the $\lambda$-sets from \cite{Sh:161}.
\bn
\definition{Definition}  We define simultaneously by induction on $\lambda >
\mu$ what is a partial $(\lambda,\mu)$-index system. \nl
1) A partial $\mu$-index system $\Gamma$ is a pair $(S,\bar \lambda) = 
(S^\Gamma,\bar \lambda^\Gamma)$ such that:
\mr
\widestnumber\item{$(f)(\alpha)$}
\item "{$(a)$}"  $\Gamma \subseteq {}^{\omega >}\text{Ord}$
\sn
\item "{$(b)$}"  $\Gamma$ is closed under initial segments,
\sn
\item "{$(c)$}"  $<> \in \Gamma$
\sn
\item "{$(d)$}"  $\bar \lambda = \langle \lambda_\eta:\eta \in S \rangle$
and $\lambda_\eta \ge \mu$
\sn
\item "{$(e)$}"  for each $\eta \in S$ for some $\alpha = \alpha(\eta,\Gamma)
< \text{ cf}(\lambda_\eta)$ we have
$$
\eta \char 94 \langle \beta \rangle \in \Gamma \text{ iff } \beta < \alpha
$$
\sn
\item "{$(f)(\alpha)$}"  if $\lambda_\eta$ is a limit cardinal $> \mu$ then
$\langle \lambda_{\eta \char 94 <\beta>}:\beta < \alpha(\eta,\Gamma) \rangle$
is strictly increasing with limit $\mu$
\sn
\item "{$(\beta)$}"  if $\lambda_\eta$ is a successor cardinal $> \mu$ then
$\lambda_\eta = \lambda^+_{\eta \char 94 <\beta>}$ for $\beta < \alpha(\eta,
\Gamma)$
\sn
\item "{$(\gamma)$}"  if $\lambda_\eta = \mu$, \ub{then} $\alpha(\eta,\Gamma)
=0$, i.e. $\eta$ is $\Delta$-maximal in $S$ (so $\eta \triangleleft \nu \in 
S \Rightarrow \lambda_\eta > \lambda_\nu \ge \mu$
\sn
\item "{$(g)$}"  the set $\{\eta \in S:\lambda_\eta > \mu \text{ but }
\alpha(\eta,\Gamma) < \lambda_\eta\}$ has the form $\{\nu \restriction \ell:
\ell < \ell(*)\}$, where $\nu = \text{ Max}(S)$ is the maximal member of $S$
in the lexicographic order and $\ell^* \le \ell g(\nu)$, let $\nu \restriction
\ell^*$ be called $\nu(\Gamma)$[if $\ell^*=0$ we stipulate $\nu(\Gamma) =
<>^-$ and $\eta \in \Gamma \Rightarrow \neg(\eta \triangleleft \nu)$?].
\endroster
\enddefinition
\bigskip

\definition{Definition}  1) A full $\mu$-index system $\Gamma$ is a partial
$\mu$-index system such that $\lambda^\Gamma_\eta > \mu \Rightarrow \alpha
(\eta,\Gamma) = \lambda^\Gamma_\eta$. \nl
2) For $\Gamma = (S,\Gamma)$ a partial $\mu$-index system, for $\eta \in
S^\Gamma$ let $\Gamma^{[\eta]} = \langle S^{[\eta]},\bar \lambda
^{<\eta>} \rangle$ where

$$
S^{<\eta>} = \{\nu:\eta \char 94 \nu \in S\}
$$

$$
\bar \lambda^{<\eta>} = \langle \lambda_{\eta \char 94 <\nu>}:\nu \in 
S^{<\eta>} \rangle.
$$
\mn
We write also $\Gamma^{<\eta>} = (S^{\Gamma,<\eta>},
\bar \lambda^{\Gamma,<\eta>})$. \nl
3) $\eta^+ = \nu$ if $\ell g(\eta) = \ell +1,1 = \ell g(\nu),\eta
\restriction \ell = \nu \restriction \ell$ and $(\eta(\ell)+1 = \nu(\ell)$.
\enddefinition
\bn
\ub{Fact}:  If $\Gamma = (S,\bar \lambda)$ a partial $\mu$-index system we
say $\bar N$ is $\Gamma$-decomposition (in ${\Cal H}(\chi)$)
\mr
\item "{$(a)$}"  $\bar N = \langle N_\eta:\eta \in \Gamma \rangle$
\sn
\item "{$(b)$}"  $N_\eta \prec ({\Cal H}(\chi),\in,<^*_\chi)$
\sn
\item "{$(c)$}"  $N_\eta$ has cardinality $\lambda_\eta$
\endroster
\bn
\ub{Assignment}:  1) GCH tail forcing? \nl
2) \scite{3.10}? write a complete proof also larger $\mu$ so revise
\scite{t.2}.

\newpage
    
REFERENCES.  
\bibliographystyle{lit-plain}
\bibliography{lista,listb,listx,listf,liste}

\enddocument

\bye